%% file: dbva.tex
\input amstexl
\input lamstex
\catcode`\"=12
\font\sanss=cmss10
\font\ssanss=cmss8 scaled 900
\font\sssanss=cmss8 scaled 600

\font\teneurm=eurm10
\font\ssanss=cmss8 scaled 900

\newfam\ssfam
\textfont\ssfam=\sanss
\scriptfont\ssfam=\ssanss
\scriptscriptfont\ssfam=\sssanss
\def\ss#1{{\fam\ssfam\relax#1}}
\def\sss#1{\hbox{\ssanss #1}}

\mathchardef\za="710B  
\mathchardef\zb="710C  
\mathchardef\zg="710D  
\mathchardef\zd="710E  
\mathchardef\zve="710F 
\mathchardef\zz="7110  
\mathchardef\zh="7111  
\mathchardef\zvy="7112 
\mathchardef\zi="7113  
\mathchardef\zk="7114  
\mathchardef\zl="7115  
\mathchardef\zm="7116  
\mathchardef\zn="7117  
\mathchardef\zx="7118  
\mathchardef\zp="7119  
\mathchardef\zr="711A  
\mathchardef\zs="711B  
\mathchardef\zt="711C  
\mathchardef\zu="711D  
\mathchardef\zvf="711E 
\mathchardef\zq="711F  
\mathchardef\zc="7120  
\mathchardef\zw="7121  
\mathchardef\ze="7122  
\mathchardef\zy="7123  
\mathchardef\zvp="7124  
\mathchardef\zvr="7125 
\mathchardef\zvs="7126 
\mathchardef\zf="7127  
\mathchardef\zG="7000  
\mathchardef\zD="7001  
\mathchardef\zY="7002  
\mathchardef\zL="7003  
\mathchardef\zX="7004  
\mathchardef\zP="7005  
\mathchardef\zS="7006  
\mathchardef\zU="7007  
\mathchardef\zF="7008  
\mathchardef\zC="7009  
\mathchardef\zW="700A  
\catcode`\"=\active

\TagsOnRight
\loadmsam
\input paper.st\relax
\hsize=33pc
\hoffset=40pt
\vsize=54pc
\voffset=6pt
\def\*{{\textstyle *}}
\def\s*{{\scriptstyle *}}
\font\pmib=cmmib10
\def\proof{\demo{Proof}}
\def\endproof{\hfill \vrule height4pt width6pt depth2pt \enddemo}

\def\bK{{\bold K}}
\def\bE{{\bold E}}
\def\bF{{\bold F}}
\def\bC{{\bold C}}
\def\bL{{\bold L}}
\def\cK{\Cal L^\pm_{K^{\s*_r}}}
\def\cL{\Cal L^\pm_{L^{\s*_r}}}

\def\sT{\ss T}
\def\ssT{\sss T}
\def\m{{\slanted m}}
\def\xd{\text{\tenrm d}}
\def\id{\operatorname{id}}
\def\dt{\xd_{\sss T}}
\overfullrule=0pt
\def\ezt{{\teneurm\zt}}
\def\ezp{{\teneurm\zp}}
\def\ezk{{\teneurm\zk}}
\def\eza{{\teneurm\za}}
\def\ezb{{\teneurm\zb}}

        \newcounter\bibno
        \def\bb{\bib\key \bibno}

        \title Double vector bundles and duality
        \endtitle

    \author
                Katarzyna Konieczna \\
                Pawe\l\ Urba\'nski \\
       Division of Mathematical Methods in Physics,
               University of Warsaw \\
               Ho\D{z}a 74, 00-682 Warsaw, Poland\\
       {\sanss email:\ urbanski\@fuw.edu.pl}
    \endauthor
  
    \subjclass{ }

      \thanks{
    Supported by KBN, grant No 2 PO3A 074 10 }

   \abstract

The notions of the dual double vector bundle and the dual double
vector bundle morphism are defined. Theorems on canonical
isomorphisms are formulated and proved. Several examples are given.

   \endabstract
        \document
        \maketitle

        \hl1{Introduction}
	The notion of a double vector bundle was introduced by
Pradines in [{\bf 3}]. The most important examples  of double
vector bundles are  given by iterated tangent and cotangent
functors applied to a manifold $M$: $\sT \sT M$, $\sT^\* \sT M$,
$\sT \sT^\* M$, and $\sT^\* \sT^\* M$. The double vector bundle
structure on $\sT \sT M$ makes possible the Lagrangian
formulation of the dynamics in classical mechanics ([{\bf 5}]).
 	It appears that the framework of double vector bundles is very
convenient for many important constructions like linear 1-forms,
linear Poisson structures, special symplectic manifolds, linear
connections etc.

	The paper is organized as follows. In Sections 2 and 3 we
present in detail  the definitions of a double vector bundle
following Mackenzie [{\bf 2}], and morphisms of double vector
bundles. The tangent $\sT \bE$ to a vector bundle $\bE$ is an
example of a double vector bundle. In Section 4 we show that the
bundle, dual to a double vector bundle with respect to the
right (or left) vector bundle structure, is, in a natural way, a
double vector bundle. In particular, the  cotangent functor
applied to a vector bundle $\bE$ provides the double vector
bundle $\sT^\* E$.  Also the mapping, dual to  a morphism of
double vector bundles, is a morphism of double vector bundles
(Section~5). The main result of the paper is contained in
Section 6. We show that the third right dual to a double vector
bundle $\bK$ is canonically isomorphic to $\bK$. In the case of
$\bK = \sT^* \bE$, the third dual can be identified with $\sT^\*
\bE^\*$. The isomorphism  $\Cal R_K
\colon \bK \rightarrow  \bK^{\*_r\*_r\*_r}$ of Theorem~17 gives in this case
the canonical isomorphism of $\sT^\* \bE$ and $\sT^\* \bE^\*$.
The graph of this isomorphism is a Lagrangian submanifold of
$\sT^\*(E\times E^\*) \simeq \sT^\* E \times \sT^\* E^\*$
generated by the pairing between $\bE$ and $\bE^\*$.

	With the canonical identifications of $\bK$ with
$\bK^{\*_r\*_r\*_r}$ and of $\bK'$ with ${\bK'}^{\*_r\*_r\*_r}$,
the third dual to a double vector bundle morphism $\bold \zF
\colon \bK \rightarrow \bK'$ is equal $\bold \zF^{-1}$.

	In Section~7 we show that  the framework of double vector
bundles is suitable for concepts like linear vector fields and
forms, linear Poisson and symplectic structures, vertical and
complete lifts, and linear connections. We discuss in detail the
case of linear connections, compatible with a metric on $\bE$,
and symmetric connections on the tangent and cotangent bundles.      


	The main results of this paper were presented at the
conference in Vietri sul Mare, October 97, and one can find
their formulations in the proceedings of this conference ([{\bf 8}]).  

        \hl1{Double vector bundles}

        Let $\bK$ be a system $(\bK_r, \bK_l,
\bE, \bF)$ of vector bundles, where $\bK_r=(K,\zt_r, E)$,
$\bK_l=(K,\zt_l, F)$, $\bE =(E,\bar\zt_l,M)$, and $\bF= (F,\bar
\zt_r,M)$, such that the diagram

                $$ \Cgaps{.8}
        \CD  & K @()\L{\zt_l} @(-1,-1) @() \L{\zt_r}  @(1,-1)&   \\
       F@() \L{\bar\zt_r} @(1,-1) &    & E @() \L{\bar\zt_l} @(-1,-1) \\
                &  M &
                \endCD
							\tag\label{Fdb1}        $$
        is commutative.

        We introduce the following notation:
        \list
        \item ${\slanted m}_r$, ${\slanted m}_l$, $\bar
{\slanted m}_r$, and $\bar {\slanted m}_l$ will denote the
operation of addition in $\bK_r$, $\bK_l$, $\bE$, and $\bF$
respectively.
        \item we use also $v+_r w$ for ${\slanted m}_r(v,w)$,
$v+_lw $ for ${\slanted m}_l(v,w) $, and simply $+$ for all
other additions,
        \item $\bold 0_r$, $\bold 0_l$, $\bar\bold 0_r $, $\bar
\bold 0_l$ will denote the zero sections of $\zt_r$, $\zt_l$,
$\bar \zt_r$, and $\bar \zt_l$ respectively.
        \endlist

        Let us suppose that the pair $(\zt_r, \bar\zt_r)$ is a
vector bundle morphism $\bK_l \rightarrow \bE$. It follows that
$K\times _E K$ is a vector subbundle of $\bK_l\oplus \bK_l$ with
$(\zt_l\times \zt_l) (K\times _E K) = F\times _M F $. We denote
this subbundle by $\bK_l \oplus _E \bK_l$.  Moreover, the
addition ${\slanted m}_r \colon K \times _E K \rightarrow K$ is a fiber
bundle morphism which projects to $\bar {\slanted m}_r \colon F\times _M F
\rightarrow F$.

        \claim{Definition }{}
        A double vector bundle $\bK$ is a system $(\bK_r, \bK_l,
\bE, \bF)$ of vector bundles $\bK_r=(K,\zt_r, E)$,
$\bK_l=(K,\zt_l, F)$, $\bE =(E,\bar\zt_l,M)$ and $\bF= (F,\bar
\zt_r,M)$ such that the diagram \Ref{Fdb1} is commutative and
the following conditions are satisfied:

        \list
        \item  pairs $(\zt_l,\bar \zt_l)$, $(\zt_r,\bar \zt_r)$         are
vector bundle morphisms,
        \item pairs of additions $({\slanted m}_l, \bar
{\slanted m}_l)$  and $({\slanted m}_r, \bar {\slanted m}_r)$
are vector bundle morphisms  $\bK_r\times _F\bK_r\rightarrow
\bK_r$  and $\bK_l\times _E\bK_l\rightarrow \bK_l$ respectively,
        \item zero sections $\bold 0_r\colon \bE\rightarrow \bK_l$,
$\bold 0_l\colon \bF\rightarrow \bK_r$ are vector bundle
morphisms.
        \endlist
        \endclaim
        In the following we use the diagram \Ref{Fdb1} to
represent the double vector bundle $\bK$.
        \claim{Proposition}{} \label{Cdb1}
         The vector bundle structures of $\bK_r$ and $\bK_l$
coincide on the intersection $C$ of $\ker \zt_l$ and  $\ker
\zt_r$.
        \endclaim
        \proof
         Since $\zt_r$ and $\zt_l$ are linear in fibers it follows
that for each $m\in M$, we have
                $$ \bold 0_r\circ \bar\bold  0_l (m) \in C, \ \ \bold
0_l\circ \bar \bold 0_r (m) \in C .
                                                                                $$
        Let us denote $0_r(m) = \bold 0_r\circ \bar \bold 0_l (m)$,
$0_l(m)= \bold 0_l\circ \bar \bold  0_r (m) $.
        Since ${\slanted m}_l\colon \bK_l \times _F \bK_l \rightarrow \bK_l$
is a vector bundle morphism, we have
                $$ {\slanted m}_l(0_r(m),0_r(m)) = 0_r(m).
                                                                                $$
         It follows that $0_r(m)$ is a neutral element of ${\slanted m}_l$,
i.e,  $0_l = 0_r$.
        For $k,k'\in C$, such that $\bar \zt_l(\zt_r(k))
= \bar \zt_l(\zt_r(k')) =m$, we have
                $$ (k,k') = (0_r(m),k')+_r(k,0_r(m)) =
(0_l(m),k')+_r(k,0_l(m))
                                                                      $$
and, consequently,
                $$ {\slanted m}_l(k,k') =  {\slanted
m}_l(k,0_l(m))+_r {\slanted m}_l(0_l(m),k') = k+_rk'.
                                                                     $$
        \endproof
        Thus, we have a vector bundle $\bC =(C, \zt, M) $, where
$\zt = \bar \zt_l\circ \zt_r = \bar \zt_r\circ \zt_l$. This vector bundle is
called {\it the core } of $\bK$.

        \claim{Proposition }{} \label{Cdb2}
        \list
        \item $\ker \zt_r$ with the vector bundle structure induced
from $\bK_r$ is canonically isomorphic to the  Whitney sum $
\bF\oplus _M \bC$.
        \item $\ker \zt_r$ with the vector bundle structure induced
from $\bK_l$ is canonically isomorphic to the  vector bundle $
F\times _M\bC$, i.e., to the pull-back of $\bC$ by
the projection $\bar \zt_r$.
        \item $\ker \zt_l$ with the vector bundle structure induced
from $\bK_l$ is canonically isomorphic to the Whitney sum
$\bE\oplus _M \bC$.
        \item $\ker \zt_l$ with the vector bundle structure induced
from $\bK_r$ is canonically isomorphic to the  vector bundle $
E\times _M \bC$, i.e.,  to the pull-back of $\bC$ by the
projection $\bar \zt_l$.
        \endlist
        \endclaim
        \proof
        \list
        \item
        Since the zero section $\bold 0_l$ is a vector bundle morphism
                $$ \bold 0_l \colon \bF\rightarrow \bK_r,
                                                                                $$
        its image $\bold 0_l (F)$ is a vector  subbundle of $
\bK_r$, contained in  $\ker \zt_r$ and  isomorphic to $\bF$.

Let $v\in \ker \zt_r$ and let $m = \bar \zt_l\circ \zt_r(v) = \bar \zt_r \circ
\zt_l(v)$. We have $\zt_r(v) = \bar \bold 0_l$ and, since $v\in \ker \zt_r$,
                $$\zt_r(\bold 0_l (\zt_l(v))) = \bar \bold 0_l(\bar
\zt_r \circ \zt_l(v)) = \bar \bold 0_l(m).
                                                                                                                                $$
        It follows, that the pair $(v, \bold 0_l(\zt_l(v)))$ is in
the domain of $\m_r$. Now, we can define two mappings
                $$ \P_C, \, \P_F \colon \ker \zt_r \rightarrow \ker \zt_r
                                                                                                                                $$
         by the formulae
                $$ \split \P_F(v) &= \bold 0_l(\zt_l(v)), \ \
m=(\bar \zt_r\circ \zt_l) (v), \\
        \P_C(v)  &= v -_r \bold 0_l(\zt_l(v)).
                                                        \endsplit                                               \tag \label{Fdb2}$$

        It is evident that $\P_F^2 = \P_F$, $\P_F \P_C =0$, $\P_C
+\P_F = id$ and
                $$\split \P_C^2(v) & = (v-_r \bold 0_l(\zt_l(v))) -_r
\bold 0_l(\zt_l(v - _r \bold 0_l(\zt_l(v)))) \\
                & = (v -_r \bold 0_l(\zt_l (v))) -_r 0_r(m) = v -_r
\bold 0_l(\zt_l (v)) \\
                & = \P_C(v).
                                                        \endsplit                                                       $$

        It follows that $\P_F, \P_C$ are projectors, which define a
splitting of $\ker\zt_r$. We have that $\P_C(v) \in \ker \zt_r$ and
                $$ \zt_l(\P_C(v)) = \zt_l(v) - \zt_l(\bold 0_l(\zt_l(v)))
= \bar \bold 0_l(m),
                                                                                                                                $$
and, consequently, $\P_C(v) \in \bC$.

It is obvious that, for each $m\in M$, the intersection of fibers
over $m$ of  $\bold 0_l (F)$ and  $\bC$ is trival and equal to $\{
0_r(m)\}$. Since the image of $\P_F$ is canonically isomorphic
to $\bF$, we conclude, finally, that $\ker \zt_r$ is canonically
isomorphic to $\bF \oplus_M \bC $.

        \item From (1) we have that   $\ker \zt_r$ can be
identified (as a manifold) with $F\times _M C$. With this
identification $\zt_l$ is the canonical projection on $F$ and the
zero section $\bold 0_l$ is given by
                $$ \bold 0_l \colon F\rightarrow F\times _M C \colon
f\mapsto (f,0).
                                                                  $$
        The addition ${\slanted m}_l$ defines a vector bundle morphism
                $${\slanted m}_l\colon \bK_r\times _F\bK_r \supset \ker
\zt_r\times _F\ker \zt_r = \bF\oplus\bC\oplus\bC\rightarrow
\bF\oplus \bC\subset \bK_r,
                                                                  $$
        hence
                $$\multline
        (f,k)+_l (f,k') = {\slanted m}_l((f,k,k')) = {\slanted m}_l((f,0,0)+_r(0,k,k'))=\\
        = (f,0)+_r(0,k+k') = (f,k+k').
                                \endmultline                     $$
        \item "(3,4)" The proof is analogous.
        \endlist
        \endproof

        The following two propositions will be useful in the
next section.
                \claim{Proposition}{} \label{Cdb3}
        Let $v,w\in K$ and let
                $$\split \zt_r(v) = \zt_r(w) = e\in E,\\
                        \zt_l(v) = \zt_l(w) = f\in F.
                        \endsplit                                  $$
        There exists $k\in C$ such that $\zt(k) = \bar \zt_l(e) =
\bar \zt_r(f)$ and, using the identifications of Proposition~\Ref{Cdb2},
                $$\split  v = w +_r (e,k),\\
                                v = w +_l (f,k).
                                        \endsplit                $$
        \endclaim
        \proof
        We have $\zt_l(v-_r w) = f-f = 0$ and $\zt_r(v-_l w) = e- e
= 0$. It follows from Proposition~\Ref{Cdb2} that there exist
$k,k'\in C$ such that
                $$ v-_rw = (e,k),\ \ v-_lw = (f,k').
                                                                $$
          We show that $k=k'$. Indeed, since $(v,w)\in K\times _EK$ and
                $$(v,w) = (v-_lw,w-_lv) +_l (w,v),
														$$
we have (linearity of ${\slanted m}_r$)
                $$ v-_rw =((v-_lw)-_r (w-_lv)) +_l(w-_rv).
                                                     $$
Hence (Proposition~\Ref{Cdb2})
                $$(e,k)= ((f,k')-_r(f,-k')) +_l(e,k) = (0,2k') +_l(e,-k)
                                                    $$
and
                $$ (0,2k') = (e,k)-_l(e,-k)= (0,2k).
                                                   $$
        \endproof
        \claim{Proposition}{} \label{Cdb4}
        Let $v,v',w,w'\in K$ be such that
                $$\split \zt_r(v) =& \zt_r(w) = e\in E,\\
                        \zt_r(v') =& \zt_r(w') = e'\in E,\\
                \zt_l(v) = \zt_l(w) = \zt_l(v')& = \zt_l(w')= f\in F,\\
                                v+_lv'& =w+_lw'.
                        \endsplit
													$$
        Then
                $$\split v-_rw&= (e,k)\in \ker \zt_l\\
                                 v'-_rw'&= (e',-k)\in \ker \zt_l
                        \endsplit
												$$
        \endclaim
        \proof
        Since $v-_rw\in \ker \zt_l$ and $v'-_rw'\in \ker \zt_l$,
there exist $k,k'\in C$ such that
                $$ v-_rw = (e,k),\ \ \ v'-_rw'  = (e',k').
                                                   $$
        From Proposition~\Ref{Cdb3} we have also
                $$v-_lw =(f,k),\ \ \ v'-_lw' = (f,k').
                                                   $$
        The equality $v+_lv'=w+_lw'$ implies
                $$ (f,k)= v-_lw = w'-_lv' = (f,-k')
                                                   $$
         and $k' =-k$.
        \endproof
        \hl2""{Local coordinates}
        Let $(x^i)_{i=1}^n$ be a coordinate system on $M$. In
the bundles $\bE, \bF $, we introduce   coordinate systems
$((x^i)_{i=1}^n, (e^a)_{a=1}^{n_E})$ and $((x^i)_{i=1}^n,
(f^A)_{A=1}^{n_F})$.  We denote also by $x^i, e^a, f^A$ their
pull-backs to the coordinates on $K$. It follows from
Proposition~\Ref{Cdb3} that  we can introduce coordinates
$(c^\za)^{n_C}_{\za =1}$  such that $(x^i, e^a, f^A, c^\za)$ is
a local coordinate system on $K$ and
        \list
        \item"" $c^\za(v+_rw) = c^\za(v) + c^\za(w)$,
        \item"" $c^\za(v+_lw) = c^\za(v) + c^\za(w)$,
        \item"" $ c^\za \circ \bold 0_r = 0$,
        \item"" $ c^\za \circ \bold 0_l = 0$,
        \endlist
         It follows that $(x^i, c^\za)$ is a vector bundle
coordinate system in $\bC$. The operation of addition $+_r$ is
characterized by the following equalities

                $$\split
        c^\za(v+_r w) &= c^\za(v) + c^\za(w), \\
        f^A(v+_r w)   &= f^A(v) + f^A(w), \\
        e^a(v+_r w)   &= e^a(v) = e^a(w), \\
        x^i(v+_r w)   &= x^i(v) = x^i(w) .
                                \endsplit
												$$

The operation of addition $+_l$ is
characterized by
                $$\split
        c^\za(v+_l w) &= c^\za(v) + c^\za(w), \\
        f^A(v+_l w)   &= f^A(v) = f^A(w), \\
        e^a(v+_l w)   &= e^a(v) + e^a(w), \\
        x^i(v+_l w)   &= x^i(v) = x^i(w) .
                                \endsplit
											$$

        \hl2""{Examples}

{ 1.} Let $\bE =(E,\bar\zt_l,M)$, $\bF= (F,\bar\zt_r,M)$, $\bC
=(C,\zt,M)$ be vector bundles and let $K=F\times _M C\times _M
E$. By $\bK(\bF,\bC,\bE)$ we denote a double vector bundle
represented by the diagram
                $$ \Cgaps{.8}
        \CD  & K @()\L{\zt_l} @(-1,-1) @() \L{\zt_r}  @(1,-1)&   \\
	F@() \L{\bar\zt_r} @(1,-1) &    & E @() \L{\bar\zt_l} @(-1,-1) \\
                &  M &
                \endCD \ \ \ ,
									\tag\label{Fdb3}        $$
        where $\zt_r \colon K=F\times _M C\times _M  E \rightarrow
E$ and $\zt_l \colon  K=F\times _M C \times _M E \rightarrow F$
are the canonical projections. The right and left vector bundle
structures are obvious:
                $$ \split
                (f,k,e) +_r (f',k',e) &= (f+f', k+k',e) \\
                (f,k,e) +_l (f,k',e') &= (f, k+k',e+ e') .
                \endsplit                                                                               $$
        The core of $\bK(\bF,\bC,\bE)$ can be identified with $\bC$.

        { 2.}
        Let $\bE = (E,\zt,M)$ be a vector bundle. The tangent
manifold $\sT E$ has two vector bundle structures ([{\bf 5}]):
        \list
        \item"" the canonical vector bundle structure of the tangent
bundle, on the canonical fibration $\ezt_E  \colon \sT E
\rightarrow E $,
        \item"" the tangent vector bundle structure on the
tangent fibration  $\sT\zt \colon  \sT E  \rightarrow
\sT M$. 
        \endlist
        It is easy task to verify that the diagram
                $$ \Cgaps{.8}
        \CD  & K @()\L{\sT\zt} @(-1,-1) @() \L{\ezt_E}  @(1,-1)&   \\
                \sT M@() \L{\ezt_M} @(1,-1) &    & E @() \L{\zt} @(-1,-1) \\
                &  M &
                \endCD
										\tag\label{Fdb4}        $$
represents a double vector bundle. We denote this double  vector
bundle by $\sT\bE$. The core consists of vertical tangent vectors
at the zero section of $\bE$. Thus, it can be identified, in an
obvious way, with $\bE$.

        { 3.}
        Let  $\bK$ be a double vector bundle represented by the
diagram
                $$ \Cgaps{.8}
        \CD  & K @()\L{\zt_l} @(-1,-1) @() \L{\zt_r}  @(1,-1)&   \\
        F@() \L{\bar\zt_r} @(1,-1) &    & E @() \L{\bar\zt_l} @(-1,-1) \\
                &  M &
                \endCD
										\tag\label{Fdb5}        $$
then also the diagram
                $$ \Cgaps{.8}
        \CD  & K @()\L{\zt_r} @(-1,-1) @() \L{\zt_l}  @(1,-1)&   \\
                E@() \L{\bar\zt_l} @(1,-1) &    & F @() \L{\bar\zt_r} @(-1,-1) \\
                &  M &                \endCD
									\tag\label{Fdb6}        $$
represents a double vector bundle. We denote it by $\ss J(\bK)$.

        In  examples  (1) and (2), we identified canonically the
core  with a  certain vector bundle.  We shall write in the
following a diagram
                $$ \Cgaps{.8}
        \CD  & K @()\L{\zt_l} @(-1,-1) @() \L{\zt_r}  @(1,-1)&   \\
                F@() \L{\bar\zt_r} @(1,-1) & C @() \0( @(0,1) @()
\L{\zt} @(0,-1)    & E @() \L{\bar\zt_l} @(-1,-1) \\
                &  M &
        \endCD
                                                                \tag \label{Fdb7}$$
        instead of the diagram \Ref{Fdb1},  if we identify the core of
\Ref{Fdb1} with the vector bundle $\bC$. In the case of $\bK =\sT
\bE$ we write then the diagram

                $$ \Cgaps{.8}
        \CD  & \sT E @()\L{\sT\zt} @(-1,-1) @() \L{\ezt_E}  @(1,-1)&   \\
                \sT M@() \L{\ezt_M} @(1,-1) & E @() \0( @(0,1) @()
\L{\zt} @(0,-1)    & E @() \L{\zt} @(-1,-1) \\
                &  M &
        \endCD
                                                               \tag \label{Fdb8}$$

        \hl1{Morphisms of double vector bundles \label{s1}}
        Let $\bK =(\bK_r, \bK_l,\bE, \bF)$ and $\bK' =(\bK'_r,
\bK'_l,\bE', \bF')$ be double vector bundles with cores $\bC$ and
$\bC'$ respectively.  A morphism $\bold \zF\colon \bK
\rightarrow \bK'$ of double vector bundles is a family $\bold
\zF =(\zF, \zF_r, \zF_l, \bar \zF)$ of mappings
                $$ \zF \colon K \rightarrow K', \qquad  \zF_r \colon E
\rightarrow E',\qquad
        \zF_l \colon F \rightarrow F', \qquad \text{and} \quad \bar  \zF \colon M \rightarrow M'
                                                                     $$
         such that $\bold \zF_r = (\zF,\zF_r), \ \bold \zF_l  = (\zF,\zF_l)$,
$\bar \bold \zF_r= (\zF_r, \bar \zF)$, and $\bar \bold \zF_l =(\zF_l, \bar
\zF)$ are morphisms of vector bundles
                $$ \bold \zF_r\colon \bK_r\rightarrow \bK'_r, \qquad  \bold
\zF_l\colon \bK_l\rightarrow \bK'_l, \qquad
        \bar \bold \zF_r\colon \bE \rightarrow \bE', \qquad
\text{and} \quad \bar \bold \zF_l\colon \bF\rightarrow \bF'.
                                                               $$

        We have thus a commutative diagram

                $$ \cgaps{.8;.8;1.2;.8;.8} \rgaps{.8;.4;.8}
        \CD  & K @()\L{\zt_l} @(-1,-2) @() \L{\zt_r} @(1,-1) @() \L{\zF}
 @(3,0) & & & K' @() \L{\zt_l'} @(-1,-2) @() \L{\zt_r'} @(1,-1) & \\
        & & E @()\dl{-2}\l{\bar \zt_l} @(-1,-2) @()\dL{-7} \L{\zF_r} @(3,0) & & &
E' @()\L{\bar \zt_l'} @(-1,-2) \\
        F  @()\L{\bar \zt_r} @(1,-1) @()\dL{7}\L{\zF_l} @(3,0) & & & F'
@()\L{\bar \zt_r'} @(1,-1)  & &  \\
        & M @()\L{\bar \zF} @(3,0) & & & M' &
                                           \endCD\ \ .  \tag \label{Fdb9}$$

        \claim{Proposition}{}
        Let  $\bold \zF\colon \bK \rightarrow \bK'$ be a morphism of
double vector bundles. Then
        \list
        \item $\zF(\ker \zt_r)\subset \ker \zt_r'$,
        \item $\zF(\ker \zt_l)\subset \ker \zt_l'$,
        \item $\zF (C) \subset C'$.
        \endlist
        \endclaim
        \proof
        Since $\bar \bold \zF_r \colon  \bE\rightarrow \bE'$ is a morphism of
vector bundles, it maps the zero section of $\bE$ into the zero section
of $\bE'$. It follows that       $\zF(\ker \zt_r)\subset \ker \zt_r'$.
Similarly, $\zF(\ker \zt_l)\subset \ker \zt_l'$ and,
consequently,   $\zF (C) \subset C'$.
        \endproof
        We denote by $\bold \zF_c = (\zF_c, \bar \zF)$ the
morphism of vector bundles $\bC$ and $\bC'$  induced by $\zF$.

        Let  $(x^i, e^a, f^A, c^\za)$  be an adapted local
coordinate system on $K$ and let $(\bar x^{\bar i}, \bar e^{\bar
a}, \bar f^{\bar A}, \bar c^{\bar \za})$ be an adapted
coordinate system in $K'$.  We have
                $$\split  \bar x^{\bar i}\circ \zF & =  \zF^{\bar i},\\
                    \bar e^{\bar a} \circ \zF  & =  \zF^{\bar a}_be^b,\\
                   \bar f^{\bar A} \circ \zF  & =  \zF^{\bar A}_Bf^B,\\
        \bar c^{\bar \za} \circ \zF  & = \zF^{\bar \za}_\zb
c^\zb + \zF^{\bar \za}_{aA}e^af^A,
                                \endsplit
                                                 \tag\label{Fdb10} $$
         where  $ \zF^{\bar i},\, \zF^{\bar a}_b$, $ \zF^{\bar A}_B,
\, \zF^{\bar \za}_\zb, \, \zF^{\bar \za}_{aA}$ are functions on
the domain of $(x^i)$ in $M$.

        \hl2""{Examples}

        { 1.} Let $\bE =(E,\zt, M)$ and $\bE' =(E',\zt',M')$ be
vector bundles and let $\bold \zF = (\zF, \bar \zF)$ be a vector
bundle morphism
                $$ \zF \colon  \bE \rightarrow \bE'.
                                                                                                                                $$

        The quadruple $\sT\bold \zF =  (\sT \zF, \zF, \sT \bar \zF,
\bar \zF)$ defines a morphism of double vector bundles
                $$ \sT \bold \zF \colon  \sT\bE \rightarrow \sT \bE'
                                                                                                                                $$
and, with the identification of cores as in the diagram
\Ref{Fdb8}, we have $\bold \zF_c \colon \bE \rightarrow \bE' $,
$\zF_c = \zF$.

        { 2.} An essential role in the Lagrangian formulation of the
classical mechanical system is played by the isomorphism
$\ezk_M$, which relates $\sT \sT M$ with $\ss J(\sT \sT M)$
([{\bf 5}]). Here, $\sT M$ is the vector bundle of tangent
vectors. All three vector bundle morphisms $\bar \bold \zF_r,
\bar \bold \zF_l, \bold \zF_c \colon  \sT M \rightarrow \sT M$
are indentities.

        { 3.} Let  $\bK = \bK(\bF,\bC,\bE)$ and  $\bK'= \bK(
\bF', \bC',\bE')$. If $\bold
\zF =(\zF, \zF_r, \zF_l, \bar \zF)$ is a morphism of double vector bundles,
                        $$\bold \zF\colon \bK \rightarrow \bK',
                                         \tag\label{Fdb11}       $$
then
                $$ \zF (f,c,e) = ( \zF_l(f), \zF_c(c) + \zC(f,e),\zF_r(e)),
                                                             $$
where the mapping
                $$ \zC \colon \bF\times _M \bE \rightarrow \bC'
                                                                                                                                $$
is bilinear.

        \hl1{The right dual}
        Let $\bK^\*_r$be the vector bundle, dual to $\bK_r$.  We
denote by $K^{\*_r}$ the total fiber bundle space  and by
$\zp_l$ the projection
                $$\zp_l\colon K^{*_r}\rightarrow E.
                                                             $$

        Let  $a\in K^{\*_r}$ and $k \in C$ satisfy $\zt (k) =
\bar \zt_l(\zp_l(a))$. We can evaluate $a$ on a vector
$(\zp_l(a),k)$  of $\ker \zt_l$. We define a mapping
$\zp_r\colon K^{\*_r}\rightarrow C^\*$  by the formula
                $$\langle k,\zp_r(a)\rangle = \langle (\zp_l(a),k),
a\rangle .
                                       \tag \label{Fdb12}    $$

It follows directly from this construction that
        \claim{Proposition}{}
        The mapping  $\zp_r\colon K^{\*_r} \rightarrow
C^\*$  is a morphism of vector bundles
                $$\zp_r \colon \bK_r^\* \rightarrow \bC^\*.
                                                             $$
        \endclaim

         We define a relation
        $$ P_r({\slanted m}_l)\colon K^{\*_r}\times
K^{\*_r}\rightarrow K^{\*_r}
                                        $$
        in the following way: $c\in P_r({\slanted m}_l)(a,b)$ if
        \list
        \item $\zp_l(c) = \zp_l(a) + \zp_l(b)$,
        \item $\langle w,c\rangle  = \langle v,a\rangle + \langle
v',b\rangle$  for each $w,v,v'\in K$ such that $\zt_r(w) =
\zp_l(c)$, $\zt_r(v) = \zp_l(a)$, $\zt_r(v') = \zp_l(b)$, and $w
={\slanted m}_l(v,v')$.
        \endlist

        \claim{Proposition}{}
        $P({\slanted m}_l)$ is univalent and, if $c = P({\slanted m}_l)(a,b)$, then
$\zp_r(a)= \zp_r(b)= \zp_r(c)$.
        \endclaim
        \proof
        $P({\slanted m}_l)$ is univalent, because for each $w$ we can find $v,v'$
such that ${\slanted m}_l(v,v') = w$ and, consequently, $c$ is
completely determined by the condition
                $$ \langle w,c\rangle = \langle v,c \rangle  + \langle v',c\rangle .
                                                                                                                                $$

        Let  $k\in C$ and $c = P({\slanted m}_l)(a,b)$. From the definition
of $\zp_r$ we have, using identifications of Proposition~\Ref{Cdb2},
                $$\multline \langle k, \zp_r(a)\rangle  - \langle k,
\zp_r(b)\rangle = \langle (\zp_l(a),k), a \rangle  +
\langle(\zp_l(b), -k), b\rangle= \\
        =\langle ((\zp_l(a),k)+_l(\zp_l(b),-k)),c\rangle   =\langle
(\zp_l(a)+ \zp_l(b),0), c\rangle =0.
                        \endmultline                        $$
         It follows that $\zp_r(a) = \zp_r(b)$. Moreover,
                $$ \multline
        \langle k, \zp_r(c)\rangle = \langle(\zp_l(c),
k),c\rangle =\langle ((\zp_l(a),k)+_l(\zp_l(b),0)),c\rangle  =\\
        =  \langle (\zp_l(a),k), a\rangle  + \langle
(\zp_l(b),0),b\rangle =  \langle (\zp_l(a),k), a\rangle =
\langle k, \zp_r(a)\rangle .
                                \endmultline                  $$
        Hence, $\zp_r(a) = \zp_r(c)$.

        \endproof

        \claim{Proposition}{}
        A pair $(a,b) \in K^{\*_r}\times  K^{\*_r}$ is in the
domain of $P({\slanted m}_l)$ if and only if $\zp_l(a) = \zp_l(b)$.
        \endclaim
        \proof
         It is enough to show that, if
$\zp_l(a) =\zp_l(b) $, then there exists $c\in K^{\*_r}$  such
that $c = P({\slanted m}_l)(a,b)$.

        Let  $\zp_r(a) =\zp_r(b) $. We define $c \in
\zp_l^{-1}(\zp_l(a) + \zp_l(b)) $ by
                $$ \langle w,c\rangle  = \langle v,a\rangle + \langle
v',b\rangle
                                           \tag\label{Fdb13}       $$
        where $w, v,v'\in K$ are such that  $\zt_r(v) = \zp_l(a)
=e$, $\zt_r(v') = \zp_l(b)=e'$ and $w ={\slanted m}_l(v,v')$. Let
$w={\slanted m}_l(u,u')$ be another representation of $w$, such that
$\zt_r(u) = e, \zt_r(u') = e'$.
        We have also $\zt_l(u) = \zt_l(u') = \zt_l(w)$. It follows
from Proposition~\Ref{Cdb4} that
                $$\split u-_rv & =(e,k)\in \ker \zt_l \\
                                u'-_rv' & =(e',-k)\in \ker \zt_l
                        \endsplit                                 $$
and, consequently,
                $$ \multline
        \langle u,a\rangle + \langle u',b\rangle  =  \langle
u-_rv,a\rangle +\langle u'-_rv',b\rangle  +  \langle v,a\rangle
+ \langle v',b\rangle
        = \langle (e,k),a\rangle + \langle (e',-k),b\rangle + \\
+ \langle u,a\rangle + \langle u',b\rangle = \langle
k,\zp_r(a)\rangle -\langle k,\zp_r(b)\rangle  +\langle
u,a\rangle + \langle u',b\rangle = \langle u,a\rangle + \langle
u',b\rangle
                        \endmultline                              $$
        This proves that $c$ is a well defined function on
$(\zt_r)^{-1}(e+e')$. Now, we have to show that this function is
linear. Let be $w = w_1 +_r w_2$ and $w_1 = u+_lu'$, $w_2 =
v+_lv'$. Since $+_r$ is linear on  $\bK_l\times _E\bK_l$, we have
                $$(u+_lu')+_r(v+_lv') = (u+_rv) +_l (u'+_r v')
                                                                 $$
and
                $$ \multline
\langle w_1+_rw_2,c\rangle = \langle u+_rv,a\rangle
+ \langle u'+_rv',b\rangle = \\
        = \langle u,a\rangle +\langle v,b\rangle  + \langle
u',b\rangle + \langle v',b\rangle = \langle w_1,c\rangle +
\langle w_2,c\rangle .
                                \endmultline                       $$
        The function $c$ is additive and, consequently, linear.
        \endproof
        We have then the operation  $P({\slanted m}_l)$ in
fibers of $\zp_r$. In the following, we shall use also $+_r$
instead of $P({\slanted m}_l)$ to denote this opeation. The
multiplication $\cdot _r$ is defined by the formula
                $$ \langle r\cdot _lv, r\cdot _ra\rangle  =r\langle v,a\rangle.
                                                                   $$
        We show in the following proposition that with these
operations $K^{\*_r}$ becomes  a vector bundle over
$C^{\textstyle *}$.
        \claim{Proposition}{}
On each fiber of $\zp_r$ the operation $P({\slanted m}_l)$ defines the
structure of an Abelian group.
        \endclaim
        \proof
 {\bf Associativity.} Let  $a,b,c \in K^{\*_r}$ be such that
$\zp_r(a) = \zp_r(b) = \zp_r(c)$ and let  $u,v,w \in K$ be such
that $\zt_l(u) = \zt_l(v) = \zt_l(w)$ and $\zt_r(u) = \zp_l(a)$,
$\zt_r(v) = \zp_l(b)$, $\zt_r(w) = \zp_l(c)$. We have
                $$\langle w+_lv+_lu,c+_r(b+_ra) \rangle = \langle w,c\rangle
+\langle v+_lu,b+_ra\rangle = \langle w,c\rangle +\langle
v,b\rangle +\langle u,a\rangle .
                                                                     $$

{\bf Commutativity.} Obvious.

{\bf Neutral element.} Let us choose an element $\za \in
\bC^{\textstyle *} $. Since $\ker \zt_r$, with the vector bundle structure of
$\bK_r$, is  identified with $ \bF\oplus_M \bC$, we can define an
element $e_\za \in K^{\*_r}$
by
                $$\zp_l(e_\za) = \bar \bold 0_l(\zt(\za)),\ \ \
\langle (f,k), e_\za \rangle = \langle k,\za\rangle .
                                                                     $$
          For each $a\in K^{\*_r}$ such that $\zp_r(a) = \za$ and
for each $w\in K$ such that $\zt_r(w) = \zp_l(a) $, we have
                $$ \langle w, a+_re_\za \rangle = \langle
w+_l(f,0),a+_re_\za \rangle = \langle w,a\rangle + \langle
(f,0),e_\za \rangle  = \langle w,a\rangle ,
                                                                      $$
        where $f = \zt_l(w)$.

        {\bf Inverse element.} Is equal $(-1)\cdot _ra$.
        \endproof

\medskip
        \hl2""{Local coordinates}
        Let $(x^i,e^a, f^A, c^\za)$ be an adapted coordinate system
on $K$ and let $(x^i, e^a, p_A, q_\za)$ be the adopted coordinate
system on the dual bundle $K^{\*_r}$, i. e., the canonical
evaluation is given by the formula
                $$ \langle v, a\rangle = \sum_A p_A(a)f^A(v) + \sum_\za
q_\za(a) c^\za(v).
                                                       \tag \label{Fdb14}$$
        We use $(x^i,c^\za)$ as a coordinate system on $C$ and
$(x^i, q_\za)$ as a coordinate system on $C^\*$. In these
coordinate systems, we have
                $$\split x_i\circ \zp_r(a) &= x_i(a), \\
                                q_\za\circ \zp_r(a) &= q_\za(a),
                                \endsplit             \tag \label{Fdb15}$$
and
                $$\split x^i(a +_r b) &= x^i (a) = x^i(b), \\
                                e^a(a +_r b) &= e^a (a) + e^a(b), \\
                                p_A(a +_r b) &= p_A (a) + p_A(b), \\
                                q_\za(a +_r b) &= q_\za (a) =  q_\za(b) .
                                \endsplit             \tag \label{Fdb16}$$

        It follows that $(K^{\*_r},\zp_r, C^\*)$ is a vector bundle.
We denote it by $\bK^{\*_r}_r$. The vector bundle $(K^{\*_r},
\zp_l,E)$ we denote by $\bK^{\*_r}_l$.

        \claim{Theorem}{} The system
$\bK^{\*_r} = (\bK^{\*_r}_r,\bK^{\*_r}_l,\bC^{\textstyle *} ,\bE)$ is a double
vector bundle.
                \label{Cdb5} \endclaim
        \proof
        We have to show that the projection, zero section and
the operation of addition related to one vector bundle structure is
linear with respect to another one.

{\bf Projections.} Linearity of $\zp_r\colon
\bK^{\*_r}_l \rightarrow \bC^{\textstyle *} $ follows directly
from its definition. Linearity of the projection $\zp_l\colon
\bK^{\*_r}_r \rightarrow \bE $ is contained in the definition of
$P({\slanted m}_l)$.

{\bf Zero sections.} Linearity of ${\bold 0}_l\colon
\bE\rightarrow \bK^{\*_r}_r$ follows directly from the definition
of $P({\slanted m}_l)$. Now, let
                $$ \bold 0_r\colon \bC^{\textstyle *} \rightarrow
\bK^{\*_r}_l = (\bK_r)^{\textstyle *}
                                                                    $$
        be the zero section. For $\za\in  C^{\textstyle *} $,  \
$e_\za = \bold 0_r(\za) $ is an element of the space dual to
$\zt_r^{-1}(0)$, defined by
                $$\zp_l(e_\za) = \bar \bold 0_l(\zt(\za)),\ \ \
\langle (f,k), e_\za \rangle = \langle k,\za\rangle .
                                                 $$
Hence linearity.

{\bf Additions.}
        First, we show that $+_l$ is a morphism of vector bundles
                $$ +_l\colon \bK^{\*_r}_r\times _E\bK^{\*_r}_r
\rightarrow \bK^{\*_r}_r.
                                                \tag\label{Fdb17}   $$
Let $a,b,a',b'\in K^{\*_r}$ be such that
                $$ \zp_l(a)= \zp_l(b) = e, \ \zp_l(a') =
\zp_l(b')=e',\ \zp_r(a) = \zp_r(a'),\ \ \text{and} \ \zp_r(b)=\zp_r(b').
                                                        $$
        Since $\zp_l \colon \bK^{\*_r}_r\rightarrow \bE$ is linear,
we have
                $$\zp_l(a +_ra') = \zp_l(b+_lb') =e +e'
                                                                    $$
and, consequently, $(a+_ra', b+_rb')\in \bK^{\*_r}_r\times
_E\bK^{\*_r}_r$. The operation
        $+_l$ is additive with respect to $+_r$ and $+_r$ is
additive with respect to $+_l$ if for each  quadruple
$(a,b,a',b')$
                $$ (a+_ra')+_l(b+_rb') = (a+_lb)+_r(a'+_lb').
                                                                                $$
         To show this equality let us take $w,v,v'\in K$ such that
$\zt_r(v)=e$, $\zt_r(v) = e'$, $\zt_l(v) = \zt_l(v') $ and $w= v
+_lv'$. We have
                $$\multline \langle w, (a+_ra')+_l(b+_rb')\rangle =
\langle w,a+_ra'\rangle  + \langle w,b+_rb'\rangle =\\
        = \langle v,a\rangle + \langle v',a'\rangle +\langle v,b\rangle +
\langle v',b'\rangle = \langle v,a+_lb \rangle + \langle v',a'+_lb'\rangle 
\rangle = \langle w,(a+_lb)+_r(a'+_lb').
                                \endmultline                            $$

        \endproof
        We identify the kernel $\ker \zt_r$  with $\bF \oplus_M \bold C$ and,
consequently, the kernel of $\zp_l$ with $\bC^{\textstyle *} \oplus_M
\bF^{\textstyle *} $. With this identifications, we have that
                $$ \langle (f,c),(c^{\textstyle *} ,f^{\textstyle *} )\rangle =
\langle f,f^{\textstyle *} \rangle + \langle c,c^{\textstyle *} \rangle
                                        \tag\label{Fdb18}               $$
and that $\zp_r$ is the canonical projection $\bC^\*\oplus\bF^\* \rightarrow
\bC^\*$. It follows that the core of $\bK^{\*_r}$ can be identified with
$\bF^\*$ and that
                $$ \ker \zp_r = \bE\oplus \bF^\*.
                                                                                $$
                \claim{Proposition}{}
        Let  $a=(e,f^\*)\in \ker \zp_r$ and $b \in K^{\*_r}$. Then, for
$v,w\in K$  such that $\zt_r(v) = e $ and $w = (e,c)\in \ker \zt_l,
e=\zp_l(b)$, we have
                $$\split \langle v,a\rangle =& \langle \zt_l (v),f^\* \rangle \\
                                \langle w,b\rangle =& \langle c, \zp_r(b) \rangle
                        \endsplit                                                       $$
                \label{Cdb6} \endclaim                                                          
        \proof
        Let $f = \zt_lv$. We have $(f,0)\in F\times_M C = \ker \zt_r$,
 $v=(v -_l(f,0))+_l(f,0)$, and $a=(e,f^\*) = (e,0)+_r(0,f^\*)$. But, from the
definition of $+_r$ in $K^{\*_r}$ and from \Ref{Fdb18}, it follows that
                $$\multline \langle v,a\rangle = 
\langle (v -_l(f,0))+_l(f,0), (e,0)+_r(0,f^\*)\rangle \\ =
\langle (v-_l(f,0)),(e,0)\rangle + \langle (f,0), (0,f^\*)\rangle =
\langle f,f^\* \rangle .
                                     \endmultline
												$$
         The proof of the second equality is analogous.
        \endproof

        Now, we can write the diagram for $\bK^{\*_r}$. If $\bK$ is
represented by the diagram
                $$ \Cgaps{.8}
        \CD  & K @()\L{\zt_l} @(-1,-1) @() \L{\zt_r}  @(1,-1)&   \\
                F@() \L{\bar\zt_r} @(1,-1) & C @() \0( @(0,1) @()
\L{\zt} @(0,-1)    & E @() \L{\bar\zt_l} @(-1,-1) \\
                &  M &
        \endCD
                                                                                                \tag \label{Fdb19}$$

then the right dual  $\bK^{\*_r}$ is represented by the diagram
        $$ \Cgaps{.8}
        \CD  & K^{\*_r} @()\L{\zp_r} @(1,-1) @() \L{\zp_l}@(-1,-1) &    \\
        E @() \L{\bar\zp_r}@(1,-1) & F^\* @() \0(\ds(-2;-2) @(0,1)
@()\ds(2;2) \L{\zp} @(0,-1) &  C^\*  @()\l{\bar \zp_r} @(-1,-1)
\\
                            &M &
        \endCD \qquad .\tag\label{Fdb20}$$

        \hl2""{The left dual} A  construction, similar to that
of $\bK^{\*_r}$, can be applied to the left dual $K^{\*_l}$. As the
result, we obtain   the structure of a double
vector bundle on $K^{\*_l}$. We denote it by $\bK^{\*_l} = (\bK^{\*_l}_r,
\bK^{\*_l}_l, \bF, \bC^\*) $. The core of $\bK^{\*_l} $ is
$\bE^\*$.
        There is  an obvious identity
                $$  \ss J(\bK^{\*_l}) = (\ss J(\bK))^{\*_r} .
                                                                                                                        \tag \label{Fdb21}$$

        \hl2""{Examples}

        { 1.} Let $\bK = \bK(\bF, \bC, \bE)$. The right dual
can be canonically identified with $\bK(\bE,\bF^\*,\bC^\*)$
and the left dual with $\bK( \bC^\*, \bE^\*,\bF)$.

        { 2.} If $\bE = (E, \zt, M)$ is a vector bundle, then
$\sT \bE$ is a double vector bundle with the diagram

                $$ \Cgaps{.8}
        \CD  & \sT E @()\L{\sT\zt} @(-1,-1) @() \L{\ezt_E}  @(1,-1)&   \\
                \sT M@() \L{\ezt_M} @(1,-1) & E @() \0( @(0,1) @()
\L{\zt} @(0,-1)    & E @() \L{\zt} @(-1,-1) \\
                &  M &
        \endCD \qquad .
                                                        \tag \label{Fdb22}$$
 Its right dual $(\sT \bE)^{\*_r}$ is represented by the diagram

                $$ \Cgaps{.8}
        \CD  & \sT^\* E @()\L{\sT^\*\zt} @(1,-1) @() \L{\ezp_E}  @(-1,-1)&   \\
                E @() \l{\zt} @(1,-1) & \sT^\* M @() \0( @(0,1) @()
\L{\ezp_M} @(0,-1)    & E^\* @() \L{\zp} @(-1,-1) \\
                &  M &
        \endCD \qquad .
                                                                                                \tag \label{Fdb23}$$

        We see that the manifold of cotangent vectors to a vector
bundle has two compatible vector bundle structures.
The double vector bundle $\ss J ((\sT \bE)^{\*_r})$, represented
by the diagram

                $$ \Cgaps{.8}
        \CD  & \sT^\* E @()\L{\sT^\*\zt} @(-1,-1) @() \L{\ezp_E}  @(1,-1)&   \\
                E^\* @() \l{\zp} @(1,-1) & \sT^\* M @() \0( @(0,1) @()
\L{\ezp_M} @(0,-1)    & E @() \L{\zt} @(-1,-1) \\
                &  M &
        \endCD \qquad  ,\tag\label{Fdb24}$$
will be denoted by $\sT^\* \bE$. In particular, for $\bE = \sT M$,
the diagram \Ref{Fdb24} assumes the form

                $$ \Cgaps{.8}
        \CD  & \sT^\* \sT M @()\L{\sT^\*\zt_M} @(-1,-1) @()
\L{\ezp_{\ssT M}}  @(1,-1)&   \\
                \sT^\* M @() \l{\zp_M} @(1,-1) & \sT^\* M @() \0( @(0,1) @()
\L{\ezp_M} @(0,-1)    & \sT M @() \L{\zt_M} @(-1,-1) \\
                &  M &
        \endCD \qquad  .\tag\label{Fdb25}$$

        { 3.} In order to identify the structure of the left
dual of $\sT \bE$ let us recall that the dual to the vector
bundle $(\sT E, \sT \zt, \sT M )$ can be canonically identified
with $(\sT E^\*, \sT \zp, \sT M)$ ([{\bf 5}]). It follows that
the left dual to $\sT \bE$ is canonically isomorphic to  $ \ss
J(\sT \bE^\*)$ with the diagram

                $$ \Cgaps{.8}
        \CD  & \sT E^\* @()\L{\sT\zp} @(1,-1) @() \L{\ezt_{E^\s*}}  @(-1,-1)&   \\
        E^\* @() \L{\zp}         @(1,-1) & E^\* @() \0( @(0,1) @()
\L{\zp} @(0,-1)    & \sT M@() \L{\ezt_M} @(-1,-1) \\
                &  M &
        \endCD \qquad .
                                                                                                \tag \label{Fdb26}$$

        \hl1{Dual morphisms}

Let $\bold A = (A,\zt, M)$ and $\bold A' = (A',\zt', M')$ be
vector bundles and let $\zf \colon A\rightarrow A'$ and
$\overline{\zf} \colon M\rightarrow M'$ be mappings such that the
pair $\hbox{\pmib \char'047} = (\zf, \overline{\zf})$ is a morphism of vector
bundles. The dual morphism $\hbox{\pmib \char '047}^\* \colon \bold A'
\rightarrow \bold A$ is not composed of mappings, but relations,
unless it is an isomorphism.

        In order to avoid the use of relations instead of
mappings, we restrict further considerations to the case of
isomorphism only.  It follows from the formulae \Ref{Fdb10} that
$\bold \zF$ is an isomorphism of double vector bundles if and
only if $\bar \bold \zF_r, \bar \bold \zF_l $ and $\bar \bold
\zF_c$ are isomorphisms of vector bundles.

        Let $ \bold \zF \colon \bK \rightarrow \bK'$ be an
isomorphism of double vector bundles, $\bold \zF = (\zF, \zF_r,
\zF_l, \bar \zF)$. Consequently, $\bold \zF_r = (\zF, \zF_r)$ is
an isomorphism of vector bundles,
                $$ \bold \zF_r \colon \bK_r \rightarrow \bK_r'.
                                                                                $$
The dual vector bundle morphism
                $$ (\bold \zF_r)^\* \colon (\bK'_r)^\* \rightarrow (\bK_r)^\*
                                                                                $$
is also an isomorphism.  Let us denote by $\zF^{\*_r}$  the  mapping of bundle
spaces,  $\zF^{\*_r} \colon  {K'}^{\*_r} \rightarrow K^{\*_r}$.
We have
                $$ (\bold \zF_r)^\*  = (\zF^{\*_r}, \zF_r^{-1}).
                                                                                $$
        \claim{Proposition}{}
        $\zF^{\*_r}$ defines an isomorphism of vector bundles ${\bK'}^{\*_r}_r$
and $\bK^{\*_r}_r$.
        \endclaim
        \proof
        First,  we show that $\zF^{\*_r}$ respects fibrations
$\zp'_r$ and $\zp_r$. Let  $a,b \in {\bK'}^{\*_r} $  be such
that $\zp'_r(a) = \zp'_r(b)$. The kernel $\ker \zt'_l$ is, with
the right bundle structure, isomorphic to $E'\times _M\bC'$ and,
with the left bundle structure, it is isomorphic  to
$\bE'\oplus_M\bC'$. The equality $\zp_r'(a) = \zp'_r(b)$ means that
for  $k'\in C $, we have
                $$ \langle (\zp_l'(a),k'), a\rangle =  \langle (\zp_l'(b),k'), b
\rangle
                                                \tag\label{Fdb27}               $$
        Since $\zF$ is linear with respect to both, the right and left, vector bundle
structures, we have
                $$ \zF((e,k)) = \zF((e,0)) + \zF((0,k)) = \zF_c(k)
                                                                         $$
        for $(e,k) \in \ker \zt_l$.
        It follows that
                $$ \multline \langle (\zp_l(\zF^{\*_r}(a)),k), \zF^{\*_r}(a)
\rangle  =  \langle \zF (\bar \zF^{-1}(a),k), a\rangle = \\
        = \langle (\bar \zF^{-1}(a),\zF_c(k)),a\rangle  =\langle (\bar
\zF^{-1}(b),\zF_c(k)),b\rangle =
\langle (\zp_l(\zF^{\*_r}(b)),k), \zF^{\*_r}(b)\rangle
                                \endmultline
                                                                \tag \label{Fdb28}$$
and
                $$ \zp_r(\zF^{\*_r}(a)) = \zp_r(\zF^{\*_r}(b)).
                                                                                $$

        Linearity of $\zF^{\*_r}$ with respect to the right vector bundle
structure follows  from
                $$\multline \langle v+_lw, \zF^{\*_r}(a+_rb) \rangle = \langle
\zF(v+_lw), a+_rb\rangle  = \langle \zF(v)+_l\zF(w),a+_rb\rangle \\
        = \langle \zF(v),a\rangle +\langle \zF(w),b\rangle  =\langle
v,\zF^{\*_r}(a)\rangle +\langle w, \zF^{\*_r}(b)\rangle = \langle  v+_lw,
\zF^{\*_r}(a)+_r\zF^{\*_r}(b)\rangle .
                                \endmultline                            $$
        \endproof

                \claim{Corollary}{} \label{Cdb7}
        We have the following equalities for $\bold \zC = \bold \zF^{\*_r}$
                $$\split
                        \bar \bold\zC_r &= (\bar \bold \zF_c)^\* \\
                        \bar \bold\zC_l &= (\bar \bold \zF_r)^{-1} \\
                        \bar \bold\zC_c &= (\bar \bold \zF_l)^\*
                        \endsplit                                       \tag \label{Fdb29}$$
                                        \endclaim
        \proof
        $\ker \zt_r$ and $ \ker \zt_r'$, with the right vector
bundle structures, are isomorphic to $\bF \oplus_M \bC$ and
$\bF' \oplus_{M'} \bC' $ respectively. The existence of these
isomorphisms implies that the restriction of $\bold\zF$ to
$\ker\zt_r$  is identifiable with  $\bar \bold \zF_l \oplus \bar
\bold \zF_c $. It follows that, with these identifications, the
dual mapping
            $$ \bold \zF ^{\*_r} \colon  (\bC')^\* \oplus_{M'}
(\bF')^\* \rightarrow (\bC)^\* \oplus_{M} (\bF)^\*
                                                   $$
         is equal $ (\bar \bold \zF_c)^\* \oplus (\bar \bold
\zF_l)^\*$.

We have also that $\ker \zp_l$ and $\ker \zp'_l$, with the left
vector bundle structure,  are isomorphic to   $(\bC)^\*
\oplus_{M} (\bF)^\*$  and $(\bC')^\* \oplus_{M} (\bF')^\*$
respectively. With these isomorphisms the restriction of $\bold
\zF ^{\*_r}$ to $\ker \zp'_l$ is identifiable with $\bar \bold
\zC_r \oplus \bar \bold \zC_c$. It follows
that  $\bar \bold\zC_r = (\bar \bold
\zF_c)^\* $, $\bar \bold\zC_c = (\bar \bold \zF_l)^\* $.
The third equality in \Ref{Fdb29} is obvious.
        \endproof

        \hl2""{Examples}

        { 1.} Let  $\bK = \bK(\bF, \bC, \bE)$, $\bK' =
\bK(\bF', \bC', \bE')$, and let $\bold \zF \colon \bK \rightarrow \bK' $, with
                $$ \zF(f,c,e) = ( \zF_l(f), \zF_c(c) + \zC(f,e), \zF_r(e)).
                                           \tag \label{Fdb30}$$

                Since we identify $\bK^{\*_r}$ with $\bK( \bE,
\bF^\*, \bC^\*)$  and  $(\bK')^{\*_r}$ with $\bK( \bE',
(\bF')^\*, (\bC')^\*)$, then the morphism $\bold \zF^{\*_r}$ is
identified as a morphism
                $$\bold \zF^{\*_r} \colon \bK( \bE', (\bF')^\*,
(\bC')^\*) \rightarrow \bK( \bE, \bF^\* ,\bC^\*) .
                                                    $$
        One can easily verify the equality
                $$ \zF^{\*_r} = ( \zF^{\*_r}_l,
\zF^{\*_r}_c + \zC^{\*_r}, \zF^{\*_r}_r),
                                                                                                                                $$
where
                $$\split \zF^{\*_r}_r(q') &= \zF_c^\*(q'), \\
                                \zF^{\*_r}_l(e') &= \zF_r^{-1}(e'), \\
                                \zF^{\*_r}_c(p') &= \zF_l^\*(p'), \\
                                \zC^{\*_r}(e',q') &= \zC^{\*}(q', \zF_r^{-1}(e')),
                                        \endsplit                                               \tag \label{Fdb31}$$
and $\zC^{\*}$ is the vector bundle mapping dual to $\zC$ with respect to
the left argument, i.e., with respect to $f\in F$.

        { 2.}([{\bf 5}])  We have (Section \Ref{s1}) an
isomorphism of double vector bundles
                $$ \ezk_M \colon \sT \sT M \rightarrow \ss J (\sT\sT M).
                                                                                                                                $$
        The right dual
                $$ (\ezk_M)^{\*_r} \colon \sT \sT^\* M \rightarrow
\sT^\* \sT M
                                                                                                                                $$
is usually denoted by $\eza_M$ and plays a crucial role in the
Lagrangian formulation of the dynamics of mechanical systems.

        \hl1{Canonical isomorphisms}
        \claim{Proposition}{}
        The three double vector bundles
                $$ (\bK^{\*_r})^{\*_l},\ (\bK^{\*_l})^{\*_r}, \
\text{and}\ \bK
                                                                \tag\label{Fdb32}               $$
        are canonically isomorphic.
        \endclaim
        \proof
        It follows from the construction in Section~4 that we can identify manifolds
$(K^{\*_r})^{\*_l}$ and $K$. Also the right vector bundle structures
coincide. Let $\zF \colon K \rightarrow (K^{\*_r})^{\*_l}$ be the canonical
diffeomorphism and let $\zy_r, \zy_l$ be projections in
$(\bK^{\*_r})^{\*_l}$. For $f^\* \in F^\*  = \ker \zp_r \cap \ker \zp_l$,
we have
                $$\langle f^\*,\zy_l(\zF(v))\rangle  = \langle
(\zt_r(v),f^\*),\zF(v)\rangle  = \langle v,(e,f^\*)\rangle = \langle
\zt_l(v), f^\* \rangle ,
                                                                                $$
        hence $\zy_l(\zF(v)) = \zt_l(v)$.
        Equality of the left vector bundle structures follows from
                $$ \langle a +_rb,\zF(v)+_l \zF(w) \rangle  = \langle
a, \zF(v)\rangle + \langle b,\zF(w) \rangle =\langle v,a \rangle + \langle
w, b\rangle = \langle v+_lw , a+_rb\rangle .
                                                                                $$
        \endproof

        In the following, we show that there is a canonical
isomorphism of double vector bundles $\bK$ and $
((\bK^{\*_r})^{\*_r})^{\*_r}$.  Through this section we denote
by $(\zt_r,\zt_l),\, (\zp_r,\zp_l),\, (\zx_r,\zx_l),$ and
$(\zy_r,\zy_l)$  the right and left projections in $\bK$,
$\bK^{\*_r}$, $(\bK^{\*_r})^{\*_r}$, and
$((\bK^{\*_r})^{\*_r})^{\*_r}$ respectively.
        Identifying  vector bundles with their  second duals, we have

                $$\split \zx_r \colon (K^{\*_r})^{\*_r} \rightarrow F &\qquad \zx_l
\colon (K^{\*_r})^{\*_r} \rightarrow C^\* \\
        \zy_r \colon ((K^{\*_r})^{\*_r})^{\*_r} \rightarrow E &\qquad \zy_l \colon
((K^{\*_r})^{\*_r})^{\*_r} \rightarrow F.
                        \endsplit                                                       $$
        The core of $(\bK^{\*_r})^{\*_r}$  is $\bE^\*$ and the core of
$((\bK^{ \*_r})^{\*_r})^{ \*_r}$ is $(\bC^\*)^\* = \bC$.

        We define a relation $\Cal R_K \subset K \times
((K^{\*_r})^{\*_r})^{\*_r}$ in the following way:

Let  $v\in K,\ \zf \in ((K^{\*_r})^{\*_r})^{\*_r}$ be such that $\bar
\zt_l(\zt_r(v)) = \bar \zy_l(\zy_r( \zf))$. We say that $(v,\zf) \in \Cal
R_K$ if for each $a\in K^{\*_r}$ and $\za \in (K^{\*_r})^{\*_r}$ such that
                $$ \zt_r(v) =\zp_l (a) ,\ \zp_r (a) = \zx_l (\za), \ \zx_r (\za)
= \zy_l(\zf),
                                                                                $$
we have
                $$ \langle a, \za\rangle = \langle v, a \rangle + \langle \za ,
\zf\rangle .
                                                                \tag \label{Fdb33}$$
        \claim{Theorem}{}
        The relation $\Cal R_K$ is a mapping which defines an
isomorphism of double vector bundles.
                \label{Cdb8} \endclaim
        \proof
        Let  $(v,\zf)\in \Cal R_K$, $e = \zt_r(v)$, and $f = \zy_l(\zf)$. We
show first that $\zt_r(v) = \zy_r(\zf)$ and $\zt_l(v) = \zy_l(\zf)$.  For
$a=(e,f^\*) \in \ker \zp_r$ and $\za = (f,e^\*) \in \ker \zx_l$ equality
\Ref{Fdb33} assumes the form
                $$\langle a,\za\rangle = \langle f,f^\*\rangle + \langle e,e^\*
\rangle  = \langle v,a\rangle + \langle \za,\zf\rangle = \langle \zt_l(v),
f^\* \rangle + \langle \zy_r(\zf), e^\* \rangle .
                                                                                                                        $$
        in view of \Ref{Fdb18} and Proposition~\Ref{Cdb6}.
         It follows that
                $$ e= \zy_r(\zf),\ \ \ f = \zt_l(v),
                                                                                $$
        since  $e^\*$ and $f^\*$ are arbitrary.

        In the next step we show that $\Cal R_K$ restricted to either right or
left fiber is a linear relation.
        Let  $(v,\zf), (v',\zf')\in \Cal R_K$ be such that
                $$\split \zt_r(v) = \zt_r(v') \ &\text{and}\ \
\zy_r(\zf) = \zy_r(\zf').
                        \endsplit                                                       $$
         Let  $a, \za$ be such that $\zp_l(a) = \zt_r(v)$, $\zp_r(a) =
\zx_l(\za)$ and $\zx_r(\za) = \zy_l(\zf+_r \zf').$
	We can find $\zb,\zb' \in K^{\*_r\*_r}$ such that
$\zx_r(\zb) = \zy_l(\zf)$, $\zx_r(\zb') =
\zy_l(\zf')$, and $\za = \zb +_l \zb'$. We have then
		$$ \langle v,a\rangle + \langle v',a\rangle  = \langle v
+_r v', a\rangle 
														$$
	and, by the definition of the right vector bundle structure
on $(K^{\*_r\*_r})^{\*_r}$,
		$$\langle \zb +_l \zb', \zf +_r \zf'\rangle  = \langle
\zb,\zf \rangle  + \langle \zb', \zf'\rangle. 
																$$
	 It follows that 
                $$\multline \langle a, \za\rangle  = \langle a, \zb+_l\zb'\rangle  = \langle a,\zb\rangle
+\langle a,\zb'\rangle \\
= \langle v,a\rangle + \langle \zb,\zf\rangle  +
\langle v',a\rangle +\langle \zb',\zf'\rangle = 
        =\langle v+_rv',a\rangle +\langle \za, \zf+_r\zf'\rangle,
                                \endmultline            \tag \label{Fdb34}$$
        and, consequently, $(v+_rv',\zf+_r\zf') \in  \Cal R_K$.
        Similar arguments show that $\Cal R_K$ is invariant with
respect to the left addition and also with respect to the right
and left  multiplications by a number.
        The dimensions of $K$ and $((K^{\*_r})^{\*_r})^{\*_r}$ are equal $n +
n_E +n_F + n_C$. Thus it remains to show that the kernel and cokernel of 
$\Cal R_K$ are trivial, and that the domain of $\Cal R_K$ is the whole $K$.

It will prove that $\Cal R_K$ is an injective mapping and,
consequently, an isomorphism of double vector bundles.

        Let  $(v,\zf)\in  \Cal R_K$ and $(v', \zf) \in \Cal R_K$. Since
$\zt_r(v) = \zy_r(\zf) = \zt_r(v')$ and $\zt_l(v) = \zy_l(\zf) = \zt_l(v')$,
we have
         $$ v-_rv' \in \ker \zt_l \ \ \text{and}\ \  v-_lv' \in \ker \zt_r.
                                                      $$
	Since we identify $\ker \zt_l$ and $\ker \zy_l$ with
$E\times _M C$, we can write
		$$ v -_r v' = (e,k) \ \ \text{and} \ \ \zf -_r\zf =(e,0).
																$$
	Let $a\in K^{\*_r}$ and $\za\in K^{\*_r\*_r} $ be such that
		$$\split
			\zp_l(a) &= \zt_r(v-_rv') = e,\\
			\zp_r(a) &= \zx_l(\za),\\
			\zx_r(\za) &= \zy_l(\zf -_r\zf).
			\endsplit								$$
	It follows that $\za \in \ker \zx_r = C^\*\times _M E^\*$.
Let $\za = (k^\*, e^\*)$ in this representation. We have from
Proposition~\Ref{Cdb6} that 
		$$\split
		\langle a, \za\rangle &=\langle e, e^\* \rangle,\\
		\langle v-_rv',a\rangle &= \langle k, \zp_r(a) \rangle =
\langle k, \zx_l(\za) \rangle=  \langle k,k^\* \rangle, \\
		\langle \za, \zf-_r\zf\rangle &= \langle
\zy_r(\zf-_r\zf), e^\* \rangle = \langle e,e^\*\rangle . 
		\endsplit												$$
	It follows that the equality 
		$$ \langle a,\za\rangle = \langle v-_rv',a\rangle
+\langle \za, \zf-_r\zf\rangle 
																$$
assumes the form 
		$$ \langle e,e^\*\rangle  = \langle k,k^\* \rangle  +
\langle e, e^\*\rangle.
													$$
	Hence  $ \langle k,k^\*\rangle =0$ for each $k^\*$ and,
consequently, $k=0$. The kernel of $\Cal R_K$ is trivial.
Similarly, we prove that the cokernel is trivial. 

	Now, let $v\in K$ and let $\za \in K^{\*_r\*_r}$ with
$\zx_r(\za) = \zt_l(v)$. We have to show that the formula 
		$$\langle \za, \zf \rangle = \langle a,\za\rangle  -
\langle v, a \rangle, 
																$$
	where $\zp_l(a)=\zt_r(v)$ and $\zp_r(a) = \zx_l(\za)$,
defines an element $\zf \in K^{\*_r\*_r\*_r}$. It is enough to
prove that the right hand side of this formula does not depend on
the choice of $a$.

	Let $a,a'$ be such that $\zp_l(a)=\zp_l(a')=\zt_r(v)$ and
$\zp_r(a) = \zp_r(a') =\zx_l(\za)$. Then 
		$$ a-_ra' \in \ker \zp_l = C^\*\times_ M F^\*
																$$
	and
		$$ a-_la' \in \ker \zp_r = E\times_ M F^\*.
																$$
	Let $a-_ra' =(k^\*,f^\*)$ and $a-_la' =(e, {f'}^\*)$.
The linearity of the left vector bundle structure on $K^{\*_r}$,
with respect to the right vector bundle structure implies that
		$$\multline
	(0,f^\*) = (k^\*,f^\*) -_l(k^\*,0) = (a-_ra') -_l(a'-_ra')\\
	= (a-_la') -_r(a'-_la') = (e,{f'}^\*) -_r(e,0) = (0,{f'}^\*)
			\endmultline					$$
	and, consequently, $f^\*= {f'}^\*$.
	We get from Proposition~\Ref{Cdb6}
		$$\multline
	\langle a,\za\rangle - \langle v,\za\rangle  -(\langle
a',\za\rangle - \langle v,\za'\rangle ) \\
	= \langle a-_ra',\za\rangle - \langle v,a-_l a'\rangle  =
\langle \zx_r(\za), f^\* \rangle - \langle \zt_l(v),
f^\*\rangle = \langle \zt_l(v), f^\*\rangle  - \langle \zt_l(v),
f^\*\rangle =0. 
				\endmultline					$$
	 We conclude that $\zf$ is well defined and, consequently,
$(v, \zf) \in \Cal R_K$. 
        \endproof

        If we replace the right-hand side of \Ref{Fdb33} by a different
combination of $\langle v, a\rangle $ and $\langle \za,
\zf\rangle $, we obtain another isomorphism. The isomorphism
corresponding to $\langle v,a \rangle - \langle \za,\zf\rangle $
will be denoted by $\Cal R_K^\pm$, the isomorphism
corresponding to $-\langle v,a \rangle + \langle \za,\zf\rangle $
will be denoted by $\Cal R_K^\mp$, and the isomorphism
corresponding to $-\langle v,a \rangle - \langle \za,\zf\rangle $
will be denoted by $\Cal R_K^=$.

        \claim{Proposition}{}
        \list
        \item $\zf = \Cal R^\pm_K(v)$ if and only if $\zf = \Cal
R_K((-1)\cdot _lv)$ or equivalently, $(-1)\cdot _l \zf = \Cal R_K(v)$,
        \item $\zf = \Cal R^\mp_K(v)$ if and only if $\zf = \Cal
R_K((-1)\cdot _rv)$ or equivalently, $(-1)\cdot _r \zf = \Cal R_K(v)$,
        \item $\zf = \Cal R^=_K(v)$ if and only if $(-1)\cdot _l\zf = \Cal
R_K((-1)\cdot _rv)$ or equivalently, $(-1)\cdot _l((-1)\cdot _r \zf) = \Cal R_K(v)$.
        \endlist
                \label{Cdb9} \endclaim
        \proof
        The proof is an immediate consequence of the equalities
                $$ -\langle v,a \rangle = \langle (-1)\cdot _r v,
a\rangle =\langle v, (-1)\cdot _l a\rangle ,
                                                                                                                        \tag \label{Fdb35}$$
                $$ -\langle \za, \zf \rangle = \langle (-1)\cdot _r \za,
\zf \rangle =\langle \za, (-1)\cdot _l \zf\rangle ,
                                                                                                                        \tag \label{Fdb36}$$
          and the fact that $\Cal R_K$ defines  an isomorphism of
double vector bundles.
        \endproof

        In an analogous way we introduce isomorphisms  $\Cal L_K,
\Cal L_K^\pm, \Cal L_K^\mp, \Cal L_K^=$ of $\bK^{\*_l\*_l\*_l}$
and $K$.
        \hl2""{Examples}

        { 1.}
        Let  $\bK =\bK( \bF, \bC,\bE)$. Using the identification
$\bK^{\*_r} = \bK(\bE, \bF^\*,\bC^\* )$ and also the identifications $\bE^{\*\*} =\bE$,
$\bF^{\*\*}= \bF$, and $\bC^{\*\*} = \bC$, we get
                $$ \bK^{\*_r\*_r} = \bK( \bC^\*, \bE^\*, \bF)
                                                                                                                                $$
and
                $$ \bK^{\*_r\*_r\*_r} = \bK(\bF, \bC, \bE).
                                                                                                                                $$
        Thus, we have obtained another identification of $\bK$ and
$\bK^{\*_r\*_r\*_r}$. With this identification the formula
\Ref{Fdb33} assumes the form
                $$  \split
        \langle (e,\zf,\zg), ( \zg,\ze,\bar f)\rangle  = \langle (f,c,e),
(e,\zf,\zg) \rangle  + \langle ( \zg, \ze,\bar f), ( \bar
f, \bar c,\bar e)\rangle ,\\
        \langle e,\ze \rangle  +\langle \bar f, \zf \rangle  = \langle f,
\zf\rangle + \langle c,\zg\rangle + \langle \bar c, \zg\rangle +
\langle \bar e, \ze\rangle ,
                                \endsplit                               \tag\label{Fdb37}$$
        where $e, \bar e\in E, f, \bar f\in F, c, \bar
c\in C, \ze \in E^\*, \zf\in F^\*, \zc \in C^\*$. Hence, $e=\bar
e, f =\bar f, c = -\bar c$, and, consequently,
                $$ \Cal R_K (f,c,e) = (f,-c,e).
                                                                                                                        \tag \label{Fdb38}$$
        Analogously,
                $$\split
                \Cal R^\pm_K(f,c,e) & = (f,c,-e) \\
                \Cal R^\mp_K(f,c,e) & = (-f,c,e) \\
                \Cal R^=_K(f,c,e) & = (-f,-c,-e).
                                \endsplit                                                                       \tag \label{Fdb39}$$

        { 2.}
Let  $\bK$ be the double vector bundle $\sT^\* \bE$ represented
by the diagram

                $$ \Cgaps{.8}
        \CD  & \sT^\* E @()\L{\sT^\*\zt} @(-1,-1) @() \L{\ezp_E}  @(1,-1)&   \\
                E^\* @() \l{\zp} @(1,-1) & \sT^\* M @() \0( @(0,1) @()
\L{\ezp_M} @(0,-1)    & E @() \L{\zt} @(-1,-1) \\
                &  M &
        \endCD \qquad  .\tag\label{Fdb40}$$
        Then the first, second, and third right duals can be identified
with double vector bundles $\ss J(\sT E),\ \sT E^\* $ and $\ss
J(\sT^\* E)$, represented by  diagrams
                        $$ \Cgaps{.8}
        \CD  & \sT E @()\L{\sT\zt} @(1,-1) @() \L{\ezt_E}  @(-1,-1)&   \\
                E @() \l{\zt} @(1,-1) & E  @() \0( @(0,1) @()
\L{\zt} @(0,-1)    & \sT M @() \L{\ezt_M} @(-1,-1) \\
                &  M &
        \endCD \ \
        \CD  & \sT E^\* @()\L{\sT\zp} @(-1,-1) @() \L{\ezt_{E^\s*}}  @(1,-1)&   \\
                \sT M@() \l{\ezt_M} @(1,-1) & E^\* @() \0( @(0,1) @()
\L{\zp} @(0,-1)    & E^\* @() \L{\zp} @(-1,-1) \\
                &  M &
        \endCD \ \
        \CD  & \sT^\* E^\* @()\L{\sT^\*\zp} @(1,-1) @() \L{\ezp_{E^\s*}}  @(-1,-1)&   \\
        E^\* @() \l{\zp} @(1,-1) & \sT^\* M @() \0( @(0,1) @()
\L{\ezp_M} @(0,-1)    & E @() \l{\zt}  @(-1,-1) \\
                &  M &
        \endCD \quad  .\tag\label{Fdb41}$$
        Canonical isomorphisms $\Cal R_K,\ \Cal R_K^\pm, \ \Cal
R_K^\mp,\ \Cal R_K^= $ define   diffeomorphisms from $\sT^\*
E^\*$ to $\sT^\* E$.  These diffeomorphisms are
antisymplectomorphisms with respect to the canonical symplectic
structure of the cotangent bundle for $\Cal R_K,\ \Cal R_K^=$
and symplectomorphisms for  $\Cal R_K^\pm,\ \Cal R_K^\mp$.

\hl2""{Identification of isomorphisms}
        Let $\bold \zF \colon \bK \rightarrow \bK' $ be an
isomorphism of double vector bundles. We have
                $$ \bold \zF^{\*_r\*_r\*_r} \colon  (\bK')^{\*_r\*_r\*_r}
\rightarrow \bK^{\*_r\*_r\*_r}.
                                                                                                $$
Using one of the introduced isomorphisms of a double vector
bundle and its third right dual, we can compare  $\bold \zF$ and
its third right dual.

        \claim{Proposition}{} \label{Cdb10}
We have the following equality
                $$  \Cal R^{-1}_K \circ  \zF^{\*_r\*_r\*_r} \circ \Cal
R_{K'} =  \zF^{-1}.
                                                           \tag \label{Fdb42}$$
\endclaim
        \proof Let  $\zf\in K^{\*_r\*_r\*_r}$, $\zf'\in
(K')^{\*_r\*_r\*_r}$, $v\in K$,  and $v'\in K'$ be such that
                $$ \zf = \zF^{\*_r\*_r\*_r}(\zf') ,\ \ \zf =\Cal R_K(v),
\ \ \zf' = \Cal R_{K'}(v').
                        \tag\label{Fdb43}                                       $$
Then, as in \Ref{Fdb33}, we have
                $$ \langle a', \za'\rangle  = \langle v', a'\rangle +
\langle \za', v\rangle
                                                                                                                        $$
          and, for $a = \zF^{\*_r}(a'),\ \za = (\zF^{\*_r\*_r})^{-1}
(\za')$, we have
                $$ \langle (\zF^{\*_r})^{-1}(a), \zF^{\*_r\*_r} (\za)
\rangle = \langle v', (\zF^{\*_r})^{-1}(a) \rangle  + \langle
\zF^{\*_r\*_r}(\za), \zf'\rangle ,
                                    \tag \label{Fdb44}      $$

                $$ \langle \zF^{\*_r\*_r\*_l} \circ (\zF^{\*_r})^{-1}
(a),\za \rangle  = \langle ((\zF^{\*_r})^{-1})^{\*_l} (v'), a\rangle
+ \langle \za, \zF^{\*_r\*_r\*_r} (\zf')\rangle .
                                                                                                \tag \label{Fdb45}$$
        From the formula \Ref{Fdb33} and from \Ref{Fdb45}, we derive the
following identities
                $$\split \zF^{\*_r\*_r\*_l} \circ (\zF^{\*_r})^{-1} (a) &=
a, \\
                (\zF^{\*_r})^{-1})^{\*_l} (v') &= \zF^{-1} (v'),\\
                \zF^{\*_r\*_r\*_r} (\zf') &= \zf
                                \endsplit                                                               $$
which make the formula \Ref{Fdb45} equivalent to
                $$ \langle a,\za \rangle  = \langle \zF^{-1} (v'), a\rangle
+ \langle \za, \zf\rangle .
                                                         $$
        It follows that
                $$ \Cal R_K (\zF^{-1}(v')) =\zf = \Cal R_K(v).
                                                                                                                $$
        Consequently,
                $$ \zF^{-1}(v') = v.
                                                          $$
        \endproof
        Equalities similar to \Ref{Fdb42} hold for pairs of relations
$(\Cal R_K^\mp, \Cal R_{K'}^\mp)$, $(\Cal R_K^\pm, \Cal
R_{K'}^\pm)$, and $(\Cal R_K^=, \Cal R_{K'}^=)$ in place of
$(\Cal R_K, \Cal R_{K'})$.

        \hl2""{Remark} In the case of $\bK = \bK(\bF,\bC,\bE)$
and $\bK' = \bK(\bF',\bC',\bE')$ we  have another isomorphisms
of $\bK$ and $\bK^{\*_r\*_r\*_r}$, $\bK'$ and
$(\bK')^{\*_r\*_r\*_r}$ (see Example 1 of this section).  In
contrast to \Ref{Fdb42}, $\zF^{\*_r\*_r\*_r}$ does not
correspond, with respect to these isomorphisms, to $\zF^{-1}$.

        \hl1{Examples and applications}

        \hl2{Vector and co-vector fields on a vector bundle}
        Let $\bE = (E,\zt, M)$ be a vector bundle and let $X$ be a
vector field on $E$. By $\widetilde{X}$ we denote the
corresponding function on $\sT^\* E$.

        \claim{Theorem}{} \label{Cdb11}
        The following three conditions are equivalent.
        \list
        \item For each  function $f$ on $E$, linear on fibers, the
function  $\langle X, \xd f\rangle $ is linear on fibers.
        \item The mapping $X\colon  E\rightarrow \sT E$ is a vector bundle
morphism from $\bE$ to $(\sT E, \sT \zt, \sT M)$.
        \item The function $\widetilde{X}$ is linear with respect to
the right and left vector bundle structures on $\sT^\* E$.
        \endlist
                                        \endclaim
        \proof
 Let $f$ be a linear function on $E$ and let $v,v'\in \sT E$ be
such that $\sT \zt v = \sT \zt v' $. We can choose  curves
$\zg,\zg'$ on $E$ which represent $v,v'$, respectively, and
satisfy the following condition: $\zt\circ \zg = \zt\circ \zg'$.
The curve $\zg+\zg'\colon t\mapsto \zg(t) +\zg'(t)$ represents
the vector  $v+_l v'$. Since $f$ is linear, we have $f\circ (\zg
+\zg')(t) = f\circ \zg(t) + f\circ \zg'(t)$ and, consequently,

                $$ 
        \langle v +_l v', \xd f (e+e')\rangle =\langle v,\xd
f(e)\rangle  + \langle v',\xd f(e')\rangle = 
\langle v +_l v', \xd f (e)+_l \xd f(e') \rangle ,
                                       $$

     i.e.,
                $$ \xd f(e +e') = \xd f(e) +_l \xd f(e').
                                                                                                                        \tag \label{F1}$$

\list
 \item"($1 \Rightarrow 2$)" First, we have to show that, if
      $\zt(e)=\zt(e')$, then $\sT\zt(X(e))=\sT\zt(X(e'))$.  It is
enough to show that for each function $g$ on $E$, which is constant
on fibers,
                $$\langle X(e),\xd g(e)\rangle=
     \langle X(e'),\xd g(e')\rangle .
                                                                \tag \label{Fdb46}$$

Let  $f$ be a linear function on $E$ and and let $g$ be a
function on $E$, constant on fibers. The function $fg$ is linear  on
fibers, hence   $X(fg)$ is also linear. Since
     $$ X(fg)=fX(g)+gX(f)$$ and  $gX(f)$ is  linear on fibers, it
follows that $fX(g)$ is   linear and, consequently,  $X(g)$ is
constant on fibers and
     $$ X(g)(e)=X(g)(e').$$
        Since
        $$ X(g)(e)=\langle X(e),\xd g(e)\rangle,
                                                $$
we get \Ref{Fdb46}.

     Now, we show that $X$ is linear on fibers. Let $\zm$ be a
covector from  $\sT^{\*}_{e+e'}E$, where $e+e' \neq 0$. There
exists a function $F_{\zm}$ on $E$ which is linear on fibers and
such that $\xd F_{\zm}(e+e')=\zm$.

We have from \Ref{F1}
     $$\multline
     \langle X(e+e'),\xd F_{\zm}(e+e')\rangle=
     X(F_{\zm})(e+e')=X(F_{\zm})(e)+X(F_{\zm})(e')= \\
     =\langle X(e), \xd F_{\zm}(e)\rangle+\langle X(e'),\xd F_{\zm}(e')\rangle=
     \langle X(e) +_lX(e'),\xd F_{\zm}(e+e')\rangle
                                             \endmultline $$
     The above calculation gives, for every $\zm$,
     $$\langle X(e+e'),\zm\rangle=\langle X(e) +_lX(e'),\zm\rangle
                                        $$
        and, consequently $X(e+e') = X(e) + X(e')$ for $e+e'\neq 0$. By
the continuity argument, we get the desired equality for all $e,e'$.

     Similar arguments show that $\langle X(\za e)= \za\cdot_l X(e)$.

     \item"($2\Rightarrow 3$)" Let  $\zt(e) = \zt(e')$. Since
$X$ is a vector bundle morphism, we have
     $$\align &\zt(e)=\zt(e')\Rightarrow\sT\zt(X(e))=\sT\zt(X(e')) \\
     \intertext{and}
     & X(e+e')=X(e) +_l X(e'), \\
     & X(\za e)=\za \cdot_l X(e).
     \endalign $$
     Let $\zm$ and $\zn$ be two covectors on $E$,  $\zm\in
\sT^{\*}_e E$ and $\zn \in \sT^{\*}_{e'} E$, such that
$\sT^\*\zt \zm =\sT^\*\zt \zn$.
From the definition of the vector bundle structure on $(\sT^\*
E, \sT^\* \zt, E^\*)$, we have
     $$
     \langle X(e+e'),\zm +_l\zn\rangle=\langle X(e) +_l X(e'),
     \zm +_l\zn\rangle =\langle X(e),\zm\rangle+\langle X(e'),\zn\rangle
     $$
     and
     $$
     \langle X(\za p),\za\cdot_l\zm\rangle=
     \langle \za\cdot_l X(e),\za\cdot_l\zm\rangle= 
     \za\langle X(e),\zm\rangle
     $$
     The above calculation shows that $X$, treated as a function on
     $\sT^{\*}E$, is linear with respect to the vector bundle
structure $(\sT^\* E, \sT^\* \zt, E^\*)$. It is obviously linear
with respect to the canonical vector bundle structure  on
$(\sT^\* E, \ezp_E, E)$.

     \item"($3\Rightarrow 1$)"

        It follows from \Ref{F1} that for a linear function $f$
     $$     \langle X(e+e'),\xd f(e+e')\rangle= \langle X(e)+_l X(e'),
     \xd f(e) +_l \xd f(e')\rangle= \langle X(e),\xd f(e)\rangle +
\langle X(e'),\xd f(e')\rangle.
                                                                 $$
\endlist
\endproof

We say that a vector field $X$ is {\it  of degree zero} if
one of the conditions of this theorem is satisfied.

        Let $(X^i,e^a)$ be an adapted coordinate system on $E$. A
vector field $X$ on $E$ is  of degree zero if, in local
coordinates,
                $$ X = X^i \frac{\partial }{\partial x^i} +
X^b_a e^a \frac{\partial }{\partial e^b},
                                                                \tag \label{Fdb47}$$
where $X^i,\, X^i_a$ are functions of $(x^i)$ only.

        Now, let $\zvy$ be a 1-form on $E$ and let $\widetilde{\zvy}$
be the corresponding function on $\sT E$.

        \claim{Theorem}{} \label{Cdb12}
        The following three conditions are equivalent.
        \list
        \item For each  vector field $X$ of degree zero, the
function  $\langle X, \zvy \rangle $ is linear on fibers.
        \item The mapping $\zvy \colon  E\rightarrow \sT^\* E$
is a vector bundle morphism from $\bE$ to $(\sT^\* E, \sT^\* \zt,
E^\*)$.
        \item The function $\widetilde{\zvy}$ is linear with respect to
the right and left vector bundle structures on $\sT E$.
        \endlist
                                        \endclaim
        \proof
     \list
     \item"($1\Rightarrow 2$)" Let $X$ be a vertical
     vector field on $E$, constant on fibers.
     For every linear function $f$ on $E$ the vector field $fX$
is polynomial of degree 0 (see \Ref{Fdb47}).
     The function $\langle fX,\zvy\rangle$ is then linear  on $E$
and, since
     $$ \langle fX,\zvy\rangle=f\langle X,\zvy\rangle, $$
the function  $\langle X,\zvy\rangle(\cdot)$ is  constant on fibers. It follows that
$\sT^\*\zt \zvy (e)$ and $\sT^\*\zt\zvy(e')$ are equal if $\zt(e)=\zt(e')$
and, consequently, that $\zvy$ is a fiber preserving mapping from $\bE$
to $(\sT^\* E, \sT^\* \zt, E^\*)$. In order to  prove the linearity
of $\zvy$ on  fibers, i.e., that
     $$\align
     &\zvy(e) +_l \zvy(e')=\zvy(e+e') \\
\intertext{and}
     &\za\cdot_l\zvy(e)=\zvy(\za e),
     \endalign $$
it is enough to prove the first equality for $e,e'$   such that $e+e'\neq 0$.

        Let $e+e'\neq 0$. For every vector $v\in\sT_{e+e'}E$
there   exists a  vector field $X$ of degree 0,
such that $X(e+e')=v$. We have then
     $$\align &\langle v,\zvy(e+e')\rangle=\langle X(e+e'),\zvy(e+e')\rangle= \\
     &=\langle X(e),\zvy(e)\rangle+\langle X(e'),\zvy(e')\rangle=
     \langle X(e)+_lX(e'),\zvy(e)+_l\zvy(e')\rangle= \\
     &=\langle X(e+e'),\zvy(e)+_l\zvy(e')\rangle=
     \langle v,\zvy(e)+_l\zvy(e')\rangle.
     \endalign $$
     Since it holds for every $v$ in $\sT_{e+e'}E$, we have
     $$\zvy(e+e')=\zvy(e)+_l\zvy(e').$$

     \item"($2\Rightarrow 3$)" Let $v\in \sT_eE$, $w\in \sT_{e'}E$,
     $\zt(e)=\zt(e')$ and $\sT\zt(v)=\sT\zt(w)$. We have
     $$
     \widetilde{\zvy} (v+_l w) = \langle v +_l w,\zvy(e+e')\rangle=
     \langle v +_l w,\zvy(e) +_l \zvy(e')\rangle= 
     \langle v,\zvy(e)\rangle+\langle w,\zvy(e')\rangle =
\widetilde{\zvy}(v) +\widetilde{\zvy}(w)
     $$
     and
     $$
     \langle \za\cdot_l v,\zvy(\za p)\rangle=
     \langle \za\cdot_l v,\za\cdot_l\zvy(e)\rangle= 
     \za\langle v,\zvy(e)\rangle,
     $$
     i.e., $\zvy$ is a linear function on $\sT E$ with respect to
     the tangent vector bundle structure on $\sT E$. Linearity
with respect to the canonical vector bundle structure on $(\sT
E, \ezt_E, E)$  is obvious.

     \item"($3\Rightarrow 1$)" Let $X$ be a  vector field on
$E$ of degree 0. It follows from Theorem~\Ref{Cdb11} that
     $$\multline
     \langle X,\zvy\rangle(e+e')=\langle X(e+e'),\zvy(e+e')\rangle=
     \langle X(e) +_l X(e'),\zvy(e+e')\rangle=
\widetilde{\zvy}(X(e) +_l X(e')) = \\
     =\widetilde{\zvy (X(e))} + \widetilde{\zvy (X(e'))}
=\langle X(e),\zvy(e)\rangle+\langle X(e'),\zvy(e')\rangle=
     \langle X,\zvy\rangle(e)+\langle X,\zvy\rangle(e')
     \endmultline $$
     Similarly,
     $$
     \langle X,\zvy\rangle(\za e)=\langle X(\za e),\zvy(\za e)\rangle=
     \langle \za\cdot_l X(e),\zvy(\za e)\rangle= 
     \za\langle X(e),\zvy(e)\rangle=\za\langle X,\zvy\rangle(e).
     $$

        \endproof
\endlist
        We say that $\zvy $ is  of degree 1 (linear)
if one of the equivalent conditions of the theorem above is satisfied.

In a local coordinate system, a linear 1-form  $\zvy$ has the following form
                $$ \zvy = \zvy_a\xd e^a + \zvy_{ia}e^a \xd x^i ,
                                                                \tag \label{Fdb48}$$
         where $\zvy_a,\, \zvy_{ia} $  are functions of $(x^i)$ only.

        \hl2{Linear  Poisson structures}
        A Poisson structure on a vector bundle $\bE$ is
called linear if,  for every two functions $f,g$, linear on
fibers of $\bE$, the Poisson bracket $\{f,g\}$ is also linear on
fibers. It follows that for $f$ linear on fibers and $g,g'$
constant on fibers the bracket $\{f,g\}$ is costant on fibers and
$\{g,g'\} =0$.

        Let $\zL$ be a Poisson bivector field on $\bE$ and let
$\widetilde{\zL} \colon \sT^\* E \rightarrow  \sT E$ be the
correspondig mapping of vector bundles.
        \claim{Proposition}{}
\label{Cdb13} $\zL$ defines  a linear Poisson structure on $\bE$
if and only if $\widetilde{\zL}$ defines a morphism of double
vector bundles $\sT^{\*} \bE\rightarrow\sT\bE$.
        \endclaim
\proof

Let  $\zL$ be the bivector field of a linear Poisson structure.
In the proof of Theorem~\Ref{Cdb11} we
have shown that for a linear function $f$ on $E$ the differential
$\xd f$ is linear, i.e.,
                $$ \xd f(e+e') = \xd f(e) +_l \xd f(e').
                                                                \tag \label{Fdb49}$$

The vector field  $\widetilde{\zL}(\xd f)$ is a linear vector field because it satisfies
the condition (1) from the theorem \Ref{Cdb11}. It means that
            $$ \widetilde{\zL}(\xd f(e) + \xd f(e')) =
\widetilde{\zL}(\xd f(e+e'))=\widetilde{\zL} (\xd f(e))+_l
\widetilde{\zL}(\xd f(e')).
                                                                \tag \label{Fdb50}$$
        On the other hand, for a function $g$, constant on
fibers, the vector field $\widetilde{\zL}(\xd g)$ is vertical
and constant on fibers (we identify spaces of vertical vectors
at different points in a fiber).
        For every pair of covectors $\zm,\zn$ such that
$\zm\in\sT^{\*}_e E$, $\zn\in\sT^{\*}_{e'} E$ and they have
the same projection on $\bE^{\*}$, there exist a linear function
$f$ and a function $g$, costant on fibers, such that $\xd f(e)=\zm$
and  $\xd f(e') + \xd g(e')=\zn$. Therefore, we have, in view of \Ref{Fdb50}
            $$\multline
                \widetilde{\zL}(\zm+_l\zn)= \widetilde{\zL} (\xd
f(e) +_l(\xd f(e') + \xd g(e')) ) = \widetilde{\zL} ((\xd f(e)
+_l(\xd f(e')) + \xd g(e+ e')) )\\ = \widetilde{\zL} (\xd f(e)
+_l(\xd f(e')) + \widetilde{\zL}(\xd g(e'+e) ) =
        \left(\widetilde{\zL} (\xd f(e))+_l \widetilde{\zL}(\xd f(e'))
\right) + \widetilde{\zL}(\xd g(e'+e) )\\ = \widetilde{\zL}
(\xd f(e))+_l\left( \widetilde{\zL}(\xd f(e'))  + \widetilde{\zL}(\xd
g(e') )\right) = \widetilde{\zL}(\zm)+_l \widetilde{\zL}(\zn).
                \endmultline $$

        Now, let the Poisson bivector  $\widetilde{\zL}$ be a
morphism of double vector bundles. It follows that, for every
pair of covectors $\zm, \zn$ such that their left projections
are equal, we have
                $$ \zL(\zm +_l \zn)=\zL(\zm)+_l\zL(\zn).$$
        Let $f,g $ be linear on fibers. Consequently,
$\xd f, \xd g$ are linear one forms.
         We have then
                $$\split
              \{f,g\}(e+e')&= \langle\widetilde{\zL}(\xd f(e+e')), \xd g(e+e')\rangle \\
              &=\langle\widetilde{\zL}(\xd f(e)+_l\xd f(e')),\xd g(e)+_l\xd g(e')\rangle \\
                &=\langle\widetilde{\zL}( \xd f(e)) +_l
\widetilde{\zL} ( \xd f(e')),\xd g(e)+_l\xd g(e')\rangle\\
              &= \langle\widetilde{\zL}(\xd f(e)),\xd g(e)\rangle +
              \langle\widetilde{\zL}(\xd f(e')),\xd g(e')\rangle \\
              &= \{f,g\}(e)+\{f,g\}(e'). 
                           \endsplit $$
\endproof

        \hl2{Special symplectic manifolds}
        Let $\bE=(E,\zt,M)$ be a vector bundle and  let $\zw$ be a
2-form on $E$. By $\widetilde{\zw}$ we denote the corresponding
vector bundle morphism
                $$ \widetilde{\zw} \colon \sT E\rightarrow \sT^\* E.
                                                        \tag \label{Fdb51}$$
          We say that $\zw$ is linear with respect to the vector
bundle structure $\bE$ if $\widetilde{\zw}$ is a morphism  of double
vector bundles
                $$ \widetilde{\bold \omega} \colon \sT\bE \rightarrow \sT^\*
\bE.
                                                               \tag \label{Fdb52}$$

        If $\zw$ is linear, then there are  three derived vector bundle
morphisms:
        $$ \split
         \widetilde{\zw}_r  &\colon E \rightarrow  E ,\\
        \widetilde{\zw}_l  & \colon \sT M \rightarrow  E^\* ,\\
        \widetilde{\zw}_c  & \colon E\rightarrow \sT^\* M .
                                \endsplit               $$
        Of course, $\widetilde{\zw}_r = \id_E$ and, because
$\widetilde{\zw}$ is skew-symmetric, we have, from \Ref{Fdb29},
                $$ \widetilde{\zw}_l = - \widetilde{\zw}_c^\*.
                                                                    $$

        \claim{Proposition}{} \label{Cdb14}
        $\zw$ is closed if and only if the pull-back of the
canonical symplectic form $\zw_M$ on $\sT^\* M$  by
$\widetilde{\zw}_c$ is equal $\zw$:
                $$ \zw = \widetilde{\zw}_c^\* \zw_M.
                                                                                \tag \label{Fdb53}$$

                                        \endclaim
        \proof
        Let $(x^i, e^a)$ be a local coordinate system on $E$ and
let $(x^i,e^a,\dot{x}^j,\dot{e}^b)$, $(x^i,e^a,p_j,f_b)$   be
adopted coordinate systems on $\sT E$, $\sT^\* E$ respectively.
For a 2-form  $\zw$ on $E$,
        $$ \zw=\frac{1}{2} \zw_{ij}\xd x^i\wedge \xd x^j+\zw_{ia}\xd x^i\wedge \xd e^a+
        \frac{1}{2} \zw_{ab}\xd e^a\wedge \xd e^b,
                                                        $$
we have
        $$\split
       & p_j\circ\widetilde\zw=\zw_{ij}\dot{x}^i-\zw_{ja}\dot{e}^a, \\
       & f_b\circ\widetilde\zw=\zw_{ib}\dot{x}^i+ \zw_{ab}\dot{e}^a.
                                              \endsplit $$

        The linearity of $\zw$ implies, in view of  \Ref{Fdb10},
        $$ \split
        &\zw_{ij}(x,e)=\zw_{ija}(x)e^a ,\\
        &\zw_{ia}(x,e)=\zw_{ia}(x) ,\\
        &\zw_{ab}(x,e)=0 .
                                \endsplit               $$

The exterior derivative of $\zw$ assumes then the form:
        $$\xd\zw= \frac{1}{2} \zw_{ija,k}e^a \xd x^k\wedge \xd x^i \wedge \xd x^j+
        (\frac{1}{2} \zw_{ija}+ \zw_{ja,i}) \xd x^i\wedge \xd x^j\wedge \xd e^a.
                                                        $$
Therefore the external derivative $d\zw$ equals $0$ if and only if the following
two conditions are satisfied:
        $$ \split
        &\zw_{ija}(x)=\zw_{ia,j}(x)-\zw_{ja,i}(x) ,\\
        &\zw_{ija,k}(x)+\zw_{jka,i}(x)+\zw_{kia,j}(x)=0 .
                                \endsplit
                                                        \tag\label{Fdb54}$$
 The second condition is an immediate consequence of the first one.

On the other hand, the mapping $\widetilde{\zw}_c$ is given in
the coordinate system, by the formula
                $$ p_i = -\zw_{ia}e^a .
                                                                                $$
        Consequently, the  pull-back of the canonical symplectic
form $\zw_M = \xd p_i \wedge \xd x^i$ by $\zw_c$ is given by

       $$\widetilde\zw_c^{\*}\zw_M=\xd(-\zw_{i}e^a)\wedge \xd x^i=
       \zw_{ia}\xd x^i\wedge \xd e^a -\zw_{ia,j}e^a \xd x^j\wedge \xd x^i.
                                                        $$
It follows that  $\zw =\widetilde\zw_c^{\*}\zw_M$ if and only if
the condition
        $$\zw_{ija}(x)=\zw_{ia,j}(x)-\zw_{ja,i}(x),
                                                      \tag\label{Fdb55} $$
which is equivalent to  \Ref{Fdb54}, is satisfied.
                                        \endproof

        If $\zw $ is nondegenerate, i.e., if $\widetilde{\zw}$ is an
isomorphism of vector bundles, then also $\widetilde{ \zw}_c$ is
an isomorphism. In that case $\widetilde{ \zw}_c$ is a
symplectomorphism. Thus, we can consider the pair $(\bE, \zw)$
as a {\it special symplectic manifold} ([{\bf 6}], [{\bf 7}]).

        \hl2{Vertical lifts and complete lifts}
        In this section, we present concepts of vertical and
complete lifts of a vector field ([{\bf 9}], [{\bf 1}]) in the
general framework of  double vector bundles.

        Let $\bK=(\bK_r, \bK_l,\bE, \bF)$ be a double vector bundle with
the core $\bC$. Let $\zg$ be a section of the core. Using
the double vector bundle structure of $\bK$ we  assign to
$\zg$ two sections $\Cal V_r \zg$ and $\Cal V_l \zg$ of $\zt_r$ and
$\zt_l$ respectively.

        From Proposition~\Ref{Cdb2} we have $\ker
\zt_r=F\times_M\bC$ and $\ker \zt_l = E\times _M \bC$.
Sections  $\Cal V_r \zg$ and $\Cal V_l \zg$ are defined by the following formulae:
               $$ \Cal V_l \zg \colon  F\rightarrow K \colon  f\mapsto
               (f,\zg(\bar \zt_r(f)))\in \ker \zt_r\subset K $$
        and
               $$ \Cal V_r \zg \colon  E\rightarrow K \colon  e\mapsto
               (e,\zg(\bar \zt_l(e)))\in \ker \zt_l\subset K. $$
Sections $\Cal V_r \zg$ and $\Cal V_l \zg$ are called {\it vertical lifts
of $\zg$} with respect to the right and left projections.

 \hl3""{Examples}
        \list
        \item
 Let $\bK=\sT\bE$. The core of $\sT\bE$ is isomorphic to $\bE$
     We can therefore lift the section of $\zt$ to the section of $\zt_E$
     and $\sT\zt$. In local coordinate system $(x^i, e^a, \dot{x}^j,
     \dot{e}^b)$ vertical lifts are given by the fomulae
              $$\split
              &\Cal V_r \zg=\zg^a\frac{\partial}{\partial e^a}, \\
              &\Cal V_l \zg=\dot{x}^j\frac{\partial}{\partial x^j}
              +\zg^a\frac{\partial}{\partial e^a}.
                                                  \endsplit $$
        \item
        Let $\bK = \sT^\* \bE$. The core we identify with
$\sT^\* M$. The right vertical lift of a 1-form $\zvy$ is the
pull-back of $\zvy$ by the  projection $\zt$.
        \endlist

        Now, let $X \colon F\rightarrow K$ be a section of
$\zt_l$. We say that this section is {\it linear} if it projects
to a mapping $\overline{X} \colon M\rightarrow E$ and the pair
$\bold X =(X,\overline{X})$ is a vector bundle morphism
                $$  \bold X \colon \bF \rightarrow \bK_r.
                                                                                $$
           In a similar way we define linear sections of $\zt_r$.

\claim{Proposition}{} \label{Cdb15}
         There exist a unique linear section $Y$ of the right vector
         bundle structure of $\bK^{\*_r}$ such that $\langle X,Y\rangle=0$
\endclaim

\proof For each point $m\in M$ the image of $\bar\zt_r^{-1}(m)$
under $X$ is a vector subspace of $\zt_r^{-1}(\overline{X}(m))$.
        We denote by $X^\circ(m)$ the anihilator of this subspace in
$\zp_l^{-1} (\overline{X}(m)) \subset K^{\*_r}$. If $n_E, n_F, n_C$ are the
dimensions of fibers of $E,F,C$ respectively then the dimension
of $X^\circ(m)$ is equal $n_C$.

        We show that $X^\circ(m)$ projects to the whole fiber of $\bC^\*$
(which is of dimension $n_C$). Since the projection $\zp_r$ is
linear with respect to the left vector bundle structure,  it is
enough to show that $\zp_r$, restricted to the anihilator
$X^\circ(m)$, is an injection.

Let $a$ be an element of $\ker \zp_r \cap X^\circ(m)$. We
represent $a$ by a pair $(e,\zf)\in E\times_M F^{\*}$, where $e=\overline{X}(m)$.
Since  $(e,\zf)$ is an element of the anihilator $X^\circ(m)$ then
$\langle X(f),(\zx(m),\zf)\rangle$ should be equal to zero for all $f$.
The following calculation shows that, in this case, $\zf$ must be zero:
     $$ 0=\langle X(f),(\zx(m),\zf)\rangle=\langle \zt_l(X(f)),\zf\rangle=
     \langle f,\zf\rangle .$$
It follows that $X^\circ(m)$ is the image of a section of $\zp_r$
over $\zp^{-1}_c(m)$.
Collecting the anihilators of $X(\bar\zt_r^{-1}(m))$ point by point in $M$
we get a section $Y$ of $\zp_r$. The uniqueness of $Y$ is obvious.
	\phantom{aa} \endproof

Let $(x^i,e^a,f^A,k^\za)$ be a local coordinate system on $K$ and
let $(x^i,e^a,\zf_A,\zk_\za)$ be the adopted local coordinate system on
$K^{\*_r}$. A  section of $\zt_l$  is linear if it is of the form
$$ \split
    &e^a=\zx^a(x^i) \\
    &k^\za=X^\za_B(x^i)f^B\ .
   \endsplit $$
A linear section $Y$ of $\zp_r$ which have the same projection onto section of
$ \bar\zt_r$ can be written as:
$$ \split
    &e^a=\zx^a(x^i) \\
    &\zf_A=Y_A^\zb(x^i) \zf_\zb \ .
   \endsplit $$
The condition $\langle X,Y\rangle=0$ gives
$$ Y_A^\zb(x^i)=-X^\za_A(x^i). $$

        Of special interest is the case of $\bK = \ss J(\sT
\bE)$. The left projection is the canonical projection $\ezt_E
\colon \sT E \rightarrow E$ and a section $X$ of this
projection is a vector field on  $E$. The right dual to $\ss
J(\sT \bE)$ we identify with  $\sT \bE^\*$  (Example 3 of
Section 4) and a section of the right projection is  a vector
field on $E^\*$. Linear sections of $\ezt_E,\ \ezt_{E^\s*}$ are
vector fields of degree 0.
        Proposition~\Ref{Cdb15} establishes the ono-to-one
correspondence between  vector fields of degree 0 on
$E$ and  vector fields of degree 0 on $E^\*$. In
particular, for $\bE =\sT M$, the complete tangent lift  $\dt X$
of a vector field on $M$ is a  vector fields of degree
0 on $\sT M$ (see [{\bf 9}], [{\bf 1}]). One can easily
recognize the corresponding vector field $Y$ on the cotangent bundle
$\sT^\* M$  as the complete cotangent lift of $X$.
        In local coordinates, for $X = X^i\dfrac{\partial
}{\partial x^i} $, we have
                $$ \dt X = X^i\frac{\partial }{\partial x^i} +
X^i_{,j}{\dot x}^j \frac{\partial }{\partial {\dot x}^i}
                                                                                $$
        and
                $$ Y = X^i\frac{\partial }{\partial x^i}
- X^i_{,j} p_i\frac{\partial }{\partial p_j}.
                                                                                $$

        \hl2{Linear connections and the dual connections}
        A connection on a vector fibration $\bE$ is given by the
horizontal distribution and can be represented by a section $\zd$
of the fibration  $\ss j^1(\bE) \rightarrow \ss j^0(\bE) =\bE$.
The connection  is linear if $\zd$ is a morphism of vector
bundles
                $$ \zd \colon \bE \rightarrow j^1(\bE),
                                                                                                                                $$
where $j^1(\bE)$ is a vector fibration over $M$.

        A connection on $\bE$ defines a splitting of the tangent bundle
$\sT E$ into the vertical and horizontal parts. Since the bundle
$\ss V E$ of vertical vectors can be identified with the product
$\bE\times _M E$, we can look at the splitting  map as an
isomorphism of vector bundles
                $$ D\colon  \sT E \rightarrow (\sT M \oplus_M\bE)\times_M
E
                                    \tag \label{Fdb56}$$
        over the identity of $E$.
                \claim{Proposition}{} \label{Cdb16}
        A mapping $D\colon  \sT E \rightarrow (\sT M \oplus_M\bE)
\times_M E$ is  the splitting related to a linear connection if
and only if $D$ defines a double vector bundle morphism
                $$  \bold D \colon \sT \bE \rightarrow \bK(\sT M,\bE,\bE),
                                      \tag\label{Fdb57}       $$
such that the corresponding mappings
                $$ \split D_r \colon & E\rightarrow E\\
                                D_l \colon & \sT M\rightarrow \sT M\\
                                D_c \colon & E\rightarrow E\\
                                \endsplit                               \tag \label{Fdb58}$$
are identities.
                                        \endclaim
\medskip

        Let $D$ be the the splitting of a linear connection on $\bE$.
The transposed left dual to $D$ defines an isomorphism
                $$\bold D^\* \colon \sT \bE^\* \rightarrow \bK(\sT M,\bE^\*, \bE^\*)
                                                                                                                \tag \label{Fdb59}$$
        and, because of \Ref{Fdb58} and \Ref{Fdb45}, $D^\*_r,\,
D^\*_l,\, D^\*_c$ are identities.
        Thus $\bold D^\*$ is the splitting of a linear connection on
$\bE^\*$. We call it {\it the dual connection}.

        Let $g\colon \bE \rightarrow \bE^\*$ be a metric on $\bE$
($g$ is a self-adjoint isomorphism of vector bundles).
        The splitting $D$ is the splitting of a metric connection if
the following diagram is commutative
                $$ \CD \sT E  @()\L{\sT g} @(0,-1)  @()\L{\bold D} @(1,0) &
\bK(\sT M, \bE, \bE) @()\l{\id_{\ssT M}\times  g \times g} @(0,-1) \\
                                \sT E^\*  @()\L{ \bold D^\*} @(1,0) & \bK(\sT M,
\bE^\*, \bE^\*)
                        \endCD  \quad .                                         \tag \label{Fdb60}$$

        \hl2{Symmetric connections}
        In this section $\bE = \sT M$. We have then the canonical
isomorphism
                $$  \ezk_M \colon \sT \sT M \rightarrow \ss J (\sT\sT M).
                                                                                                                                $$
        We introduce also an isomorphism
                $$ \ezk \colon \bK(\sT M,\sT M,\sT M) \rightarrow \ss J(\bK(\sT M,
\sT M, \sT M))
                                                                                                                                $$
        by
                $$  \ezk(v,w,u) = (u,w,v).
                                                                                                                                $$

        \claim{Proposition}{} \label{Cdb17}
        A connection $D$ is  symmetric (torsion-free) if and only if
                $$  \ezk \circ  \bold D = \ss J (\bold D)\circ \ezk_M
                                                                                                                        \tag \label{Fdb61}$$
        i.~e., if the following diagram is commutative
                $$ \CD \sT \sT M  @()\L{\ezk_M} @(0,-1)  @()\L{\bold D} @(1,0) &
\bK(\sT M, \sT M,\sT M) @()\l{\zk} @(0,-1) \\
                                \ss J(\sT \sT M) @()\L{ \ss J(\bold D)} @(1,0) & \ss
J(\bK(\sT M,\sT M,\sT M))
                        \endCD  \quad .                                         \tag \label{Fdb62}$$

                                        \endclaim
        \proof
        Let $(x^i, {\dot x}^j, {x'}^k, {\dot x}'{}^l)$  be an adopted
coordinate system on $\sT \sT M$ and let  $(y^i, {\dot y}_l^j,
{\dot y}_c^k, {\dot y}_r^l )$ be a coordinate system on $\bK(\sT
M, \sT M, \sT M)$. We have then
                $$ \split
                        y^i(\zk\circ D) &=x^i \\
                        {\dot y}_l^j(\zk\circ D) &= {\dot x}^j \\
                        {\dot y}_r^j(\zk\circ D) &= { x'}^j \\
                        {\dot y}_c^l(\zk\circ D) &= {\dot x}'{}^l + \zG^l_{ij}
{\dot x}^i {x'}^j .
                        \endsplit
                                                                                \tag \label{Fdb63}$$
        On the other hand,

                $$ \split
                        y^i( \ss J(D)\circ \zk_M) &=x^i \\
                        {\dot y}_l^j( \ss J(D)\circ \zk_M) &= {\dot x}^j \\
                        {\dot y}_r^j( \ss J(D)\circ \zk_M) &= { x'}^j \\
                        {\dot y}_c^l( \ss J(D)\circ \zk_M) &= {\dot x}'{}^l + \zG^l_{ij}
{\dot x}^j {x'}^i .
                        \endsplit
                                                                                \tag \label{Fdb64}$$
        The diagram \Ref{Fdb62} is commutative if and only if
$\zG^l_{ij} = \zG^l_{ji}$, i.~e., if the connection is symmetric.
\phantom{aa}        \endproof
        In order to obtain conditions for a connection to be
symmetric, in terms of the dual connection, let us consider first a
more general  commutative diagram of isomorphisms of double
vector bundles.

                $$ \CD
                        \bK @() \L{\bold \zF} @(1,0) @() \L{\zk_K} @(0,-1) & \bL @()
\L{\zk_L} @(0,-1) \\
                        \ss J(\bK) @() \L{\ss J(\bold \zF)} @(1,0) & \ss J(\bL)
        \endCD \quad .
                                                                                        \tag \label{Fdb65}$$

        The left dual to this diagram is the commutative diagram
(see \Ref{Fdb21})
                $$ \CD
                        \bK^{\*_l}   & \bL^{\*_l} @()
\L{\bold \zF^{\*_l}} @(-1,0)  \\
                        \ss J(\bK^{\*_r}) @() \L{\zk_K^{\*_l}} @(0,1)  & \ss
J(\bL^{\*_r}) @() \L{\zk_L^{\*_l}} @(0,1) @() \L{\ss J(\bold \zF
^{\*_r})} @(-1,0)
        \endCD \quad .
                                                                                        \tag \label{Fdb66}$$

        Using the canonical isomorphisms of $\bK, \bL$ and
$\bK^{\*_r\*_l} , \bL^{\*_r\*_l} $,  we obtain
                $$ \bK^{\*_r} = \cK(\bK^{\*_r\*_l\*_l
\*_l }) = \cK (\bK^{\*_l \*_l})
                                                                                        \tag \label{Fdb67}$$
         and
                $$ \bL^{\*_r} = \cL(\bL^{\*_r\*_l\*_l
\*_l }) = \cL (\bL^{\*_l \*_l}) .
                                                                                        \tag \label{Fdb68}$$
        With these identifications, we can replace the diagram \Ref{Fdb66}
by an equivalent one:
                $$
                        \cgaps{\w{\ss J(\bold \zF ^{\*_l\*_l})^{-1}}}\CD
                        \bK^{\*_l}   & \bL^{\*_l} @()
\L{\bold \zF^{\*_l}} @(-1,0)  \\
                        \ss J(\bK^{\*_l\*_l}) @() \L{\zb_K^{-1}} @(0,1)  & \ss
J(\bL^{\*_l\*_l}) @() \L{\zb_L^{-1}} @(0,1) @() \L{\ss J(\bold \zF
^{\*_l\*_l})^{-1}} @(-1,0)
        \endCD \quad ,
                                                                                        \tag \label{Fdb69}$$
or
                $$      \cgaps{\w{\ss J(\bold \zF ^{\*_l\*_l})^{-1}}}\CD
                        \bK^{\*_l} @() \L{\zb_K} @(0,-1) @()\L{(\bold
\zF^{\*_l})^{-1}} @(1,0) & \bL^{\*_l} @() \L{\zb_L} @(0,-1)  \\
                        \ss J(\bK^{\*_l\*_l})    & \ss J(\bL^{\*_l\*_l})
@() \L{\ss J(\bold \zF ^{\*_l\*_l})^{-1}} @(-1,0)
        \endCD \quad ,
                                                                                        \tag \label{Fdb70}$$
where
                $$ \zb_K^{-1} = \zk^{\*_l}_K\circ \ss J(\cK) ,\
\zb_L^{-1} = \zk^{\*_l}_L\circ \ss J(\cL) .
                                                                                                                                $$

        In the case of $\bK = \sT \sT M$ and $\bL = \bK(\sT M, \sT
M, \sT M)$, we have, as in \Ref{Fdb41},
                $$ \split
                \bK^{\*_l} = \ss J(\sT\sT ^\* M), &\quad \bL^{\*_l} = \ss
J(\bK(\sT M, \sT^\* M, \sT^\* M)), \\
                \bK^{\*_r} = \ss J(\sT ^\* \sT M), &\quad \bL^{\*_r} = \ss
\bK(\sT M, \sT^\* M, \sT^\* M) \\
                \bK^{\*_l\*_l} = \sT^\*\sT ^\* M &\quad \bL^{\*_l\*_l} = \ss
\bK(\sT M, \sT^\* M, \sT^\* M),
                        \endsplit
                                                                        \tag \label{Fdb71}$$
        and the canonical evaluation between $\bL^{\*_l}$ and
$\bL^{\*_l\*_l}$ is given by the formula
                $$  \langle (v,f,g), (w,h,g)\rangle  = \langle v,h
\rangle  + \langle w,f\rangle .
                                                                        \tag \label{Fdb72}$$
        Moreover, for $\zk_K = \ezk_M$, $\zk_L =\zk$, $\zF =D$, we have
                $$\split
                (\zk_K)^{\*_l} = \ss J(\eza_M^{-1}),  &\quad
(\zk_L)^{\*_l} =\id,\\
        \bold \zF^{\*_l} = (\ss J(\bold D^\*))^{-1}, &\quad (\ss
J(\bold \zF^{\*_l \*_l}))^{-1}= \ss J ((\ss J(\bold
D^\*))^{\*_l}) = (\bold D^\*)^{\*_r}.
                \endsplit
                                                                        \tag \label{Fdb73}$$
        The commutative diagram \Ref{Fdb70} is then  equivalent to
\medskip

                $$ 
                \CD
                \ss J(\sT \sT^\* M) @()\L{\ss J(\bold D^\*)} @(1,0)
@()\L{\zb_K} @(0,-1)  & \ss J(\bK(\sT M, \sT^\* M, \sT^\* M))  @()
\L{\zb_L} @(0,-1)\\
                \ss J(\sT^\* \sT^\* M) & \ss J(\bK(\sT M, \sT^\* M,
\sT^\* M)) @()\L{(\bold D^\*)^{\*_r}} @(-1,0)
                \endCD\quad ,
                                                                                \tag \label{Fdb74}$$
\medskip
and
                $$ \cgaps{\w{\ss J((\bold D^\*)^{\*_r})}}
                \CD
                \sT \sT^\* M @()\L{\bold D^\*} @(1,0)
@()\L{\ss J(\zb_K)} @(0,-1)  & \bK(\sT M, \sT^\* M, \sT^\* M)  @()
\L{\ss J(\zb_L)} @(0,-1)\\
                \sT^\* \sT^\* M & \bK(\sT M, \sT^\* M,
\sT^\* M) @()\L{\ss J((\bold D^\*)^{\*_r})} @(-1,0)
                \endCD\quad .
                                                                                \tag \label{Fdb75}$$
\medskip

        The isomorphism $\cL$ is given by the formula (\Ref{Fdb39})
                $$ \cL \colon (v,f,g) \mapsto  (-v,f,g)
                                                                        $$
        and $\cK \colon \sT^\*\sT^\* M\rightarrow \ss J (\sT^\*\sT
M)$ is a symplectomorphism such that it projects to the identity
on $\sT^\* M$ and the core isomorphism is also the identity on
$\sT ^\* M$. The isomorphism $\eza_M \colon \sT\sT^\* M
\rightarrow \sT^\* \sT M$ is a symplectomorphism between the
tangent canonical symplectic structure on $\sT \sT^\* M$
and the canonical symplectic structure on
$\sT^\*\sT M$ which projects to the identity
on $\sT^\* M$ and the core isomorphism is also the identity on
$\sT ^\* M$ ([{\bf 1}], [{\bf 5}]) .
        It follows that
                $$ \zb = \ss J(\zb_L) \colon (v,f,g) \mapsto (-v,f,g)
                                                                $$
        and that $\zb_K$ is a symplectomorphism such that it
projects to the identity on $\sT^\* M$ and the core isomorphism
is also the identity on $\sT ^\* M$.
        We conclude  that $\zb_K$ is the canonical
symplectic structure on $\sT^\* M$.
        We get the diagram
                $$ \cgaps{\w{\ss J((\bold D^\*)^{\*_r})}}
                \CD
                \sT \sT^\* M @()\L{\bold D^\*} @(1,0)
@()\L{\ezb_M} @(0,-1)  & \bK(\sT M, \sT^\* M, \sT^\* M)  @()
\L{\zb} @(0,-1)\\
                \sT^\* \sT^\* M & \bK(\sT M, \sT^\* M,
\sT^\* M) @()\L{\ss J((\bold D^\*)^{\*_r})} @(-1,0)
                \endCD\quad .
                                                                                \tag \label{Fdb76}$$

$(\bold D^\*)^{-1}  $ is completely determined by its values on
                $$\bK(\sT M, \sT^\* M, \sT ^\* M)\supset V =\{
(v,f,g)\colon  f=0\}
                                                                                                $$
        and, consequently, by $ W = (\bold D^\*)^{-1} (V)  $
        Of course, $(\bold D^\*)^{-1}(V)  $ is the horizontal
distribution of $\bold D^\*$, and $\zb (V) =V$.
        Thus, the diagram \Ref{Fdb76} is commutative if and only if
                $$ \ss J((\bold D^\*)^{\*_r})\circ \zb \circ \bold D^\*
(W) = \ezb_M(W).
                                                                                        \tag \label{Fdb77}$$
        Let $X,Y$ be vector spaces, $T\colon X\rightarrow Y$ a linear
mapping and $Z\subset X$ a vector subspace. We have the equality
                $$ T^\*((F(Z))^\circ ) = Z^\circ .
                                                                                                                                $$
        It follows that, since $V^\circ  = V$,
                $$ \ss J((\bold D^\*)^{\*_r})(V) = W^\circ
                                                                                                \tag \label{Fdb78}$$
         and, consequently, the equality \Ref{Fdb77}
is equivalent to
                $$ W^\circ = \ezb W
                                                                                                                        $$
           We have proved the following theorem.

                \claim{Proposition}{} \label{Cdb18}
        The diagram \Ref{Fdb76} is commutative (the connection is
symmetric) if and only if the horizontal distribution of the
dual connection $\bold D^\*$ is lagrangian.
                                        \endclaim

        \bigskip
        \makebib
\bb \key{\bf 1} \by J. Grabowski and P. Urba\'nski   \pp 6743--6777
        \paper Tangent lifts of Poisson and related structures
        \yr  1995 \vol 28
        \jour J. Phys. A  \endbib

\bb \key{\bf 2} \by K. C. H.  Mackenzie  \pp 180-239
        \paper Double Lie algebroids and second-order geometry I
        \yr 1992  \vol 94
        \jour Adv. Math  \endbib

\bb \key{\bf 3} \by  J. Pradines
        \paper G\'eom\'etrie differentielle au-dessus d'un grupo\"\i de
          \pp 1194-1196\yr 1968 \vol 266
        \jour C. R. Acad. Sci. Paris, s\'erie A
           \endbib

\bb \key{\bf 4} \by  J. Pradines
        \paper Repr\'esentation des jets non holonomes par des
morphismes vectoriels doubles soud\'es
          \pp 1523-1526\yr 1974 \vol 278
        \jour C. R. Acad. Sci. Paris, s\'erie A
           \endbib

\bb \key{\bf 5} \by W. M. Tulczyjew
        \book Geometric Formulations of Physical Theories
        \publ Bibliopolis, Napoli \yr 1989
         \endbib
\bb \key{\bf 6} \by W. M. Tulczyjew \pp 675-678
                \paper Les sous-varie\'et\'es lagrangiennes et la
dynamique lagrangienne
                \yr 1976 \vol 283 \jour C. R. Acad. Sci. Paris, s\'erie A
                \endbib
\bb \key{\bf 7} \by W. M. Tulczyjew \pp 15-18
                \paper Les sous-varie\'et\'es lagrangiennes et la
dynamique hamiltonienne
                \yr 1976 \vol 283 \jour C. R. Acad. Sci. Paris, s\'erie A
                \endbib
\bb \key{\bf 8} \by P. Urba\'nski   \paper Double vector bundles in
classical mechanics \yr 1997 \vol \jour Rend. Sem. Mat. Torino
				\endbib
\bb \key{\bf 9}  \by K. Yano and S. Ishihara
        \book Tangent and Cotangent Bundles
          \publ Marcel Dekker, Inc., New York \yr 1973
        \endbib
        \endmakebib
        \enddocument

%% file: lamstex.tex


\catcode`\@=11
\ifx\amstexloaded@\relax\else
 \errmessage{AmS-TeX must be loaded before LamS-TeX}\fi
\ifx\laxread@\undefined\else\catcode`\@=\active \fi
\def\err@#1{\errmessage{LamS-TeX error: #1}}
\def^^L{\par}
\let\+\tabalign
\def\newcount{\alloc@0\count\countdef\insc@unt}
\def\newdimen{\alloc@1\dimen\dimendef\insc@unt}
\def\newskip{\alloc@2\skip\skipdef\insc@unt}
\def\newmuskip{\alloc@3\muskip\muskipdef\@cclvi}
\def\newbox{\alloc@4\box\chardef\insc@unt}
\let\newtoks\relax
\def\newhelp#1#2{\newtoks#1#1\expandafter{\csname#2\endcsname}}
\def\newtoks{\alloc@5\toks\toksdef\@cclvi}
\def\newread{\alloc@6\read\chardef\sixt@@n}
\def\newwrite{\alloc@7\write\chardef\sixt@@n}
\def\newfam{\alloc@8\fam\chardef\sixt@@n}
\def\newlanguage{\alloc@9\language\chardef\@cclvi}
\def\newinsert#1{\global\advance\insc@unt by\m@ne
  \ch@ck0\insc@unt\count
  \ch@ck1\insc@unt\dimen
  \ch@ck2\insc@unt\skip
  \ch@ck4\insc@unt\box
  \allocationnumber=\insc@unt
  \global\chardef#1=\allocationnumber
  \wlog{\string#1=\string\insert\the\allocationnumber}}
\def\newif#1{\count@\escapechar \escapechar\m@ne
  \expandafter\expandafter\expandafter
   \edef\@if#1{true}{\let\noexpand#1=\noexpand\iftrue}%
  \expandafter\expandafter\expandafter
   \edef\@if#1{false}{\let\noexpand#1=\noexpand\iffalse}%
  \@if#1{false}\escapechar\count@}

\def\Err@#1{\errhelp\defaulthelp@\err@{#1}}
{\catcode`\@=\active
 \edef\next{\gdef\noexpand@{\futurelet\noexpand\next
  \csname at\string@\endcsname}}
 \next
}
\def\at@{\ifcat\noexpand\next a\let\next@\at@@\else
 \ifcat\noexpand\next0\let\next@\at@@\else
 \ifcat\noexpand\next\relax\let\next@\at@@\else
 \let\next@\at@@@\fi\fi\fi\next@}
\def\at@@@{\errhelp\athelp@\err@{Invalid use of @}}
\def\at@@#1{\expandafter
 \ifx\csname\string#1@at\endcsname\relax\let\next@\at@@@\else
 \DN@{\csname\string#1@at\endcsname}\fi\next@}
\def\atdef@#1{\expandafter\def\csname\string#1@at\endcsname}
\newif\iftest@
\def\tagin@#1{\tagin@false
 \DN@##1\tag##2##3\next@{\test@true\ifx\tagin@##2\test@false\fi}%
 \next@#1\tag\tagin@\next@\tagin@false\iftest@\tagin@true\fi}
\let\lkerns@\relax
\def\nolinebreak{\RIfM@\mathmodeerr@\nolinebreak\else
 \ifhmode\saveskip@\lastskip\unskip
 \nobreak\ifdim\saveskip@>\z@\hskip\saveskip@\fi\lkerns@
 \else\vmodeerr@\nolinebreak\fi\fi}
\def\allowlinebreak{\RIfM@\mathmodeerr@\allowlinebreak\else
 \ifhmode\saveskip@\lastskip\unskip
 \allowbreak\ifdim\saveskip@>\z@\hskip\saveskip@\fi\lkerns@
 \else\vmodeerr@\allowlinebreak\fi\fi}
\def\linebreak{\RIfM@\mathmodeerr@\linebreak\else
 \ifhmode\unskip\unkern\break\lkerns@
 \else\vmodeerr@\linebreak\fi\fi}
\let\nkerns@\relax
\def\newline{\RIfM@\mathmodeerr@\newline\else
 \ifhmode\unskip\unkern\null\hfill\break\nkerns@
 \else\vmodeerr@\newline\fi\fi}%
\def\newbox@{\alloc@@4\box\chardef\insc@unt}
\def\newcount@{\alloc@@0\count\countdef\insc@unt}
\def\accentedsymbol#1#2{\expandafter\newbox@\csname\exstring@#1@box\endcsname
 \setbox\csname\exstring@#1@box\endcsname\hbox{$\m@th#2$}%
 \define#1{\copy\csname\exstring@#1@box\endcsname{}}}
\def\rightadd@#1\to#2{\toks@{\\#1}\toks@@\expandafter{#2}\xdef#2{\the\toks@@
 \the\toks@}\toks@{}\toks@@{}}
\def\fontlist@{\\\tenrm\\\sevenrm\\\fiverm\\\teni\\\seveni\\\fivei
 \\\tensy\\\sevensy\\\fivesy\\\tenex\\\tenbf\\\sevenbf\\\fivebf
 \\\tensl\\\tenit}
\def\font@#1=#2 {\rightadd@#1\to\fontlist@\font#1=#2 }
\def\ismember@#1#2{\global\let\Next@ F\let\next@= #2%
 {\def\\##1{\let\nextii@##1\ifx\nextii@\next@\global\let\Next@ T\fi}#1}%
 \test@false\ifx\Next@ T\test@true\fi\let\next@\relax}
\def\FNSS@#1{\let\FNSS@@#1\FN@\FNSS@@@}
\def\FNSS@@@{\ifx\next\space@\def\FNSS@@@@. {\FN@\FNSS@@@}\else
 \def\FNSS@@@@.{\FNSS@@}\fi\FNSS@@@@.}
\atdef@"{\unskip
 \DN@{\ifx\next`\DN@`{\FN@\nextii@}%
  \else\ifx\next\lq\DN@\lq{\FN@\nextii@}%
  \else\DN@####1{\FN@\nextiii@}\fi\fi
  \next@}%
 \DNii@{\ifx\next`\DN@`{\sldl@``}%
  \else\ifx\next\lq\DN@\lq{\sldl@``}%
  \else\DN@{\dlsl@`}\fi\fi\next@}%
 \def\nextiii@{\ifx\next'\DN@'{\srdr@''}%
  \else\ifx\next\rq\DN@\rq{\srdr@''}%
  \else\DN@{\drsr@'}\fi\fi\next@}%
 \FNSS@\next@}
\def\root{%
  \DN@{\ifx\next\uproot\let\next@\nextii@\else
   \ifx\next\leftroot\let\next@\nextiii@\else
   \let\next@\plainroot@\fi\fi\next@}%
  \DNii@\uproot##1{\uproot@##1\relax\FNSS@\nextiv@}%
  \def\nextiv@{\ifx\next\leftroot\let\next@\nextv@\else
   \let\next@\plainroot@\fi\next@}%
  \def\nextv@\leftroot##1{\leftroot@##1\relax\plainroot@}%
  \def\nextiii@\leftroot##1{\leftroot@##1\relax\FNSS@\nextvi@}%
  \def\nextvi@{\ifx\next\uproot\let\next@\nextvii@\else
   \let\next@\plainroot@\fi\next@}%
  \def\nextvii@\uproot##1{\uproot@##1\relax\plainroot@}%
  \bgroup\uproot@\z@\leftroot@\z@
 \FNSS@\next@}
\def\loop#1\repeat{\def\iterate{#1\relax\expandafter\iterate\fi}%
 \iterate\let\iterate\relax}
\def\gloop@#1\repeat{\gdef\iterate@{#1\relax\expandafter\iterate@\fi}%
 \iterate@\global\let\iterate@\relax}
\def\printoptions{\W@{Do you want S(yntax check),
  G(alleys) or P(ages)?^^JType S, G or P, follow by <return>: }\loop
 \read\m@ne to\ans@
 \edef\next@{\def\noexpand\Ans@{\ans@}}\uppercase\expandafter{\next@}%
 \ifx\Ans@\S@\test@true\syntax\else
 \ifx\Ans@\G@\test@true\galleys\else
 \ifx\Ans@\P@\test@true\else
 \test@false\fi\fi\fi
 \iftest@\else\W@{Type S, G or P, follow by <return>: }%
 \repeat}
\expandafter\let\csname A@;\endcsname;
\expandafter\let\csname A@:\endcsname:
\expandafter\let\csname A@?\endcsname?
\expandafter\let\csname A@!\endcsname!
\def\APdef#1{\def\next@{\expandafter\let\csname A@\string#1\endcsname#1}%
 \afterassignment\next@\def#1}
\let\fextra@\,
\def\tdots@{\unskip
 \DN@{$\m@th\mathinner{\ldotp\ldotp\ldotp}\,
   \ifx\next,\,$\else\ifx\next.\,$\else
   \ifx\next;\,$\else
   \expandafter\ifx\csname A@\string;\endcsname\next\fextra@$\else
   \ifx\next:\,$\else
   \expandafter\ifx\csname A@\string:\endcsname\next\fextra@$\else
   \ifx\next?\,$\else
   \expandafter\ifx\csname A@\string?\endcsname\next\fextra@$\else
   \ifx\next!\,$\else
   \expandafter\ifx\csname A@\string!\endcsname\next\fextra@$\else
   $ \fi\fi\fi\fi\fi\fi\fi\fi\fi\fi}%
 \ \FN@\next@}
\def\extrap@#1{%
 \ifx\next,\DN@{#1\,}\else
 \ifx\next;\DN@{#1\,}\else
 \expandafter\ifx\csname A@\string;\endcsname\next\DN@{#1\fextra@}\else
 \ifx\next.\DN@{#1\,}\else\extra@
 \ifextra@\DN@{#1\,}\else
 \let\next@#1\fi\fi\fi\fi\fi\next@}
\def\dotsc{\DN@{\ifx\next;\plainldots@\,\else
 \expandafter\ifx\csname A@\string;\endcsname\next\plainldots@\fextra@\else
 \ifx\next.\plainldots@\,\else\extra@\plainldots@
 \ifextra@\,\fi\fi\fi\fi}%
 \FN@\next@}
\def\keybin@{\keybin@true
 \ifx\next+\else\ifx\next=\else\ifx\next<\else\ifx\next>\else\ifx\next-\else
 \ifx\next*\else\ifx\next:\else
 \expandafter\ifx\csname A@\string;\endcsname\next\else
 \keybin@false\fi\fi\fi\fi\fi\fi\fi\fi}
\def\boldkey#1{\ifcat\noexpand#1A%
  \ifcmmibloaded@{\fam\cmmibfam#1}\else
   \Err@{First bold symbol font not loaded}\fi
 \else
 \let\next=#1%
 \ifx#1!\mathchar"5\bffam@21 \else
 \expandafter\ifx\csname A@\string!\endcsname\next\mathchar"5\bffam@21 \else
 \ifx#1(\mathchar"4\bffam@28 \else\ifx#1)\mathchar"5\bffam@29 \else
 \ifx#1+\mathchar"2\bffam@2B \else\ifx#1:\mathchar"3\bffam@3A \else
 \expandafter\ifx\csname A@\string:\endcsname\next\mathchar"3\bffam@3A \else
 \ifx#1;\mathchar"6\bffam@3B \else
 \expandafter\ifx\csname A@\string;\endcsname\next\mathchar"6\bffam@3B \else
 \ifx#1=\mathchar"3\bffam@3D \else
 \ifx#1?\mathchar"5\bffam@3F \else
 \expandafter\ifx\csname A@\string?\endcsname\next\mathchar"5\bffam@3F \else
 \ifx#1[\mathchar"4\bffam@5B \else
 \ifx#1]\mathchar"5\bffam@5D \else
 \ifx#1,\mathchari@63B \else
 \ifx#1-\mathcharii@200 \else
 \ifx#1.\mathchari@03A \else
 \ifx#1/\mathchari@03D \else
 \ifx#1<\mathchari@33C \else
 \ifx#1>\mathchari@33E \else
 \ifx#1*\mathcharii@203 \else
 \ifx#1|\mathcharii@06A \else
 \ifx#10\bold0\else\ifx#11\bold1\else\ifx#12\bold2\else\ifx#13\bold3\else
 \ifx#14\bold4\else\ifx#15\bold5\else\ifx#16\bold6\else\ifx#17\bold7\else
 \ifx#18\bold8\else\ifx#19\bold9\else
  \Err@{\noexpand\boldkey can't be used with #1}%
 \fi\fi\fi\fi\fi\fi\fi\fi\fi\fi\fi\fi\fi\fi\fi
 \fi\fi\fi\fi\fi\fi\fi\fi\fi\fi\fi\fi\fi\fi\fi\fi\fi\fi}
\def\arabic#1{#1}
\def\alph#1{\count@#1\relax\advance\count@96 \ifnum\count@>122
 \Err@{\noexpand\alph invalid for numbers > 26}\else\char\count@\fi}
\def\Alph#1{\count@#1\relax\advance\count@64 \ifnum\count@>90
 \Err@{\noexpand\Alph invalid for numbers > 26}\else\char\count@\fi}

\def\Roman#1{\uppercase\expandafter{\romannumeral#1}}
\def\fnsymbol#1{\count@#1\relax
 \count@@\count@
 \advance\count@\m@ne\divide\count@7
 \count@@@\count@\advance\count@@@\@ne
 \multiply\count@7 \advance\count@@-\count@
 \count@\count@@@
 {\loop
  \ifcase\count@@\or*\or\dag\or\ddag\or\P\or\S\or\text{$\|$}\or\#\fi
  \advance\count@\m@ne\ifnum\count@>\z@\repeat}}
\def\cardnine@#1{\ifcase#1\or one\or two\or three\or four\or five\or
 six\or seven\or eight\or nine\fi}
\let\alloc@\alloc@@
\newcount\ten@
\ten@10
\def\cardinal#1{\count@#1\relax
 \ifnum\count@>99 \number\count@
 \else
  \ifnum\count@=\z@ zero%
  \else
   \ifnum\count@<\ten@\cardnine@\count@
   \else
    \ifnum\count@<20
     \advance\count@-\ten@
     \ifcase\count@ ten\or eleven\or twelve\or thirteen\or fourteen\or
      fifteen\or sixteen\or seventeen\or eighteen\or nineteen\fi
    \else
     \count@@\count@\count@@@\count@@
     \divide\count@\ten@\multiply\count@\ten@
     \advance\count@@@-\count@\divide\count@\ten@
     \ifcase\count@\or\or twenty\or thirty\or forty\or fifty\or sixty\or
      seventy\or eighty\or ninety\fi
     \ifnum\count@@@=\z@\else-\cardnine@\count@@@\fi
    \fi
   \fi
  \fi
 \fi}
\def\ordnine@#1{\ifcase#1\or first\or second\or third\or fourth\or fifth\or
 sixth\or seventh\or eighth\or ninth\fi}
\newcount\count@@@@
\def\ordsuffix@{\count@@@@\count@
 \divide\count@\ten@
 \count@@@\count@\count@@\count@
 \divide\count@@\ten@\multiply\count@@\ten@
 \advance\count@@@-\count@@
 \ifnum\count@@@=\@ne th%
 \else
  \count@@@\count@@@@
  \count@@\count@@@@
  \divide\count@@\ten@\multiply\count@@\ten@
  \advance\count@@@-\count@@
  \ifcase\count@@@ th\or st\or nd\or rd\else th\fi
 \fi}
\def\nordinal#1{\count@#1\relax\number\count@\ordsuffix@}
\def\spordinal#1{\count@#1\relax\number\count@$^{\text{\ordsuffix@}}$}
\def\ordinal#1{\count@#1\relax
 \ifnum\count@>99 \number\count@\ordsuffix@
 \else
   \ifnum\count@=\z@ zeroth%
  \else
    \ifnum\count@<\ten@\ordnine@\count@
    \else
     \ifnum\count@<20 \advance\count@-\ten@
      \ifcase\count@ tenth\or eleventh\or twelfth\or thirteenth\or
       fourteenth\or fifteenth\or sixteenth\or seventeenth\or eighteenth\or
       nineteenth\fi
     \else
      \count@@\count@
      \divide\count@\ten@\multiply\count@\ten@
      \count@@@\count@@\advance\count@@@-\count@
      \divide\count@\ten@
      \ifcase\count@\or\or twent\or thirt\or fort\or fift\or sixt\or sevent\or
       eight\or ninet\fi
      \ifnum\count@@@=\z@ ieth\else y-\ordnine@\count@@@\fi
     \fi
    \fi
  \fi
 \fi}
\font@\tensmc=cmcsc10
\textonlyfont@\smc\tensmc
\newtoks\noexpandtoks@
\noexpandtoks@{\let\arabic\relax\let\alph\relax\let\Alph\relax
 \let\Roman\relax\let\fnsymbol\relax\let\rm\relax
 \let\it\relax\let\bf\relax\let\sl\relax\let\smc\relax
 \let\/\relax\let\null\relax}
\def\noexpands@{\the\noexpandtoks@}
\def\Nonexpanding#1{\global\noexpandtoks@
 \expandafter{\the\noexpandtoks@\let#1\relax}}
\def\prevanish@{\saveskip@\z@\ifhmode\saveskip@\lastskip\unskip\fi}
\def\postvanish@{\ifdim\saveskip@>\z@\hskip\saveskip@\fi\FN@\postvanish@@}
\def\postvanish@@{\DN@.{}%
 \ifx\next\space@\ifdim\saveskip@>\z@\DN@. {}\fi\fi\next@.}
\def\invisible#1{\prevanish@\ignorespaces#1\unskip\postvanish@}
\def\vanishlist@{\\\invisible}
\let\noindent@\noindent
\def\noindent{\par\noindent@\FN@\pretendspace@}
\def\pretendspace@{\ismember@\vanishlist@\next
 \iftest@\nobreak\hskip-\p@\hskip\p@\fi}

\newtoks\everypartoks@
\def\noindent@@{\par\everypartoks@\expandafter{\the\everypar}\everypar{}%
 \noindent@\everypar\expandafter{\the\everypartoks@}}
\def\page{\Err@{\noexpand\page has no meaning by itself}}
\let\page@C\pageno
\let\page@P\empty
\let\page@Q\empty
\def\page@S#1{#1\/}
\def\page@F{\rm}
\def\page@N{\arabic}   
\newif\ifindexing@
\def\indexfile{\ifindexing@\else
 \alloc@@7\write\chardef\sixt@@n\ndx@
 \immediate\openout\ndx@=\jobname.ndx
 \global\indexing@true\fi}
\global\advance\insc@unt\m@ne
\ch@ck0\insc@unt\count
\ch@ck1\insc@unt\dimen
\ch@ck2\insc@unt\skip
\ch@ck4\insc@unt\box
\allocationnumber\insc@unt
\global\chardef\margin@\allocationnumber
\dimen\margin@\maxdimen
\count\margin@\z@
\skip\margin@\z@
\newif\ifindexproofing@
\def\indexproofing{\indexproofing@true}
\def\noindexproofing{\indexproofing@false}
\def\unmacro@#1:#2->#3\unmacro@{\def\macpar@{#2}\def\macdef@{#3}}
\def\starparts@#1{\def\stari@{#1}\def\starii@{#1}\let\stariii@\empty
 \test@false
 \DN@##1*##2##3\next@{\ifx\starparts@##2\test@false\else\test@true\fi}%
 \next@#1*\starparts@\next@
 \iftest@\DN@{\starparts@@#1\starparts@@}\else\let\next@\relax\fi\next@}
\def\starparts@@#1*#2\starparts@@{\def\starii@{#1}\def\stariii@{*#2}}
\def\windex@{\ifindexing@
 \expandafter\unmacro@\meaning\stari@\unmacro@
 \edef\macdef@{\string"\macdef@\string"}%
 \edef\next@{\write\ndx@{\macdef@}}\next@
 \write\ndx@{{\number\pageno}{\page@N}{\page@P}{\page@Q}}%
 \fi
 \ifindexproofing@
  \ifx\stariii@\empty\else
   \expandafter\unmacro@\meaning\stariii@\unmacro@\fi
  \insert\margin@{\hbox{\rm\vrule\height9\p@\depth2\p@\width\z@\starii@
  \ifx\stariii@\empty\else\tt\macdef@\fi}}\fi}
\catcode`\"=\active
\def"{\FN@\quote@}
\def\quote@{\ifx\next"\expandafter\quote@@\else\expandafter\quote@@@\fi}
\def\quote@@@#1"{\starparts@{#1}\starii@\windex@}
\def\quote@@"#1"{\prevanish@\starparts@{#1}\windex@\FN@\quote@@@@}
\def\quote@@@@{\ifx\next"\DN@"{\postvanish@}\else
 \let\next@\postvanish@\fi\next@}
\rightadd@"\to\vanishlist@
\def\idefine#1{\DN@{#1}\DNii@{\noexpand#1}%
 \afterassignment\idefine@\def\nextiii@}
\def\idefine@{\ifindexing@
 \expandafter\let\next@\nextiii@
 \expandafter\unmacro@\meaning\nextiii@\unmacro@
 \immediate\write\ndx@{\noexpand\define\nextii@\macpar@{\macdef@}}\fi}
\def\iabbrev*#1#2{\ifindexing@\toks@{#2}%
 \immediate\write\ndx@{\noexpand\abbrev*\noexpand#1{\the\toks@}}\fi}
\newread\laxread@
\newwrite\laxwrite@
\let\fnpages@\empty
\def\Finit@#1#2\Finit@{\let\nextii@#1\def\nextiii@{#2}}
\catcode`\~=11
\def\getparts@ @#1~#2~#3~#4~#5~#6{\def\nextiv@{#1}%
 \def\nextiii@{#2~#3~#4~#5~}\count@#6\relax}
\newif\ifdocument@
\def\document{\ifdocument@\else\global\document@true
 \let\fontlist@\empty
 \immediate\openin\laxread@=\jobname.lax\relax
 {\endlinechar\m@ne\noexpands@\catcode`\@=11 \catcode`\~=11
  \loop\ifeof\laxread@\else
   \read\laxread@ to\next@
   \ifx\next@\empty
   \else
    \expandafter\Finit@\next@\Finit@
    \if\nextii@ F%
     \expandafter\rightadd@\nextiii@\to\fnpages@
    \else
     \expandafter\getparts@\next@
     \edef\next@{\gdef\csname\nextiv@ @L\endcsname{\nextiii@\number\count@}}%
     \next@
    \fi
   \fi
  \repeat}%
 \immediate\closein\laxread@
 \immediate\openout\laxwrite@=\jobname.lax\relax\fi}
\let\thelabel@\relax
\def\thelabels@{\thelabel@ ~\thelabel@@ ~\thelabel@@@ ~\thelabel@@@@ ~}
\def\label#1{\prevanish@
 \ifx\thelabel@\relax
  \Err@{There's nothing here to be labelled}%
 \else
  {\noexpands@
  \expandafter\ifx\csname#1@L\endcsname\relax
   \expandafter\xdef\csname#1@L\endcsname{\thelabels@0}%
   \immediate\write\laxwrite@{@#1~\thelabels@1}%
  \else
   \edef\next@{@~\csname#1@L\endcsname}%
    \expandafter\getparts@\next@
    \ifodd\count@
    \expandafter\xdef\csname#1@L\endcsname{\thelabels@0}%
    \immediate\write\laxwrite@{@#1~\thelabels@1}%
   \else
    \Err@{Label #1 already used}%
   \fi
  \fi
  }%
 \fi
 \postvanish@}
\rightadd@\label\to\vanishlist@
\def\thepages@{\page@N{\number\page@C}~%
 \page@S{\page@P\page@N{\number\page@C}\page@Q}~%
 \number\page@C ~\page@P\page@N{\number\page@C}\page@Q ~}
\def\pagelabel#1{\prevanish@
 \expandafter\ifx\csname#1@L\endcsname\relax
  {\noexpands@
  \expandafter\xdef\csname#1@L\endcsname{\thepages@2}}%
  \write\laxwrite@{@#1~\thepages@3}%
 \else
  {\noexpands@
  \edef\next@{@~\csname#1@L\endcsname}%
  \expandafter\getparts@\next@
  \ifodd\count@
   \ifnum\count@=\@ne
    \expandafter\xdef\csname#1@L\endcsname{\thelabels@2}%
   \fi
   \write\laxwrite@{@#1~\thepages@3}%
  \else
   \Err@{Label #1 already used}%
  \fi
  }%
 \fi
 \postvanish@}
\rightadd@\pagelabel\to\vanishlist@
\newif\ifreferr@
\referr@true
\def\RefErrors{\global\referr@true}
\def\RefWarnings{\global\referr@false}
\setbox\z@\hbox{\global\count@=`^^30}
\ifnum\count@=48 \let\versionthree@\relax\fi
\def\nolabel@#1#2#3{\expandafter\ifx\csname#2@L\endcsname\relax
 \ifreferr@\Err@{No \noexpand\label found for #2}\else
 \W@{Warning: No \noexpand\label found for #2.}%
 \ifx\versionthree@\relax\W@{l.\number\inputlineno\space ... \string#1{#2}}\fi
 \fi#3\else}
\def\csL@#1{{\noexpands@\xdef\Next@{\csname#1@L\endcsname}}}
\def\ref#1{\nolabel@\ref{#1}\relax
 \DNii@##1~##2\nextii@{##1}%
 \csL@{#1}\expandafter\nextii@\Next@\nextii@\fi}
\def\Ref#1{\nolabel@\Ref{#1}\relax
 \DNii@##1~##2~##3\nextii@{##2}%
 \csL@{#1}\expandafter\nextii@\Next@\nextii@\fi}
\def\nref#1{\nolabel@\nref{#1}\relax
 \DNii@##1~##2~##3~##4\nextii@{##3}%
 \csL@{#1}\expandafter\nextii@\Next@\nextii@\fi}
\def\pref#1{\nolabel@\pref{#1}\relax
 \DNii@##1~##2~##3~##4~##5\nextii@{##4}%
 \csL@{#1}\expandafter\nextii@\Next@\nextii@\fi}
\let\pref@\pref
\def\Evaluatenref#1{\nolabel@\Evaluatenref{#1}{\gdef\Nref{-10000 }}%
 \DNii@##1~##2~##3~##4\nextii@{\DNii@{##3}}%
 \csL@{#1}\expandafter\nextii@\Next@\nextii@
 \xdef\Nref{\nextii@}\fi}
\def\Evaluatepref#1{\nolabel@\Evaluatepref{#1}{\global\let\Pref\empty}%
 \DNii@##1~##2~##3~##4~##5\nextii@{\DNii@{##4}}%
 \csL@{#1}\expandafter\nextii@\Next@\nextii@
 \xdef\Pref{\nextii@}\fi}
\def\readlax#1{\immediate\openin\laxread@=#1.lax\relax
 \ifeof\laxread@\W@{}\W@{File #1.lax not found.}\W@{}\fi
 {\endlinechar\m@ne\noexpands@\catcode`\@=11 \catcode`\~=11
  \loop\ifeof\laxread@\else
   \read\laxread@ to\nextv@
   \ifx\nextv@\empty
   \else
    \expandafter\Finit@\nextv@\Finit@
    \ifx\nextii@ F%
    \else
     \expandafter\getparts@\nextv@
     \expandafter\ifx\csname\nextiv@ @L\endcsname\relax
      \edef\next@{\gdef\csname\nextiv@ @L\endcsname
       {\nextiii@\ifnum\count@=\@ne0\else2\fi}}%
      \next@
     \else
      \Err@{Label \nextiv@\space in #1.lax already used}%
     \fi
    \fi
   \fi
  \repeat}%
 \immediate\closein\laxread@}
\catcode`\~=\active

\def\FNSSP@{\FNSS@\pretendspace@}
\everydisplay{\csname displaymath \endcsname}
\expandafter\def\csname displaymath \endcsname#1$${#1$$\FNSSP@}
\def\locallabel@{\let\thelabel@\Thelabel@\let\thelabel@@\Thelabel@@
 \let\thelabel@@@\Thelabel@@@\let\thelabel@@@@\Thelabel@@@@}
\newcount\tag@C
\tag@C\z@
\let\tag@P\empty
\let\tag@Q\empty
\def\tag@S#1{{\rm(}{#1\/}{\rm)}}
\let\tag@N\arabic
\def\tag@F{\rm}
\def\maketag@{\FN@\maketag@@}
\def\maketag@@{\ifx\next\relax\DN@\relax{\FN@\maketag@@}\else
 \ifx\next"\let\next@\maketag@@@\else
 \let\next@\maketag@@@@\fi\fi\next@}
\def\xdefThelabel@#1{\xdef\Thelabel@{#1{\Thelabel@@@}}}
\def\xdefThelabel@@#1{\xdef\Thelabel@@{#1{\Thelabel@@@@}}}
\def\maketag@@@@#1\maketag@{\global\advance\tag@C\@ne
 {\noexpands@
  \xdef\Thelabel@@@{\number\tag@C}%
  \xdefThelabel@\tag@N
  \xdef\Thelabel@@@@{\ifmathtags@$\tag@P\Thelabel@\tag@Q$\else
   \tag@P\Thelabel@\tag@Q\fi}%
  \xdefThelabel@@\tag@S
  }%
 \locallabel@
 \hbox{\tag@F\thelabel@@}%
 #1}
\def\Qlabel@#1{{\noexpands@\xdef\Thelabel@@{#1}%
 \let\style\empty\xdef\Thelabel@@@@{#1}%
 \let\pre\empty\let\post\empty\xdef\Thelabel@{#1}%
 \let\numstyle\empty\xdef\Thelabel@@@{#1}}}
\def\maketag@@@"#1"#2\maketag@{%
 {\let\pre\tag@P\let\post\tag@Q\let\style\tag@S\let\numstyle\tag@N
  \hbox{\tag@F#1}%
  \noexpands@
  \Qlabel@{#1}%
  }%
 \locallabel@
 #2}
\def\align@{\inalign@true\inany@true
 \vspace@\allowdisplaybreak@\displaybreak@\intertext@
 \def\tag{\global\tag@true\ifnum\and@=\z@
  \DN@{&\omit\global\rwidth@\z@&\relax}\else
  \DN@{&\relax}\fi\next@}%
 \iftagsleft@\DN@{\csname align \endcsname}\else
  \DN@{\csname align \space\endcsname}\fi\next@}
\def\noset@{\def\Offset##1##2{\prevanish@\postvanish@}%
 \def\Reset##1##2{\prevanish@\postvanish@}}
\def\measure@#1\endalign{\global\lwidth@\z@\global\rwidth@\z@
 \global\maxlwidth@\z@\global\maxrwidth@\z@
 \global\and@\z@
 \setbox\z@\vbox
  {\noset@\everycr{\noalign{\global\tag@false\global\and@\z@}}\Let@
  \halign{\setboxz@h{$\m@th\displaystyle{\@lign##}$}%
   \global\lwidth@\wdz@
   \ifdim\lwidth@>\maxlwidth@\global\maxlwidth@\lwidth@\fi
   \global\advance\and@\@ne
   &\setboxz@h{$\m@th\displaystyle{{}\@lign##}$}\global\rwidth@\wdz@
   \ifdim\rwidth@>\maxrwidth@\global\maxrwidth@\rwidth@\fi
   \global\advance\and@\@ne
   &\Tag@\eat@{##}\crcr#1\crcr}}%
 \totwidth@\maxlwidth@\advance\totwidth@\maxrwidth@}
\def\prepost@{\global\let\tag@P@\tag@P\global\let\tag@Q@\tag@Q}
\def\reprepost@{\let\tag@P\tag@P@\let\tag@Q\tag@Q@}
\expandafter\def\csname align \space\endcsname#1\endalign
 {\measure@#1\endalign\global\and@\z@
 \ifingather@\everycr{\noalign{\global\and@\z@}}\else\displ@y@\fi
 \Let@\tabskip\centering@
 \halign to\displaywidth
  {\hfil\strut@\setboxz@h{$\m@th\displaystyle{\@lign##\prepost@}$}%
  \boxz@\global\advance\and@\@ne
  \tabskip\z@skip
  &\setboxz@h{$\m@th\displaystyle{{}\@lign##\prepost@}$}%
  \global\rwidth@\wdz@\boxz@\hfil\global\advance\and@\@ne
  \tabskip\centering@
  &\setboxz@h{\@lign\strut@\reprepost@\maketag@##\maketag@}%
  \dimen@\displaywidth\advance\dimen@-\totwidth@
  \divide\dimen@\tw@\advance\dimen@\maxrwidth@\advance\dimen@-\rwidth@
  \ifdim\dimen@<\tw@\wdz@\llap{\vtop{\normalbaselines\null\boxz@}}%
  \else\llap{\boxz@}\fi
  \tabskip\z@skip
  \crcr#1\crcr
  \black@\totwidth@}}
\expandafter\def\csname align \endcsname#1\endalign{\measure@#1\endalign
 \global\and@\z@
 \ifdim\totwidth@>\displaywidth\let\displaywidth@\totwidth@\else
  \let\displaywidth@\displaywidth\fi
 \ifingather@\everycr{\noalign{\global\and@\z@}}\else\displ@y@\fi
 \Let@\tabskip\centering@\halign to\displaywidth
  {\hfil\strut@\setboxz@h{$\m@th\displaystyle{\@lign##\prepost@}$}%
  \global\lwidth@\wdz@\global\lineht@\ht\z@
  \boxz@\global\advance\and@\@ne
  \tabskip\z@skip&\setboxz@h{$\m@th\displaystyle{{}\@lign##\prepost@}$}%
  \ifdim\ht\z@>\lineht@\global\lineht@\ht\z@\fi
  \boxz@\hfil\global\advance\and@\@ne
  \tabskip\centering@&\kern-\displaywidth@
  \setboxz@h{\@lign\strut@\reprepost@\maketag@##\maketag@}%
  \dimen@\displaywidth\advance\dimen@-\totwidth@
  \divide\dimen@\tw@\advance\dimen@\maxlwidth@\advance\dimen@-\lwidth@
  \ifdim\dimen@<\tw@\wdz@
   \rlap{\vbox{\normalbaselines\boxz@\vbox to\lineht@{}}}\else
   \rlap{\boxz@}\fi
  \tabskip\displaywidth@\crcr#1\crcr\black@\totwidth@}}
\def\attag@#1{\let\Maketag@\maketag@\let\TAG@\Tag@
 \let\Prepost@\prepost@\let\Reprepost@\reprepost@
 \let\Tag@\relax\let\maketag@\relax
 \let\prepost@\relax\let\reprepost@\relax
 \ifmeasuring@
  \def\llap@##1{\setboxz@h{##1}\hbox to\tw@\wdz@{}}%
  \def\rlap@##1{\setboxz@h{##1}\hbox to\tw@\wdz@{}}%
 \else\let\llap@\llap\let\rlap@\rlap\fi
 \toks@{\hfil\strut@
  $\m@th\displaystyle{\@lign\the\hashtoks@\prepost@}$%
  \tabskip\z@skip\global\advance\and@\@ne&
  $\m@th\displaystyle{{}\@lign\the\hashtoks@\prepost@}$\hfil
  \ifxat@\tabskip\centering@\fi\global\advance\and@\@ne}%
 \iftagsleft@
  \toks@@{\tabskip\centering@&\Tag@\kern-\displaywidth
   \rlap@{\@lign\reprepost@\maketag@\the\hashtoks@\maketag@}%
   \global\advance\and@\@ne\tabskip\displaywidth}\else
  \toks@@{\tabskip\centering@&\Tag@\llap@{\@lign\reprepost@\maketag@
   \the\hashtoks@\maketag@}\global\advance\and@\@ne\tabskip\z@skip}\fi
 \atcount@#1\relax\advance\atcount@\m@ne
 \loop\ifnum\atcount@>\z@
  \toks@\expandafter{\the\toks@&\hfil$\m@th\displaystyle{\@lign
  \the\hashtoks@\prepost@}$\global\advance\and@\@ne
  \tabskip\z@skip
  &$\m@th\displaystyle{{}\@lign\the\hashtoks@\prepost@}$\hfil\ifxat@
  \tabskip\centering@\fi\global\advance\and@\@ne}\advance\atcount@\m@ne
 \repeat
 \edef\preamble@{\the\toks@\the\toks@@}%
 \edef\preamble@@{\preamble@}%
 \let\maketag@\Maketag@\let\Tag@\TAG@
 \let\prepost@\Prepost@\let\reprepost@\Reprepost@}
\def\unlabel@{\def\label##1{\prevanish@\postvanish@}%
 \def\pagelabel##1{\prevanish@\postvanish@}}
\newcount\tag@CC
\expandafter\def\csname alignat \endcsname#1#2\endalignat
 {\inany@true\xat@false
 \def\tag{\global\tag@true
  \count@#1\relax\multiply\count@\tw@\advance\count@\m@ne
  \gdef\tag@{&}%
  \loop\ifnum\count@>\and@\xdef\tag@{&\omit\tag@}%
  \advance\count@\m@ne\repeat
  \tag@\relax}%
 \vspace@\allowdisplaybreak@\displaybreak@\intertext@
 \displ@y@\measuring@true\tag@CC\tag@C
 \setbox\savealignat@\hbox{\noset@\unlabel@$\m@th\displaystyle\Let@
  \attag@{#1}\vbox{\halign{\span\preamble@@\crcr#2\crcr}}$}%
 \measuring@false
 \Let@\attag@{#1}\tag@C\tag@CC
 \tabskip\centering@\halign to\displaywidth
  {\span\preamble@@\crcr#2\crcr\black@{\wd\savealignat@}}}
\expandafter\def\csname xalignat \endcsname#1#2\endxalignat
 {\inany@true\xat@true
 \def\tag{\global\tag@true
  \count@#1\relax\multiply\count@\tw@\advance\count@\m@ne
  \gdef\tag@{&}%
  \loop\ifnum\count@>\and@\xdef\tag@{&\omit\tag@}%
  \advance\count@\m@ne\repeat
  \tag@\relax}%
 \vspace@\allowdisplaybreak@\displaybreak@\intertext@
 \displ@y@\measuring@true\tag@CC\tag@C
 \setbox\savealignat@\hbox{\noset@\unlabel@$\m@th\displaystyle\Let@
  \attag@{#1}\vbox{\halign{\span\preamble@@\crcr#2\crcr}}$}%
 \measuring@false\Let@\attag@{#1}\tag@C\tag@CC
 \tabskip\centering@\halign to\displaywidth
 {\span\preamble@@\crcr#2\crcr\black@{\wd\savealignat@}}}
\def\gather{\RIfMIfI@\DN@{\onlydmatherr@\gather}\else
 \ingather@true\inany@true\def\tag{&\relax}%
 \vspace@\allowdisplaybreak@\displaybreak@\intertext@
 \displ@y\Let@
 \iftagsleft@\DN@{\csname gather \endcsname}\else
  \DN@{\csname gather \space\endcsname}\fi\fi
 \else\DN@{\onlydmatherr@\gather}\fi\next@}
\def\exstring@{\expandafter\eat@\string}
\def\newcounter#1{\define#1{}%
 \edef\next@{\def\noexpand#1{\futurelet\noexpand\next
  \csname\exstring@#1@Z\endcsname}}\next@
 \edef\next@{\def\csname\exstring@#1@Z\endcsname
  {\global\advance\csname\exstring@#1@C\endcsname\@ne
  {\csname\exstring@#1@F\endcsname\csname\exstring@#1@S\endcsname
   {\csname\exstring@#1@P\endcsname\csname\exstring@#1@N\endcsname
   {\noexpand\number\csname\exstring@#1@C\endcsname}%
   \csname\exstring@#1@Q\endcsname}}%
  \noexpand\ifx\noexpand\next\noexpand\label
   \def\noexpand\next@\noexpand\label########1{{\noexpand\noexpands@
    \xdef\noexpand\Thelabel@{\csname\exstring@#1@N\endcsname
     {\noexpand\number\csname\exstring@#1@C\endcsname}}%
    \xdef\noexpand\Thelabel@@@{\noexpand\number
     \csname\exstring@#1@C\endcsname}%
    \xdef\noexpand\Thelabel@@{\csname\exstring@#1@S\endcsname
     {\csname\exstring@#1@P\endcsname
     \csname\exstring@#1@N\endcsname
     {\noexpand\number\csname\exstring@#1@C\endcsname}%
     \csname\exstring@#1@Q\endcsname}}%
    \xdef\noexpand\Thelabel@@@@{\csname\exstring@#1@P\endcsname
     \csname\exstring@#1@N\endcsname
     {\noexpand\number\csname\exstring@#1@C\endcsname}%
     \csname\exstring@#1@Q\endcsname}}%
    {\noexpand\locallabel@\noexpand\label{########1}}}%
   \noexpand\else\let\noexpand\next@\relax\noexpand\fi\noexpand\next@}}\next@
 \expandafter\newcount@\csname\exstring@#1@C\endcsname
 \expandafter\let\csname\exstring@#1@N\endcsname\arabic
 \expandafter\def\csname\exstring@#1@S\endcsname##1{##1\/}%
 \expandafter\let\csname\exstring@#1@P\endcsname\empty
 \expandafter\let\csname\exstring@#1@Q\endcsname\empty
 \expandafter\def\csname\exstring@#1@F\endcsname{\rm}%
 }
\def\HASH@#1#2{\ifnum#2=\z@\else
 \edef\next@{\toks@{\the\toks@\the\hashtoks@#2}%
 \toks@@{\the\toks@@{\the\hashtoks@#2}}}\next@\expandafter\HASH@\fi}
\def\HASH@@{\toks@{}\toks@@{}\expandafter\HASH@\macpar@00}
\def\usecounter#1#2{\expandafter\ifx\csname\exstring@#1@Z\endcsname
 \relax\Err@{\noexpand#1not created with \string\newcounter}\fi
 \expandafter\let\csname\exstring@#1@@Z\endcsname\relax
 \expandafter\let\csname\exstring@#1@@Z@\endcsname\relax
 \expandafter\let\csname\exstring@#1@@Z@@\endcsname\relax
 \edef\next@{\def\noexpand#2{\futurelet\noexpand\next
  \csname\exstring@#1@@Z\endcsname}}\next@
 \edef\next@{\def\csname\exstring@#1@@Z\endcsname{\noexpand\ifx
  \noexpand\next\noexpand\label\def\noexpand\next@\noexpand\label
   ########1{\csname\exstring@#1@@Z@\endcsname
   {\noexpand#1\noexpand\label{########1}}}%
   \noexpand\else\noexpand\ifx\noexpand\next
   \noexpand"\def\noexpand\next@\noexpand"########1\noexpand"%
   {\csname\exstring@#1@@Z@\endcsname{{\expandafter\noexpand
   \csname\exstring@#1@F\endcsname
   \let\noexpand\pre\expandafter\noexpand\csname\exstring@#1@P\endcsname
   \let\noexpand\post\expandafter\noexpand\csname\exstring@#1@Q\endcsname
   \let\noexpand\style\expandafter\noexpand\csname\exstring@#1@S\endcsname
   \let\noexpand\numstyle\expandafter\noexpand\csname\exstring@#1@N\endcsname
   ########1}}}\noexpand\else
   \def\noexpand\next@{\csname\exstring@#1@@Z@\endcsname{\noexpand#1}}%
   \noexpand\fi\noexpand\fi\noexpand\next@}}\next@
 \def\next@{\expandafter\expandafter\expandafter\unmacro@\expandafter
  \meaning\csname\exstring@#1@@Z@@\endcsname\unmacro@
  \HASH@@
  \edef\next@{\def\csname\exstring@#1@@Z@\endcsname\the\toks@{%
   \expandafter\noexpand\csname\exstring@#1@@Z@@\endcsname\the\toks@@
   \noexpand\FNSSP@}}\next@}%
 \afterassignment\next@
 \expandafter\def\csname\exstring@#1@@Z@@\endcsname}
\def\listbi@{\penalty50 \medskip}
\def\listbii@{\penalty100 \smallskip}
\let\listbiii@\relax
\let\listbiv@\relax
\let\listbv@\relax
\def\listmi@{\advance\leftskip30\p@\relax}
\let\listmii@\listmi@
\let\listmiii@\listmi@
\let\listmiv@\listmi@
\let\listmv@\listmi@
\def\itemi@#1{\noindent@@\llap{#1\hskip5\p@}}
\let\itemii@\itemi@
\let\itemiii@\itemi@
\let\itemiv@\itemi@
\let\itemv@\itemi@
\def\liste@{\penalty-50 \medskip}
\def\listei@{\penalty-100 \smallskip}
\let\listeii@\relax
\let\listeiii@\relax
\let\listeiv@\relax
\expandafter\newcount\csname list@C1\endcsname
\csname list@C1\endcsname\z@
\expandafter\newcount\csname list@C2\endcsname
\csname list@C2\endcsname\z@
\expandafter\newcount\csname list@C3\endcsname
\csname list@C3\endcsname\z@
\expandafter\newcount\csname list@C4\endcsname
\csname list@C4\endcsname\z@
\expandafter\newcount\csname list@C5\endcsname
\csname list@C5\endcsname\z@
\expandafter\let\csname list@P1\endcsname\empty
\expandafter\let\csname list@P2\endcsname\empty
\expandafter\let\csname list@P3\endcsname\empty
\expandafter\let\csname list@P4\endcsname\empty
\expandafter\let\csname list@P5\endcsname\empty
\expandafter\let\csname list@Q1\endcsname\empty
\expandafter\let\csname list@Q2\endcsname\empty
\expandafter\let\csname list@Q3\endcsname\empty
\expandafter\let\csname list@Q4\endcsname\empty
\expandafter\let\csname list@Q5\endcsname\empty
\expandafter\def\csname list@S1\endcsname#1{{\rm(}{#1\/}{\rm)}}
\expandafter\def\csname list@S2\endcsname#1{{\rm(}{#1\/}{\rm)}}
\expandafter\def\csname list@S3\endcsname#1{{\rm(}{#1\/}{\rm)}}
\expandafter\def\csname list@S4\endcsname#1{{\rm(}{#1\/}{\rm)}}
\expandafter\def\csname list@S5\endcsname#1{{\rm(}{#1\/}{\rm)}}
\expandafter\let\csname list@N1\endcsname\arabic
\expandafter\let\csname list@N2\endcsname\arabic
\expandafter\let\csname list@N3\endcsname\arabic
\expandafter\let\csname list@N4\endcsname\arabic
\expandafter\let\csname list@N5\endcsname\arabic
\expandafter\def\csname list@F1\endcsname{\rm}
\expandafter\def\csname list@F2\endcsname{\rm}
\expandafter\def\csname list@F3\endcsname{\rm}
\expandafter\def\csname list@F4\endcsname{\rm}
\expandafter\def\csname list@F5\endcsname{\rm}
\newcount\listlevel@
\listlevel@\z@
\def\list@@C{\csname list@C\number\listlevel@\endcsname}
\def\list@@P{\csname list@P\number\listlevel@\endcsname}
\def\list@@Q{\csname list@Q\number\listlevel@\endcsname}
\def\list@@S{\csname list@S\number\listlevel@\endcsname}
\def\list@@N{\csname list@N\number\listlevel@\endcsname}
\def\list@@F{\csname list@F\number\listlevel@\endcsname}
\newif\iffirstitemi@
\newif\iffirstitemii@
\newif\iffirstitemiii@
\newif\iffirstitemiv@
\newif\iffirstitemv@
\def\Firstitem@true{\csname firstitem\romannumeral\listlevel@
 @true\endcsname}
\def\Firstitem@false{\csname firstitem\romannumeral\listlevel@
 @false\endcsname}
\def\Listm@{\csname listm\romannumeral\listlevel@ @\endcsname}
\def\Item@{\csname item\romannumeral\listlevel@ @\endcsname}
\def\Liste@{\csname liste\romannumeral\listlevel@ @\endcsname}
\newif\iflistcontinue@
\def\keepitem{\listcontinue@true}
\newcount\list@C@
\def\list{%
 \iflistcontinue@\csname list@C1\endcsname\csname list@C@\endcsname\fi
 \global\csname list@C2\endcsname\z@
 \global\csname list@C3\endcsname\z@
 \global\csname list@C4\endcsname\z@
 \global\csname list@C5\endcsname\z@
 \begingroup
 \firstitemi@true
 \listlevel@\@ne
 \def\item{\FN@\item@}%
 \FN@\list@}
\Invalid@\runinitem
\def\list@{\ifx\next\par
 \DN@\par{\FN@\list@}\else
 \ifx\next\runinitem
  \DN@\runinitem{\FN@\runinitem@}\else
  \DN@{\par\dimen@\parskip\parskip\dimen@}\fi\fi\next@}
\newif\ifoutlevel@
\newif\ifrunin@
\def\item@{%
 \ifoutlevel@\Liste@\outlevel@false\fi
 \ifrunin@\runin@false\par
  \dimen@\parskip\parskip\dimen@
  \Listm@\fi
 \iffirstitemi@\listbi@\listmi@\firstitemi@false\else\par\fi
 \iffirstitemii@\listbii@\listmii@\firstitemii@false\else\par\fi
 \iffirstitemiii@\listbiii@\listmiii@\firstitemiii@false\else\par\fi
 \iffirstitemiv@\listbiv@\listmiv@\firstitemiv@false\else\par\fi
 \iffirstitemv@\listbv@\listmv@\firstitemv@false\else\par\fi
 \DN@"##1"{{\let\pre\list@@P\let\post\list@@Q
  \let\style\list@@S\let\numstyle\list@@N
  \vskip-\parskip
  \Item@{\list@@F##1}%
  \noexpands@
  \Qlabel@{##1}}%
  \locallabel@
  \FNSSP@}%
 \DNii@{\global\advance\list@@C\@ne
  {\noexpands@
   \xdef\Thelabel@@@{\number\list@@C}%
   \xdefThelabel@\list@@N
   \xdef\Thelabel@@@@{\list@@P\Thelabel@\list@@Q}%
   \xdefThelabel@@\list@@S
  }%
  \locallabel@
  \vskip-\parskip
  \Item@{\list@@F\thelabel@@}%
  \FN@\pretendspace@}%
 \ifx\next"\expandafter\next@\else\expandafter\nextii@\fi}
\def\runinitem@{%
  \runin@true
  \Firstitem@false
  \DN@"##1"{{\let\pre\list@@P\let\post\list@@Q
   \let\style\list@@S\let\numstyle\list@@N
   \unskip\space{\list@@F##1} %
   \noexpands@
   \Qlabel@{##1}}%
   \locallabel@
   \ignorespaces}%
  \DNii@{\global\advance\list@@C\@ne
   {\noexpands@
    \xdef\Thelabel@@@{\number\list@@C}%
    \xdefThelabel@\list@@N
    \xdef\Thelabel@@@@{\list@@P\Thelabel@\list@@Q}%
    \xdefThelabel@@\list@@S
   }%
   \locallabel@
   \unskip\space{\list@@F\thelabel@@} }%
  \ifx\next"\expandafter\next@\else\expandafter\nextii@\fi}
\def\inlevel{\ifnum\listlevel@=5
 \DN@{\Err@{Already 5 levels down}}\else
 \DN@{\begingroup\advance\listlevel@\@ne
 \Firstitem@true\FN@\inlevel@}\fi\next@}
\def\inlevel@{\ifx\next\par
 \DN@\par{\FN@\inlevel@}\else
 \ifx\next\runinitem
  \DN@\runinitem{\FN@\runinitem@}\else
  \let\next@\relax\fi\fi\next@}
\def\outlevel{\ifnum\listlevel@=\@ne
 \Err@{At top level}\else
 \par\global\list@@C\z@\endgroup\outlevel@true\fi}
\def\endlist{%
 \expandafter\global\csname list@C@\endcsname\csname list@C1\endcsname
 \par
 \global\toks\@ne{}\count@\listlevel@
 {\loop
  \ifnum\count@>\z@\global\toks\@ne\expandafter{\the\toks\@ne\endgroup}%
  \advance\count@\m@ne
  \repeat}%
 \the\toks\@ne
 \liste@
 \listcontinue@false\global\csname list@C1\endcsname\z@
 \vskip-\parskip
 \noindent@@
 \FN@\pretendspace@}
\newif\iffirstdescribe@
\def\describe{\par
 \begingroup\firstdescribe@true
 \def\item##1{%
  \iffirstdescribe@\penalty50 \medskip\vskip-\parskip
  \firstdescribe@false\else\par\fi
  \noindent@@\hangindent2pc\hangafter\@ne
  {\bf##1}\hskip.5em}}

\Invalid@\pullin
\Invalid@\pullinmore
\newif\iffirstpull@
\def\margins{\par\begingroup\firstpull@true
 \def\pullin##1##2{\par
  \iffirstpull@\firstpull@false\else\endgroup\fi
  \begingroup\DN@{##1}%
  \ifx\next@\empty\leftskip\z@\else\ifx\next@\space\leftskip\z@
  \else\leftskip##1\fi\fi
  \DN@{##2}\ifx\next@\empty\rightskip\z@\else\ifx\next@\space
  \rightskip\z@\else\rightskip##2\fi\fi\ignorespaces}%
 \def\pullinmore##1##2{\par
  \xdef\Next@{\leftskip\the\leftskip\relax\rightskip\the\rightskip\relax}%
  \iffirstpull@\firstpull@false\else\endgroup\fi
  \begingroup\Next@
  \DN@{##1}%
  \ifx\next@\empty\else\ifx\next@\space\else\advance\leftskip##1\fi\fi
  \DN@{##2}\ifx\next@\empty\else\ifx\next@\space\else
  \advance\rightskip##2\fi\fi\ignorespaces}}

\newif\ifnopunct@
\newif\ifnospace@
\newif\ifoverlong@
\let\nofrillslist@\empty
\let\overlonglist@\empty
\def\nopunct{\nopunct@true\FN@\nopunct@}
\def\nospace{\nospace@true\FN@\nospace@}
\def\overlong{\overlong@true\FN@\overlong@}
\def\nopunct@{\ifx\next\nospace
 \DN@\nospace{\nospace@true\FN@\nopnos@}\else\ifx\next\overlong
 \DN@\overlong{\overlong@true\FN@\nopol@}\else
 \let\next@\nopunct@@\fi\fi\next@}
\def\nopunct@@#1{\ismember@\nofrillslist@#1%
 \iftest@\let\next@#1\else
 \DN@{\nopunct@false\Err@{\noexpand\nopunct can't be used with
 \string#1}#1}\fi\next@}
\def\nospace@{\ifx\next\nopunct
 \DN@\nopunct{\nopunct@true\FN@\nopnos@}\else\ifx\next\overlong
 \DN@\overlong{\overlong@true\FN@\nosol@}\else
 \let\next@\nospace@@\fi\fi\next@}
\def\nospace@@#1{\ismember@\nofrillslist@#1%
 \iftest@\let\next@#1\else
 \DN@{\nospace@false\Err@{\noexpand\nospace can't be used with
 \string#1}#1}\fi\next@}
\def\overlong@{\ifx\next\nopunct
 \DN@\nopunct{\nopunct@true\FN@\nopol@}\else\ifx\next\nospace
 \DN@\nospace{\nospace@true\FN@\nosol@}\else
 \let\next@\overlong@@\fi\fi\next@}
\def\overlong@@#1{\ismember@\overlonglist@#1%
 \iftest@\let\next@#1\else
 \DN@{\overlong@false\Err@{\noexpand\overlong can't be used with
 \string#1}#1}\fi\next@}
\def\nopnos@{\ifx\next\overlong
 \DN@\overlong{\overlong@true\nopnosol@}\else
 \let\next@\nopnos@@\fi\next@}
\def\nopol@{\ifx\next\nospace
 \DN@\nospace{\nospace@true\nopnosol@}\else
 \let\next@\nopol@@\fi\next@}
\def\nosol@{\ifx\next\nopunct
 \DN@\nopunct{\nopunct@true\nopnosol@}\else
 \let\next@\nosol@@\fi\next@}
\def\nopnos@@#1{\ismember@\nofrillslist@#1%
 \iftest@\let\next@#1\else
 \DN@{\nopunct@false\nospace@false
  \Err@{\noexpand\nopunct\noexpand\nospace
   can't be used with \string#1}#1}\fi\next@}
\def\testii@#1{\ismember@\nofrillslist@#1%
 \iftest@\let\nextiii@ T\else\let\nextiii@ F\fi
 \ismember@\overlonglist@#1%
 \iftest@\let\nextiv@ T\else\let\nextiv@ F\fi
 \test@false\if\nextiii@ T\if\nextiv@ T\test@true\fi\fi}
\def\nopol@@#1{\testii@{#1}%
 \iftest@\let\next@#1%
 \else\DN@{\if\nextiii@ T\else\nopunct@false\fi
  \if\nextiv@ T\else\overlong@false\fi
  \Err@{\if\nextiii@ T\else\noexpand\nopunct\fi
  \if\nextiv@ T\else\noexpand\overlong\fi can't be used
  with \string#1}#1}\fi\next@}
\def\nosol@@#1{\testii@{#1}%
 \iftest@\let\next@#1%
 \else\DN@{\if\nextiii@ T\else\nospace@false\fi
  \if\nextiv@ T\else\overlong@false\fi
  \Err@{\if\nextiii@ T\else\noexpand\nospace\fi
  \if\nextiv@ T\else\noexpand\overlong\fi can't be used
  with \string#1}#1}\fi\next@}
\def\nopnosol@#1{\testii@{#1}%
 \iftest@\let\next@#1%
 \else\DN@{\if\nextiii@ T\else\nopunct@false\nospace@false\fi
  \if\nextiv@ T\else\overlong@false\fi
  \Err@{\if\nextiii@ T\else\noexpand\nopunct\noexpand\nospace\fi
  \if\nextiv@ T\else\noexpand\overlong\fi can't be used
  with \string#1}#1}\fi\next@}
\def\punct@#1{\ifnopunct@\else#1\fi}
\def\addspace@#1{\ifnospace@\else#1\fi}
\def\hss@{\ifoverlong@\z@ plus\@m\p@ minus\@m\p@
 \else \z@ plus\@m\p@\fi}
\rightadd@\demo\to\nofrillslist@
\newif\ifclaim@
\def\exxx@{\expandafter\expandafter\expandafter\eat@\expandafter\string}
\let\colon@:
\def\demo#1{\ifclaim@
 \Err@{Previous \expandafter\noexpand\claimtype@ has
  no matching \string\end\exxx@\claimtype@}%
 \let\next@\relax
 \else
  \par
  \ifdim\lastskip<\smallskipamount\removelastskip\smallskip\fi
  \begingroup
  \noindent@@{\smc\ignorespaces#1\unskip
   \punct@{\null\colon@}\addspace@\enspace}%
  \nopunct@false\nospace@false
  \rm
  \DN@{\FNSSP@}%
 \fi
 \next@}
\def\enddemo{\par\endgroup\nopunct@false\nospace@false\smallskip}
\rightadd@\claim\to\nofrillslist@
\def\claim@F{\smc}
\def\claim@@@F{\csname\exxx@\claimtype@ @F\endcsname}
\def\claimformat@#1#2#3{%
 \medbreak\noindent@@{\smc#1 {\claim@@@F#2} #3%
 \punct@{\null.}\addspace@\enspace}\sl}
\def\claimformat@@#1#2{\claimformat@{\ignorespaces#1\unskip}%
 {\ifx\thelabel@@\empty\unskip\else\thelabel@@\fi}%
 {\ignorespaces#2\unskip}%
 \let\Claimformat@@\claimformat@@\FNSSP@}
\let\Claimformat@@\claimformat@@
\def\claim@@@P{\csname\exxx@\claimtype@ @P\endcsname}
\def\claim@@@Q{\csname\exxx@\claimtype@ @Q\endcsname}
\def\claim@@@S{\csname\exxx@\claimtype@ @S\endcsname}
\def\claim@@@N{\csname\exxx@\claimtype@ @N\endcsname}
\def\claim@@@C{\csname claim@C\claimclass@\endcsname}
\newcount\claim@C
\claim@C\z@
\let\claim@P\empty
\let\claim@Q\empty
\def\claim@S#1{#1\/}
\let\claim@N\arabic
\def\claim{\claim@true\let\claimclass@\empty
 \def\claimtype@{\claim}\FN@\claim@}
\def\claim@{%
 \ifx\next\c
  \let\next@\claim@c
 \else
  \ifx\next"%
   \let\next@\claim@q
  \else
   \begingroup\global\advance\claim@C\@ne
   {\noexpands@
    \xdef\Thelabel@@@{\number\claim@C}%
    \xdefThelabel@\claim@N
    \xdef\Thelabel@@@@{\claim@P\Thelabel@\claim@Q}%
    \xdefThelabel@@\claim@S
   }%
   \locallabel@
   \let\next@\Claimformat@@
  \fi
 \fi
 \next@}
\def\claim@c\c#1{\claim@true\begingroup
 \expandafter
 \ifx\csname claim@C#1\endcsname\relax
  \expandafter\newcount@\csname claim@C#1\endcsname
  \global\csname claim@C#1\endcsname\@ne
 \else
  \global\advance\csname claim@C#1\endcsname\@ne
 \fi
 \def\claimclass@{#1}%
 {\noexpands@
  \xdef\Thelabel@@@{\number\claim@@@C}%
  \xdefThelabel@\claim@@@N
  \xdef\Thelabel@@@@{\claim@@@P\Thelabel@\claim@@@Q}%
  \xdefThelabel@@\claim@@@S
 }%
 \locallabel@
 \FNSS@\claim@c@}
\def\claim@q"#1"{\begingroup
 {\let\pre\claim@@@P\let\post\claim@@@Q
  \let\style\claim@@@S\let\numstyle\claim@@@N
  \noexpands@
  \Qlabel@{#1}}%
 \locallabel@
 \FNSS@\claim@q@}
\def\claim@c@{\ifx\next"%
 \global\advance\claim@@@C\m@ne\let\next@\claim@cq
 \else\let\next@\Claimformat@@\fi\next@}
\def\claim@cq"#1"{{\let\pre\claim@@@P\let\post\claim@@@Q
 \let\style\claim@@@S\let\numstyle\claim@@@N
 \noexpands@
 \Qlabel@{#1}}%
 \locallabel@
 \FNSS@\Claimformat@@}
\def\claim@q@{\ifx\next\c\expandafter\claim@qc
 \else\expandafter\Claimformat@@\fi}
\def\claim@qc\c#1{\expandafter\ifx\csname claim@C#1\endcsname\relax
 \expandafter\newcount@\csname claim@C#1\endcsname
 \global\csname claim@C#1\endcsname\z@\fi
 \FNSS@\Claimformat@@}
\def\endclaim{\endgroup\claim@false\nopunct@false\nospace@false
 \let\Claimformat@@\claimformat@@\medbreak}
\Invalid@\claimclause
\def\newclaim{\FN@\newclaim@}
\def\newclaim@{\ifx\next\claimclause
 \DN@\claimclause##1{\newclaim@@{##1}}\else
 \DN@{\newclaim@@\relax}\fi\next@}
\def\claimlist@{\\\claim}
\newtoks\claim@i
\newtoks\claim@v
\let\noclaimclause@=F
\def\newclaim@@#1#2#3\c#4#5{\define#2{}%
 \rightadd@#2\to\claimlist@\rightadd@#2\to\nofrillslist@%
 \expandafter\def\csname\exstring@#2@P\endcsname{\claim@P}%
 \expandafter\def\csname\exstring@#2@Q\endcsname{\claim@Q}%
 \expandafter\def\csname\exstring@#2@S\endcsname{\claim@S}%
 \expandafter\def\csname\exstring@#2@N\endcsname{\claim@N}%
 \expandafter\def\csname\exstring@#2@F\endcsname{\claim@F}%
 \expandafter\def\csname end\exstring@#2\endcsname{\endclaim}%
 \expandafter\ifx\csname claim@C#4\endcsname\relax
  \expandafter\newcount@\csname claim@C#4\endcsname
  \global\csname claim@C#4\endcsname\z@\fi
 \edef\next@{\let\csname\exstring@#2@C\endcsname
   \csname claim@C#4\endcsname}\next@
 \def#2{\ifx\noclaimclause@ T\else#1\fi
  \global\claim@i{#1}\gdef\claim@iv{#4}\global\claim@v{#5}%
  \def\claimtype@{#2}\def\Claimformat@@{\claimformat@@{#5}}\claim@c\c{#4}}}
\def\shortenclaim#1#2{\define#2{}%
 \ismember@\claimlist@#1%
 \iftest@
  \rightadd@#2\to\nofrillslist@%
  \expandafter\def\csname\exstring@#2@P\endcsname
   {\csname\exstring@#1@P\endcsname}%
  \expandafter\def\csname\exstring@#2@Q\endcsname
   {\csname\exstring@#1@Q\endcsname}%
  \expandafter\def\csname\exstring@#2@S\endcsname
   {\csname\exstring@#1@S\endcsname}%
  \expandafter\def\csname\exstring@#2@N\endcsname
   {\csname\exstring@#1@N\endcsname}%
  \expandafter\def\csname\exstring@#2@F\endcsname
   {\csname\exstring@#1@F\endcsname}%
  \expandafter\def\csname end\exstring@#2\endcsname{\endclaim}%
  \edef\next@{\let\csname\exstring@#2@C\endcsname
    \csname claim\exstring@#1C\endcsname}\next@
  \setbox\z@\vbox{\let\noclaimclause@ T#1""\relax\endgroup}%
  \edef#2{\the\claim@i
   \def\noexpand\claimtype@{\noexpand#2}%
   \def\noexpand\Claimformat@@{\noexpand\claimformat@@{\the\claim@v}\relax}%
   \noexpand\claim@c\noexpand\c{\claim@iv}}%
 \else
  \Err@{\noexpand#1not yet created by \string\newclaim}%
 \fi}
\def\classtest@#1{\DN@{#1}\ifx\next@\claimclass@
 \test@true\else\test@false\fi}
\def\typetest@#1{\DN@{#1}\ifx\next@\claimtype@\test@true\else
  \test@false\fi}
\newif\iftoc@
\def\tocfile{\iftoc@\else\alloc@@7\write\chardef\sixt@@n\toc@
 \immediate\openout\toc@=\jobname.toc
 \alloc@@7\write\chardef\sixt@@n\tic@
 \immediate\openout\tic@=\jobname.tic
 \global\toc@true\fi}
\rightadd@\hl\to\nofrillslist@
\rightadd@\HL\to\overlonglist@
\def\HL@@C{\csname HL@C\HLlevel@\endcsname}
\def\HL@@P{\csname HL@P\HLlevel@\endcsname}
\def\HL@@Q{\csname HL@Q\HLlevel@\endcsname}
\def\HL@@S{\csname HL@S\HLlevel@\endcsname}
\def\HL@@N{\csname HL@N\HLlevel@\endcsname}
\def\HL@@F{\csname HL@F\HLlevel@\endcsname}
\def\HL@@@C{\csname\exxx@\HLtype@ @C\endcsname}
\def\HL@@@P{\csname\exxx@\HLtype@ @P\endcsname}
\def\HL@@@Q{\csname\exxx@\HLtype@ @Q\endcsname}
\def\HL@@@S{\csname\exxx@\HLtype@ @S\endcsname}
\def\HL@@@N{\csname\exxx@\HLtype@ @N\endcsname}
\def\HL#1{\expandafter
 \ifx\csname HL@C#1\endcsname\relax
  \DN@{\Err@{\string\HL#1 not defined in this style}}%
 \else
  \DN@{\gdef\HLlevel@{#1}\def\HLname@{\HL{#1}}\let\HLtype@\relax\FNSS@\HL@}%
 \fi
 \next@}%
\newif\ifquoted@
\let\aftertoc@\relax
\def\HL@{%
 \DN@"##1"##2\endHL{\def\entry@{##2}\quoted@true
  {\noexpands@
  \ifx\HLtype@\relax
   \let\pre\HL@@P\let\post\HL@@Q\let\style\HL@@S\let\numstyle\HL@@N
  \else
   \let\pre\HL@@@P\let\post\HL@@@Q\let\style\HL@@@S\let\numstyle\HL@@@N
  \fi
  \Qlabel@{##1}\let\style\relax\xdef\Qlabel@@@@{##1}%
  \xdef\Thepref@{\Thelabel@@@@}}%
  \csname HL@\HLlevel@\endcsname##2\endHL
  \let\pref\Thepref@
  \csname HL@I\HLlevel@\endcsname
  \csname HL@J\HLlevel@\endcsname
  \let\pref\pref@
  \HLtoc@	
  \aftertoc@
  \let\aftertoc@\relax\overlong@false}%
 \DNii@##1\endHL{\def\entry@{##1}\quoted@false
  {\noexpands@
  \ifx\HLtype@\relax
   \global\advance\HL@@C\@ne
   \xdef\Thelabel@@@{\number\HL@@C}%
   \xdefThelabel@{\HL@@N}%
   \xdef\Thelabel@@@@{\HL@@P\Thelabel@\HL@@Q}%
   \xdefThelabel@@{\HL@@S}%
  \else
   \global\advance\HL@@@C\@ne
   \xdef\Thelabel@@@{\number\HL@@@C}%
   \xdefThelabel@{\HL@@@N}%
   \xdef\Thelabel@@@@{\HL@@@P\Thelabel@\HL@@@Q}%
   \xdefThelabel@@{\HL@@@S}%
  \fi
  \xdef\Thepref@{\Thelabel@@@@}}%
  \csname HL@\HLlevel@\endcsname##1\endHL
  \let\pref\Thepref@
  \csname HL@I\HLlevel@\endcsname
  \csname HL@J\HLlevel@\endcsname
  \let\pref\pref@
  \HLtoc@
  \aftertoc@
  \let\aftertoc@\relax\overlong@false}%
 \ifx\next"\expandafter\next@\else\expandafter\nextii@\fi}%
\Invalid@\endHL
\def\hl@@C{\csname hl@C\hllevel@\endcsname}
\def\hl@@P{\csname hl@P\hllevel@\endcsname}
\def\hl@@Q{\csname hl@Q\hllevel@\endcsname}
\def\hl@@S{\csname hl@S\hllevel@\endcsname}
\def\hl@@N{\csname hl@N\hllevel@\endcsname}
\def\hl@@F{\csname hl@F\hllevel@\endcsname}
\def\hl@@@C{\csname\exxx@\hltype@ @C\endcsname}
\def\hl@@@P{\csname\exxx@\hltype@ @P\endcsname}
\def\hl@@@Q{\csname\exxx@\hltype@ @Q\endcsname}
\def\hl@@@S{\csname\exxx@\hltype@ @S\endcsname}
\def\hl@@@N{\csname\exxx@\hltype@ @N\endcsname}
\def\hl#1{\expandafter
 \ifx\csname hl@C#1\endcsname\relax
  \DN@{\Err@{\string\hl#1 not defined in this style}}%
 \else
  \DN@{\gdef\hllevel@{#1}\def\hlname@{\hl{#1}}\let\hltype@\relax\FNSS@\hl@}%
 \fi
 \next@}
\def\hl@{%
 \DN@"##1"##2{\def\entry@{##2}\quoted@true
  {\noexpands@
  \ifx\hltype@\relax
   \let\pre\hl@@P\let\post\hl@@Q\let\style\hl@@S\let\numstyle\hl@@N
  \else
   \let\pre\hl@@@P\let\post\hl@@@Q\let\style\hl@@@S\let\numstyle\hl@@@N
  \fi
  \Qlabel@{##1}\let\style\relax\xdef\Qlabel@@@@{##1}%
  \xdef\Thepref@{\Thelabel@@@@}}%
  \csname hl@\hllevel@\endcsname{##2}%
  \let\pref\Thepref@
  \csname hl@I\hllevel@\endcsname
  \csname hl@J\hllevel@\endcsname
  \let\pref\pref@
  \hltoc@
  \aftertoc@
  \let\aftertoc@\relax\nopunct@false\nospace@false\FNSSP@}%
 \DNii@##1{\def\entry@{##1}\quoted@false
  {\noexpands@
  \ifx\hltype@\relax
   \global\advance\hl@@C\@ne
   \xdef\Thelabel@@@{\number\hl@@C}%
   \xdefThelabel@{\hl@@N}%
   \xdef\Thelabel@@@@{\hl@@P\Thelabel@\hl@@Q}%
   \xdefThelabel@@{\hl@@S}%
  \else
   \global\advance\hl@@@C\@ne
   \xdef\Thelabel@@@{\number\hl@@@C}%
   \xdefThelabel@{\hl@@@N}%
   \xdef\Thelabel@@@@{\hl@@@P\Thelabel@\hl@@@Q}%
   \xdefThelabel@@{\hl@@@S}%
  \fi
  \xdef\Thepref@{\Thelabel@@@@}}%
  \csname hl@\hllevel@\endcsname{##1}%
  \let\pref\Thepref@
  \csname hl@I\hllevel@\endcsname
  \csname hl@J\hllevel@\endcsname
  \let\pref\pref@
  \hltoc@
  \aftertoc@
  \let\aftertoc@\relax\nopunct@false\nospace@false\FNSSP@}%
 \ifx\next"\expandafter\next@\else\expandafter\nextii@\fi}%
\def\six@#1#2 #3 #4 #5 #6 #7 {\DN@{#2}\ifx\next@\empty
 \DN@##1\six@{}\else
 \write#1{ #2 #3 #4 #5 #6 #7}\DN@{\six@#1}\fi
 \next@}
\def\Sixtoc@{\ifx\macdef@\empty\else
 \DN@##1##2\next@{\def\macdef@{##1##2}}%
 \expandafter\next@\macdef@\next@
 \edef\next@
  {\noexpand\six@\toc@\macdef@
  \space\space\space\space\space\space\space\space\space\space\space\space
  \noexpand\six@}%
 \next@\let\macdef@\relax\fi}
\def\QorThelabel@@@@{\ifquoted@
 \noexpand\noexpand\noexpand"\Qlabel@@@@\noexpand\noexpand\noexpand"\else
 \Thelabel@@@@\fi}
\def\HLtoc@{%
 \iftoc@
 \expandafter\expandafter\expandafter\unmacro@
  \expandafter\meaning\csname HL@W\HLlevel@\endcsname\unmacro@
  {\noexpands@\let\style\relax
   \edef\next@{\write\toc@{\noexpand\noexpand\expandafter\noexpand\HLname@
   {\macdef@}{\QorThelabel@@@@}}}%
  \next@}%
  \expandafter\unmacro@\meaning\entry@\unmacro@
  \Sixtoc@
  \write\toc@{\noexpand\Page{\number\pageno}{\page@N}%
   {\page@P}{\page@Q}^^J}%
 \fi}
\def\hltoc@{%
 \iftoc@
 \expandafter\expandafter\expandafter\unmacro@
  \expandafter\meaning\csname hl@W\hllevel@\endcsname\unmacro@
  {\noexpands@\let\style\relax
  \edef\next@{\write\toc@{%
   \ifnopunct@\noexpand\noexpand\noexpand\nopunct\fi
   \ifnospace@\noexpand\noexpand\noexpand\nospace\fi
   \noexpand\noexpand\expandafter\noexpand\hlname@
   {\macdef@}{\QorThelabel@@@@}}}%
  \next@}%
  \expandafter\unmacro@\meaning\entry@\unmacro@
  \Sixtoc@
  \write\toc@{\noexpand\Page{\number\pageno}{\page@N}%
   {\page@P}{\page@Q}^^J}%
 \fi}
\def\mainfile#1{\def\mainfile@{#1}}
\def\checkmainfile@{\ifx\mainfile@\undefined
 \Err@{No \noexpand\mainfile specified}\fi}
\expandafter\newcount@\csname HL@C1\endcsname
\csname HL@C1\endcsname\z@
\expandafter\def\csname HL@S1\endcsname#1{#1\null.}
\expandafter\let\csname HL@N1\endcsname\arabic
\expandafter\let\csname HL@P1\endcsname\empty
\expandafter\let\csname HL@Q1\endcsname\empty
\expandafter\def\csname HL@F1\endcsname{\bf}
\expandafter\let\csname HL@W1\endcsname\empty
\expandafter\newcount@\csname hl@C1\endcsname
\csname hl@C1\endcsname\z@
\expandafter\def\csname hl@S1\endcsname#1{#1\/}
\expandafter\let\csname hl@N1\endcsname\arabic
\expandafter\let\csname hl@P1\endcsname\empty
\expandafter\let\csname hl@Q1\endcsname\empty
\expandafter\def\csname hl@F1\endcsname{\bf}
\expandafter\let\csname hl@W1\endcsname\empty
\expandafter\def\csname HL@1\endcsname#1\endHL{\bigbreak
 {\locallabel@
  \global\setbox\@ne\vbox{\Let@\tabskip\hss@
  \halign to\hsize{\bf\hfil\ignorespaces##\unskip\hfil\cr
  \expandafter\ifx\csname HL@W1\endcsname\empty\else
   \csname HL@W1\endcsname\space\fi
  {\HL@@F\ifx\thelabel@@\empty\else\thelabel@@\space\fi}%
  \ignorespaces#1\crcr}}%
  }%
 \unvbox\@ne\nobreak\medskip}
\expandafter\def\csname hl@1\endcsname#1{\medbreak\noindent@@
 {\locallabel@
 \bf{\hl@@F\ifx\thelabel@@\empty\else\thelabel@@\space\fi}%
 \ignorespaces#1\unskip\punct@{\null.}\addspace@\enspace}}
\expandafter\def\csname HL@I1\endcsname{\Reset\hl1{1}%
 \ifx\pref\empty\newpre\hl1{}\else\newpre\hl1{\pref.}\fi}
\def\NameHL#1#2{\define#2{}%
 \expandafter\ifx\csname HL@R#1\endcsname\relax
 \else
  \def\nextiv@{\let\nextiii@}%
  \expandafter\nextiv@\csname HL@R#1\endcsname
  \expandafter\let\nextiii@\undefined
  \expandafter\let\csname\exxx@\nextiii@ @C\endcsname\relax
  \expandafter\let\csname\exxx@\nextiii@ @P\endcsname\relax
  \expandafter\let\csname\exxx@\nextiii@ @Q\endcsname\relax
  \expandafter\let\csname\exxx@\nextiii@ @S\endcsname\relax
  \expandafter\let\csname\exxx@\nextiii@ @N\endcsname\relax
  \expandafter\let\csname\exxx@\nextiii@ @F\endcsname\relax
  \expandafter\let\csname\exxx@\nextiii@ @W\endcsname\relax
  \expandafter\let\csname end\exxx@\nextiii@\endcsname\undefined
 \fi
 \expandafter\gdef\csname HL@R#1\endcsname{#2}%
 \expandafter\gdef\csname\exstring@#2@R\endcsname{{HL}{#1}}%
 \iftoc@\write\toc@{\noexpand\NameHL#1\noexpand#2^^J}\fi
 \rightadd@#2\to\overlonglist@
 \edef\next@{\let\csname\exstring@#2@C\endcsname\expandafter\noexpand
  \csname HL@C#1\endcsname}\next@
 \edef\next@{\let\csname\exstring@#2@P\endcsname\expandafter\noexpand
  \csname HL@P#1\endcsname}\next@
 \edef\next@{\let\csname\exstring@#2@Q\endcsname\expandafter\noexpand
  \csname HL@Q#1\endcsname}\next@
 \edef\next@{\let\csname\exstring@#2@S\endcsname\expandafter\noexpand
  \csname HL@S#1\endcsname}\next@
 \edef\next@{\let\csname\exstring@#2@N\endcsname\expandafter\noexpand
  \csname HL@N#1\endcsname}\next@
 \edef\next@{\let\csname\exstring@#2@F\endcsname\expandafter\noexpand
  \csname HL@F#1\endcsname}\next@
 \edef\next@{\let\csname\exstring@#2@W\endcsname\expandafter\noexpand
  \csname HL@W#1\endcsname}\next@
 \edef\next@{\def\noexpand#2####1\expandafter\noexpand
  \csname end\exstring@#2\endcsname
  {\def\noexpand\HLtype@{\noexpand#2}%
   \def\noexpand\HLname@{\noexpand#2}%
   \gdef\noexpand\HLlevel@{#1}%
   \noexpand\FNSS@\noexpand\HL@####1\noexpand\endHL}}%
  \next@
 \edef\next@{\noexpand\Invalid@\expandafter\noexpand
  \csname end\exstring@#2\endcsname}%
 \next@}
\def\Namehl#1#2{\define#2{}%
 \expandafter\ifx\csname hl@R#1\endcsname\relax
 \else
  \def\nextiv@{\let\nextiii@}%
  \expandafter\nextiv@\csname hl@R#1\endcsname
  \expandafter\let\nextiii@\undefined
  \expandafter\let\csname\exxx@\nextiii@ @C\endcsname\relax
  \expandafter\let\csname\exxx@\nextiii@ @P\endcsname\relax
  \expandafter\let\csname\exxx@\nextiii@ @Q\endcsname\relax
  \expandafter\let\csname\exxx@\nextiii@ @S\endcsname\relax
  \expandafter\let\csname\exxx@\nextiii@ @N\endcsname\relax
  \expandafter\let\csname\exxx@\nextiii@ @F\endcsname\relax
  \expandafter\let\csname\exxx@\nextiii@ @W\endcsname\relax
 \fi
 \expandafter\gdef\csname hl@R#1\endcsname{#2}%
 \expandafter\gdef\csname\exstring@#2@R\endcsname{{hl}{#1}}%
 \iftoc@\write\toc@{\noexpand\Namehl#1\noexpand#2^^J}\fi
 \rightadd@#2\to\nofrillslist@%
 \edef\next@{\let\csname\exstring@#2@C\endcsname\expandafter\noexpand
  \csname hl@C#1\endcsname}\next@
 \edef\next@{\let\csname\exstring@#2@P\endcsname\expandafter\noexpand
  \csname hl@P#1\endcsname}\next@
 \edef\next@{\let\csname\exstring@#2@Q\endcsname\expandafter\noexpand
  \csname hl@Q#1\endcsname}\next@
 \edef\next@{\let\csname\exstring@#2@S\endcsname\expandafter\noexpand
  \csname hl@S#1\endcsname}\next@
 \edef\next@{\let\csname\exstring@#2@N\endcsname\expandafter\noexpand
  \csname hl@N#1\endcsname}\next@
 \edef\next@{\let\csname\exstring@#2@F\endcsname\expandafter\noexpand
  \csname hl@F#1\endcsname}\next@
 \edef\next@{\let\csname\exstring@#2@W\endcsname\expandafter\noexpand
  \csname hl@W#1\endcsname}\next@
 \edef\next@{\def\noexpand#2{%
  \def\noexpand\hltype@{\noexpand#2}%
  \def\noexpand\hlname@{\noexpand#2}%
  \gdef\noexpand\hllevel@{#1}%
  \noexpand\FNSS@\noexpand\hl@}}%
 \next@}%
\def\Initialize{\FN@\Init@}
\def\Init@{\ifx\next\HL\let\next@\InitH@\else\ifx\next\hl\let\next@\InitH@
  \else\let\next@\InitS@\fi\fi\next@}
\def\InitH@#1#2{\expandafter\ifx\csname\exstring@#1@C#2\endcsname\relax
 \DN@{\Err@{\noexpand#1level #2 not defined in this style}}\else
 \DN@{\expandafter\gdef\csname\exstring@#1@J#2\endcsname}\fi\next@}
\def\InitC@#1#2{\edef\nextii@{\expandafter\noexpand\csname#1\endcsname{#2}}}
\def\InitS@#1{\expandafter\ifx\csname\exstring@#1@R\endcsname\relax
 \Err@{\noexpand#1not defined in this style}\let\next@\relax\else
 \DN@{\let\next@}\expandafter\next@\csname\exstring@#1@R\endcsname
 \expandafter\InitC@\next@
 \DN@{\expandafter\InitH@\nextii@}\fi\next@}
\def\value#1{\expandafter
 \ifx\csname\exstring@#1@C\endcsname\relax
  \expandafter\ifx\csname\exstring@#1@C1\endcsname\relax
   \DN@{\Err@{\noexpand\value can't be used with \string#1}}%
  \else
   \DN@{\value@#1}%
  \fi
 \else
  \DN@{\number\csname\exstring@#1@C\endcsname\relax}%
 \fi
 \next@}
\def\value@#1#2{\expandafter
 \ifx\csname\exstring@#1@C#2\endcsname\relax
  \DN@{\Err@{\string\value\string#1 can't be followed by \string#2}}%
 \else
  \DN@{\number\csname\exstring@#1@C#2\endcsname\relax}%
 \fi
 \next@}
\newcount\Value
\def\Evaluate#1{\expandafter
 \ifx\csname\exstring@#1@C\endcsname\relax
  \expandafter\ifx\csname\exstring@#1@C1\endcsname\relax
   \DN@{\Err@{\noexpand\Evaluate can't be used with \string#1}}%
  \else
   \DN@{\Evaluate@#1}%
  \fi
 \else
  \DN@{\global\Value\csname\exstring@#1@C\endcsname}%
 \fi
 \next@}
\def\Evaluate@#1#2{\expandafter
 \ifx\csname\exstring@#1@C#2\endcsname\relax
  \DN@{\Err@{\string\Evaluate\string#1 can't be followed by \string#2}}%
 \else
  \DN@{\global\Value\csname\exstring@#1@C#2\endcsname}%
 \fi\next@}
\def\pre#1{\expandafter
 \ifx\csname\exstring@#1@P\endcsname\relax
  \expandafter\ifx\csname\exstring@#1@P1\endcsname\relax
   \DN@{\Err@{\noexpand\pre can't be used with \string#1}}%
  \else
   \DN@{\pre@#1}%
  \fi
 \else
  \DN@{{\csname\exstring@#1@P\endcsname}}%
 \fi
 \next@}
\def\pre@#1#2{\expandafter
 \ifx\csname\exstring@#1@P#2\endcsname\relax
  \DN@{\Err@{\string\pre\string#1 can't be followed by \string#2}}%
 \else
  \DN@{{\csname\exstring@#1@P#2\endcsname}}%
 \fi
 \next@}
\def\post#1{\expandafter
 \ifx\csname\exstring@#1@Q\endcsname\relax
  \expandafter\ifx\csname\exstring@#1@Q1\endcsname\relax
   \DN@{\Err@{\noexpand\post can't be used with \string#1}}%
  \else
   \DN@{\post@#1}%
  \fi
 \else
  \DN@{{\csname\exstring@#1@Q\endcsname}}%
 \fi
 \next@}
\def\post@#1#2{\expandafter
 \ifx\csname\exstring@#1@Q#2\endcsname\relax
  \DN@{\Err@{\string\post\string#1 can't be followed by \string#2}}%
 \else
  \DN@{{\csname\exstring@#1@Q#2\endcsname}}%
 \fi
 \next@}
\def\style#1{\expandafter
 \ifx\csname\exstring@#1@S\endcsname\relax
  \expandafter\ifx\csname\exstring@#1@S1\endcsname\relax
   \DN@{\Err@{\noexpand\style can't be used with \string#1}}%
  \else
   \DN@{\style@#1}%
  \fi
 \else
  \DN@{\csname\exstring@#1@S\endcsname}%
 \fi
 \next@}
\def\style@#1#2{\expandafter
 \ifx\csname\exstring@#1@S#2\endcsname\relax
  \DN@{\Err@{\string\style\string#1 can't be followed by \string#2}}%
 \else
  \DN@{\csname\exstring@#1@S#2\endcsname}%
 \fi
 \next@}
\def\fontstyle#1{\expandafter
 \ifx\csname\exstring@#1@F\endcsname\relax
  \expandafter\ifx\csname\exstring@#1@F1\endcsname\relax
   \DN@{\Err@{\noexpand\fontstyle can't be used with \string#1}}%
  \else
   \DN@{\fontstyle@#1}%
  \fi
 \else
  \DN@##1{{\csname\exstring@#1@F\endcsname##1}}%
 \fi
 \next@}
\def\fontstyle@#1#2{\expandafter
 \ifx\csname\exstring@#1@F#2\endcsname\relax
  \DN@{\Err@{\string\fontstyle\string#1 can't be followed by \string#2}}%
 \else
  \DN@##1{{\csname\exstring@#1@F#2\endcsname##1}}%
 \fi
 \next@}
\def\Reset#1{\expandafter
 \ifx\csname\exstring@#1@C\endcsname\relax
  \expandafter\ifx\csname\exstring@#1@C1\endcsname\relax
   \DN@{\Err@{\noexpand\Reset can't be used with \string#1}}%
  \else
   \DN@{\Reset@#1}%
  \fi
 \else
  \DN@##1{\count@##1\relax\ifx#1\page\else\advance\count@\m@ne\fi
   \global\csname\exstring@#1@C\endcsname\count@}%
 \fi
 \next@}
\def\Reset@#1#2{\expandafter
 \ifx\csname\exstring@#1@C#2\endcsname\relax
  \DN@{\Err@{\string\Reset\string#1 can't be followed by \string#2}}%
 \else
  \DN@##1{\count@##1\relax\advance\count@\m@ne
   \global\csname\exstring@#1@C#2\endcsname\count@}%
 \fi
 \next@}
\def\Offset#1{\expandafter
 \ifx\csname\exstring@#1@C\endcsname\relax
  \expandafter\ifx\csname\exstring@#1@C1\endcsname\relax
   \DN@{\Err@{\noexpand\Offset can't be used with \string#1}}%
  \else
   \DN@{\Offset@#1}%
  \fi
 \else
  \DN@##1{\count@##1\relax\advance\count@\m@ne\global\advance
   \csname\exstring@#1@C\endcsname\count@}%
 \fi
 \next@}
\def\Offset@#1#2{\expandafter
 \ifx\csname\exstring@#1@C#2\endcsname\relax
  \DN@{\Err@{\string\Offset\string#1 can't be followed by \string#2}}%
 \else
  \DN@##1{\count@##1\relax\advance\count@\m@ne
   \global\advance\csname\exstring@#1@C#2\endcsname\count@}%
 \fi
 \next@}
\def\getR@#1#2{\def\nextiv@{\let\nextiii@}\expandafter\nextiv@
 \csname\exstring@#1@R#2\endcsname}
\def\letR@#1#2#3{\expandafter\let\csname#1@#3#2\endcsname\Next@}
\def\letR@@#1#2{\expandafter\let\csname\exstring@#1@#2\endcsname\Next@}
\def\newpre#1{\expandafter
 \ifx\csname\exstring@#1@P\endcsname\relax
  \expandafter\ifx\csname\exstring@#1@P1\endcsname\relax
   \DN@{\Err@{\noexpand\newpre can't be used with \string#1}}%
  \else
   \DN@{\newpre@#1}%
  \fi
 \else
  \DN@{%
   \DNii@{%
    \endgroup
    \expandafter\let\csname\exstring@#1@P\endcsname\Next@
    \expandafter\ifx\csname\exstring@#1@R\endcsname\relax\else
    \getR@#1{}\expandafter\letR@\nextiii@ P\fi
    }%
   \begingroup\noexpands@\afterassignment\nextii@\xdef\Next@}%
 \fi
 \next@}
\def\newpre@#1#2{\expandafter
 \ifx\csname\exstring@#1@P#2\endcsname\relax
  \DN@{\Err@{\string\newpre\string#1 can't be followed by \string#2}}%
 \else
  \DN@{%
   \DNii@{%
    \endgroup
    \expandafter\let\csname\exstring@#1@P#2\endcsname\Next@
    \expandafter\ifx\csname\exstring@#1@R#2\endcsname\relax\else
    \getR@#1{#2}\expandafter\letR@@\nextiii@ P\fi
    }%
   \begingroup\noexpands@\afterassignment\nextii@\xdef\Next@}%
 \fi
 \next@}
\def\newpost#1{\expandafter
 \ifx\csname\exstring@#1@Q\endcsname\relax
  \expandafter\ifx\csname\exstring@#1@Q1\endcsname\relax
   \DN@{\Err@{\noexpand\newpost can't be used with \string#1}}%
  \else
   \DN@{\newpost@#1}%
  \fi
 \else
  \DN@{%
   \DNii@{%
    \endgroup
    \expandafter\let\csname\exstring@#1@Q\endcsname\Next@
    \expandafter\ifx\csname\exstring@#1@R\endcsname\relax\else
    \getR@#1{}\expandafter\letR@\nextiii@ Q\fi
    }%
   \begingroup\noexpands@\afterassignment\nextii@\xdef\Next@}%
 \fi
 \next@}
\def\newpost@#1#2{\expandafter
 \ifx\csname\exstring@#1@Q#2\endcsname\relax
  \DN@{\Err@{\string\newpost\string#1 can't be followed by \string#2}}%
 \else
  \DN@{%
   \DNii@{%
    \endgroup
    \expandafter\let\csname\exstring@#1@Q#2\endcsname\Next@
    \expandafter\ifx\csname\exstring@#1@R#2\endcsname\relax\else
    \getR@#1{#2}\expandafter\letR@@\nextiii@ Q\fi
    }%
   \begingroup\noexpands@\afterassignment\nextii@\xdef\Next@}%
 \fi
 \next@}
\def\newstyle#1{\expandafter
 \ifx\csname\exstring@#1@S\endcsname\relax
  \expandafter\ifx\csname\exstring@#1@S1\endcsname\relax
   \DN@{\Err@{\noexpand\newstyle can't be used
    with \string#1}}%
  \else
   \DN@{\newstyle@#1}%
  \fi
 \else
  \DN@{%
   \DNii@{%
    \expandafter\let\csname\exstring@#1@S\endcsname\Next@
    \expandafter\ifx\csname\exstring@#1@R\endcsname\relax\else
    \getR@#1{}\expandafter\letR@\nextiii@ S\fi
    }%
   \afterassignment\nextii@\gdef\Next@}%
 \fi
 \next@}
\def\newstyle@#1#2{\expandafter
 \ifx\csname\exstring@#1@S#2\endcsname\relax
  \DN@{\Err@{\string\newstyle\string#1 can't be followed by
   \string#2}}%
 \else
  \DN@{%
   \DNii@{%
    \expandafter\let\csname\exstring@#1@S#2\endcsname\Next@
    \expandafter\ifx\csname\exstring@#1@R#2\endcsname\relax\else
    \getR@#1{#2}\expandafter\letR@@\nextiii@ S\fi
    }%
   \afterassignment\nextii@\gdef\Next@}%
 \fi
 \next@}
\def\newnumstyle#1{\expandafter
 \ifx\csname\exstring@#1@N\endcsname\relax
  \expandafter\ifx\csname\exstring@#1@N1\endcsname\relax
   \DN@{\Err@{\noexpand\newnumstyle can't be used with
    \string#1}}%
  \else
   \DN@{\newnumstyle@#1}%
  \fi
 \else
  \DN@##1{%
   \gdef\Next@{##1}%
    \expandafter\let\csname\exstring@#1@N\endcsname\Next@
    \expandafter\ifx\csname\exstring@#1@R\endcsname\relax\else
    \getR@#1{}\expandafter\letR@\nextiii@ N\fi
    }%
 \fi
 \next@}
\def\newnumstyle@#1#2{\expandafter
 \ifx\csname\exstring@#1@N#2\endcsname\relax
  \DN@{\Err@{\string\newnumstyle\string#1 can't be followed by
   \string#2}}%
 \else
  \DN@##1{%
   \gdef\Next@{##1}%
    \expandafter\let\csname\exstring@#1@N#2\endcsname\Next@
    \expandafter\ifx\csname\exstring@#1@R#2\endcsname\relax\else
    \getR@#1{#2}\expandafter\letR@@\nextiii@ N\fi
    }%
  \fi
 \next@}
\def\newfontstyle#1{\expandafter
 \ifx\csname\exstring@#1@F\endcsname\relax
  \expandafter\ifx\csname\exstring@#1@F1\endcsname\relax
   \DN@{\Err@{\noexpand\newfontstyle can't be used with
    \string#1}}%
  \else
   \DN@{\newfontstyle@#1}%
  \fi
 \else
  \DN@##1{%
   \gdef\Next@{##1}%
    \expandafter\let\csname\exstring@#1@F\endcsname\Next@
    \expandafter\ifx\csname\exstring@#1@R\endcsname\relax\else
    \getR@#1{}\expandafter\letR@\nextiii@ F\fi
    }%
 \fi
 \next@}
\def\newfontstyle@#1#2{\expandafter
 \ifx\csname\exstring@#1@F#2\endcsname\relax
  \DN@{\Err@{\string\newfontstyle\string#1 can't be followed by
   \string#2}}%
 \else
  \DN@##1{%
   \gdef\Next@{##1}%
    \expandafter\let\csname\exstring@#1@F#2\endcsname\Next@
    \expandafter\ifx\csname\exstring@#1@R#2\endcsname\relax\else
    \getR@#1{#2}\expandafter\letR@@\nextiii@ F\fi
    }%
 \fi
 \next@}
\def\word#1{\expandafter
 \ifx\csname\exstring@#1@W\endcsname\relax
  \expandafter\ifx\csname\exstring@#1@W1\endcsname\relax
   \DN@{\Err@{\noexpand\word can't be used with \string#1}}%
  \else
   \DN@{\word@#1}%
  \fi
 \else
  \DN@{{\csname\exstring@#1@W\endcsname}}%
 \fi
 \next@}
\def\word@#1#2{\expandafter
 \ifx\csname\exstring@#1@W#2\endcsname\relax
  \DN@{\Err@{\string\word\noexpand#1can't be followed by \string#2}}%
 \else
  \DN@{{\csname\exstring@#1@W#2\endcsname}}%
 \fi
 \next@}
\def\newword#1{\expandafter
 \ifx\csname\exstring@#1@W\endcsname\relax
  \expandafter\ifx\csname\exstring@#1@W1\endcsname\relax
   \DN@{\Err@{\noexpand\newword can't be used  with \string#1}}%
  \else
   \DN@{\newword@#1}%
  \fi
 \else
  \DN@{%
   \DNii@{%
    \expandafter\let\csname\exstring@#1@W\endcsname\Next@
    \expandafter\ifx\csname\exstring@#1@R\endcsname\relax\else
     \getR@#1{}\expandafter\letR@\nextiii@ W\fi
    }%
   \afterassignment\nextii@\gdef\Next@}%
 \fi
 \next@}
\def\newword@#1#2{\expandafter
 \ifx\csname\exstring@#1@W#2\endcsname\relax
  \DN@{\Err@{\string\newword\noexpand#1can't be followed by \string#2}}%
 \else
  \DN@{%
   \DNii@{%
    \expandafter\let\csname\exstring@#1@W#2\endcsname\Next@
    \expandafter\ifx\csname\exstring@#1@R#2\endcsname\relax\else
     \getR@#1{#2}\expandafter\letR@@\nextiii@ W\fi
    }%
   \afterassignment\nextii@\gdef\Next@}%
 \fi
 \next@}
\newif\iffn@
\newcount\footmark@C
\footmark@C\z@
\def\footmark@S#1{$^{#1}$}
\let\footmark@N\arabic
\def\footmark@F{\rm}
\def\foottext@S#1{$^{#1}$}
\def\foottext@F{\rm}
\let\modifyfootnote@\relax
\def\modifyfootnote#1{\def\modifyfootnote@{#1}}
\def\vfootnote@#1{\insert\footins
 \bgroup
 \floatingpenalty\@MM\interlinepenalty\interfootnotelinepenalty
 \leftskip\z@\rightskip\z@\spaceskip\z@\xspaceskip\z@
 \rm\splittopskip\ht\strutbox\splitmaxdepth\dp\strutbox
 \locallabel@\noindent@@{\foottext@F#1}\modifyfootnote@
 \footstrut\FN@\fo@t}
\def\fo@t{\ifcat\bgroup\noexpand\next\expandafter\f@@t\else
 \expandafter\f@t\fi}
\def\f@t#1{#1\@foot}
\def\f@@t{\bgroup\aftergroup\@foot\afterassignment\FNSSP@\let\next@}
\def\@foot{\unskip\lower\dp\strutbox\vbox to\dp\strutbox{}\egroup
 \iffn@\expandafter\fn@false\else
 \expandafter\postvanish@\fi}
\newif\ifplainfn@
\plainfn@true
\def\fancyfootnotes{\plainfn@false}
\newcount\fancyfootmarkcount@
\fancyfootmarkcount@\z@
\newcount\lastfnpage@
\lastfnpage@-\@M
\let\justfootmarklist@\empty
\def\footmark{\let\@sf\empty
 \ifhmode\edef\@sf{\spacefactor\the\spacefactor}\/\fi
 \DN@{\ifx"\next\expandafter\nextii@\else\expandafter\footmark@\fi}%
 \DNii@"##1"{%
  \iffirstchoice@
   {\let\style\footmark@S\let\numstyle\footmark@N
   \footmark@F##1%
   \noexpands@
   \let\style\foottext@S
   \Qlabel@{##1}%
   }%
   \iffn@\else
    {\noexpands@
    \xdef\Next@{{\Thelabel@}{\Thelabel@@}{\Thelabel@@@}{\Thelabel@@@@}}%
    }%
    \expandafter\rightappend@\Next@\to\justfootmarklist@
   \fi
  \fi
  \@sf\relax}%
 \FN@\next@}
\def\footmark@{%
 \iffirstchoice@
  \global\advance\footmark@C\@ne
  \ifplainfn@
   \xdef\adjustedfootmark@{\number\footmark@C}%
  \else
   {\let\\\or\xdef\Next@{\ifcase\number\footmark@C\fnpages@\else
     -\@M\fi}}%
   \ifnum\Next@=-\@M
    \xdef\adjustedfootmark@{\number\footmark@C}%
   \else
    \ifnum\Next@=\lastfnpage@
     \global\advance\fancyfootmarkcount@\@ne
    \else
     \global\fancyfootmarkcount@\@ne
     \global\lastfnpage@\Next@
    \fi
    \xdef\adjustedfootmark@{\number\fancyfootmarkcount@}%
   \fi
  \fi
  {\noexpands@
  \xdef\Thelabel@@@{\adjustedfootmark@}%
  \xdefThelabel@\footmark@N
  \xdef\Thelabel@@@@{\Thelabel@}%
  \xdefThelabel@@\foottext@S
  }%
  \iffn@\else
   {\noexpands@
   \xdef\Next@{{\Thelabel@}{\Thelabel@@}{\Thelabel@@@}{\Thelabel@@@@}}%
   }%
   \expandafter\rightappend@\Next@\to\justfootmarklist@
  \fi
  \ifplainfn@
  \else
   \edef\next@{\write\laxwrite@{F\noexpand\the\pageno}}\next@
  \fi
 \fi
 \footmark@S{\footmark@N{\adjustedfootmark@}}%
 \@sf\relax}
\def\foottext{\prevanish@
 \ifx\justfootmarklist@\empty
  \Err@{There is no \noexpand\footmark for this \string\foottext}\fi
 \DN@\\##1##2\next@{\DN@{##1}\gdef\justfootmarklist@{##2}}%
 \expandafter\next@\justfootmarklist@\next@
 \expandafter\foottext@\next@}
\def\foottext@#1#2#3#4{{\noexpands@
  \xdef\Thelabel@{#1}\xdef\Thelabel@@{#2}%
  \xdef\Thelabel@@@{#3}\xdef\Thelabel@@@@{#4}}%
  \vfootnote@{\thelabel@@}}
\rightadd@\foottext\to\vanishlist@
\def\footnote{\fn@true
 \let\@sf\empty
 \ifhmode\edef\@sf{\spacefactor\the\spacefactor}\/\fi
 \DN@{\ifx"\next\expandafter\nextii@\else\expandafter\nextiii@\fi}%
 \DNii@"##1"{\footmark"##1"\vfootnote@{\let\style\foottext@S
  \let\numstyle\footmark@N##1}}%
 \def\nextiii@{\footmark\vfootnote@{\foottext@S{\footmark@N
  {\adjustedfootmark@}}}}%
 \FN@\next@}
\newdimen\litindent
\litindent20\p@
\newbox\litbox@
\newbox\Litbox@
\newcount\interlitpenalty@
\interlitpenalty@\@M
\newcount\litlines@
{\obeyspaces\gdef\defspace@{\def {\allowbreak\hskip.5emminus.15em}}}
{\obeylines\gdef\letM@{\let^^M\CtrlM@}}
\def\CtrlM@{\egroup
 \ifcase\litlines@\advance\litlines@\@ne\or
 \box\litbox@\advance\litlines@\@ne\else
 \penalty\interlitpenalty@\box\litbox@\fi
 \Lit@}
\def\Lit@{\setbox\litbox@\hbox\bgroup\litdefs@\hskip\litindent}
\newcount\littab@
\littab@8
\def\littab#1{\littab@#1\relax}
{\catcode`\^^I=\active\gdef\letTAB@{\let^^I\TAB@}}
\def\TAB@{\egroup
 \dimen@\wd\litbox@
 \advance\dimen@-\litindent
 \setboxz@h{\tt0}%
 \dimen@ii\littab@\wdz@
 \divide\dimen@\dimen@ii
 \multiply\dimen@\dimen@ii
 \advance\dimen@\littab@\wdz@
 \advance\dimen@\litindent
 \setbox\litbox@\hbox\bgroup\litdefs@\hbox to\dimen@{\unhbox\litbox@\hfil}}
{\catcode`\`=\active\gdef`{\relax\lq}}
\let\litbs@\relax
\let\litbs@@\relax
\def\litbackslash#1{%
 \edef\litbs@{\catcode`\string#1=\z@
 \def\noexpand\litbs@@{\def\expandafter\noexpand\csname\string#1\endcsname
  {\char`\string#1}}}}
\def\litcodes@{\catcode`\\=12
 \catcode`\{=12 \catcode`\}=12
 \catcode`\$=12 \catcode`\&=12
 \catcode`\#=12
 \catcode`\^=12 \catcode`\_=12
 \catcode`\@=12 \catcode`\~=12 \catcode`\"=12
 \catcode`\;=12 \catcode`\:=12 \catcode`\!=12 \catcode`\?=12
 \catcode`\%=12 \litbs@\catcode`\`=\active\obeyspaces\defspace@}
\def\activate@#1#2{{\lccode`\~=`#2%
 \lowercase{%
  \if0#1%
  \gdef\Next@{\def~{\egroup\endgroup\bigskip\vskip-\parskip
   \def\next@{\noindent@@\FN@\pretendspace@}\FNSS@\next@}}\else
  \gdef\Next@{\def~{\egroup\egroup\endgroup}}\fi
  }%
 }}
\def\litdefs@{\let\0\empty\let\1\litdelim@\def\ {\char32 }\litbs@@}%
\def\litdelimiter#1{%
 \edef\litdelim@{\char`#1}%
 \def\lit#1{\leavevmode\begingroup\litcodes@\litdefs@
  \tt\hyphenchar\tentt\m@ne\lit@}%
 \def\lit@##1#1{##1\endgroup\null}%
 \def\Lit#1{\ifhmode$$\abovedisplayskip\bigskipamount
  \abovedisplayshortskip\bigskipamount
  \belowdisplayskip\z@\belowdisplayshortskip\z@
  \postdisplaypenalty\@M
  $$\vskip-\baselineskip\else\bigskip\fi
  \begingroup\litlines@\z@
  \catcode`#1=\active\activate@0#1\Next@
  \def\displaybreak{\egroup\break\litlines@\z@\Lit@}%
  \def\allowdisplaybreak{\egroup\allowbreak\litlines@\z@\Lit@}%
  \def\allowdisplaybreaks{\egroup\allowbreak\interlitpenalty@\z@
   \litlines@\z@\Lit@}%
  \litcodes@\tt\catcode`\^^I=\active\letTAB@
  \obeylines\letM@\Lit@}%
 \def\Litbox##1=#1{\begingroup\ifodd##1\relax\aftergroup\global\fi
  \aftergroup\setbox\aftergroup##1\aftergroup\box\aftergroup\Litbox@
  \def\allowdisplaybreak{\egroup\allowbreak\litlines@\z@\Lit@}%
  \def\allowdisplaybreaks{\egroup\allowbreak\interlitpenalty@\z@
   \litlines@\z@\Lit@}%
  \catcode`#1=\active\activate@1#1\Next@
  \litcodes@\tt\catcode`\^^I=\active\letTAB@
  \obeylines\letM@\global\setbox\Litbox@\vbox\bgroup\litindent\z@%
  \litlines@\z@\Lit@}%
}
\newbox\titlebox@
\setbox\titlebox@\vbox{}
\rightadd@\title\to\overlonglist@
\def\title{\begingroup\Let@
 \global\setbox\titlebox@\vbox\bgroup\tabskip\hss@
 \halign to\hsize\bgroup\bf\hfil\ignorespaces##\unskip\hfil\cr}
\def\endtitle{\crcr\egroup\egroup\endgroup\overlong@false}
\newbox\authorbox@
\rightadd@\author\to\overlonglist@
\def\author{\begingroup\Let@
 \global\setbox\authorbox@\vbox\bgroup\tabskip\hss@
 \halign to\hsize\bgroup\rm\hfil\ignorespaces##\unskip\hfil\cr}
\def\endauthor{\crcr\egroup\egroup\endgroup\overlong@false}
\newbox\affilbox@
\def\affil{\begingroup\Let@
 \global\setbox\affilbox@\vbox\bgroup\tabskip\hss@
 \halign to\hsize\bgroup\rm\hfil\ignorespaces##\unskip\hfil\cr}%
\def\endaffil{\crcr\egroup\egroup\endgroup\overlong@false}
\let\date@\relax
\def\date#1{\gdef\date@{\ignorespaces#1\unskip}}
\def\today{\ifcase\month\or January\or February\or March\or April\or May\or
 June\or July\or August\or September\or October\or November\or December\fi
 \space\number\day, \number\year}
\def\maketitle{\hrule\height\z@\vskip-\topskip
 \vskip24\p@ plus12\p@ minus12\p@
 \unvbox\titlebox@
 \ifvoid\authorbox@\else\vskip12\p@ plus6\p@ minus3\p@\unvbox\authorbox@\fi
 \ifvoid\affilbox@\else\vskip10\p@ plus5\p@ minus2\p@\unvbox\affilbox@\fi
 \ifx\date@\relax\else\vskip6\p@ plus2\p@ minus\p@\centerline{\rm\date@}\fi
 \vskip18\p@ plus12\p@ minus6\p@}
\def\cite{%
 \DNii@(##1)##2{{\rm[}{##2}, {##1\/}{\rm]}}%
 \def\nextiii@##1{{\rm[}{##1\/}{\rm]}}%
 \DN@{\ifx\next(\expandafter\nextii@\else\expandafter\nextiii@\fi}%
 \FN@\next@}
\def\makebib@W{Bibliography}
\def\makebib{\begingroup\rm\bigbreak\centerline{\smc\makebib@W}%
 \nobreak\medskip
 \sfcode`\.=\@m\everypar{}\parindent\z@
 \def\nopunct{\nopunct@true}\def\nospace{\nospace@true}%
 \nopunct@false\nospace@false
 \def\lkerns@{\null\kern\m@ne sp\kern\@ne sp}%
 \def\nkerns@{\null\kern-\tw@ sp\kern\tw@ sp}%
}
\let\endmakebib\endgroup
\newif\ifnoprepunct@
\newif\ifnoprespace@
\newif\ifnoquotes@
\def\noprepunct{\noprepunct@true}
\def\noprespace{\noprespace@true}
\def\noquotes{\noquotes@true}
\newbox\nobox@
\newbox\keybox@
\newbox\bybox@
\newbox\paperbox@
\newbox\paperinfobox@
\newbox\jourbox@
\newbox\volbox@
\newbox\issuebox@
\newbox\yrbox@
\newbox\pgbox@
\newbox\ppbox@
\newbox\bookbox@
\newbox\inbookbox@
\newbox\bookinfobox@
\newbox\publbox@
\newbox\publaddrbox@
\newbox\edbox@
\newbox\edsbox@
\newbox\langbox@
\newbox\translbox@
\newbox\finalinfobox@
\def\setbibinfo@#1{\edef\next@{\ifnopunct@1\else0\fi
 \ifnospace@1\else0\fi\ifnoprepunct@1\else0\fi\ifnoprespace@1\else0\fi
 \ifnoquotes@1\else0\fi}%
 \DNii@{00000}%
 \ifx\next@\nextii@\else\xdef\bibinfo@{\bibinfo@\the#1,\next@}%
 \fi}
\def\getbibinfo@#1{\ifx\bibinfo@\empty
 \let\next@0\let\nextii@0\let\nextiii@0\let\nextiv@0\let\nextv@0\else
 \edef\next@{\def
  \noexpand\next@####1\the#1,####2####3####4####5####6####7\noexpand\next@
  {\let\noexpand\next@####2\let\noexpand\nextii@####3%
  \let\noexpand\nextiii@####4\let\noexpand\nextiv@####5%
  \let\noexpand\nextv@####6}%
  \noexpand\next@\bibinfo@\the#1,00000\noexpand\next@}\next@
 \fi}
\newif\ifbookinquotes@
\def\bookinquotes{\bookinquotes@true}
\newif\ifpaperinquotes@
\def\paperinquotes{\paperinquotes@true}
\newif\ifininbook@
\def\ininbook{\ininbook@true}
\newif\ifopenquotes@
\def\closequotes@{\ifopenquotes@''\openquotes@false\fi}
\newif\ifbeginbib@
\newif\ifendbib@
\newif\ifprevjour@
\newif\ifprevbook@
\newdimen\bibindent@
\bibindent@20\p@
\def\bib{\global\let\bibinfo@\empty\global\let\translinfo@\relax\beginbib@true
 \begingroup\noindent@
 \hangindent\bibindent@\hangafter\@ne\bib@}
\def\v@id#1{\setbox#1\box\voidb@x}
\def\bib@{\v@id\nobox@\v@id\keybox@\v@id\bybox@\v@id\paperbox@
 \v@id\paperinfobox@\v@id\jourbox@\v@id\volbox@\v@id\issuebox@
 \v@id\yrbox@\v@id\pgbox@\v@id\ppbox@\v@id\bookbox@\v@id\inbookbox@
 \v@id\bookinfobox@\v@id\publbox@\v@id\publaddrbox@\v@id\edbox@
 \v@id\edsbox@\v@id\langbox@\v@id\translbox@\v@id\finalinfobox@
 \bgroup}
\def\Setnonemptybox@#1#2{\unskip\setbibinfo@#1\egroup#2%
 \def\aftergroup@{\ifdim\wd#1=\z@\setbox#1\box\voidb@x\fi}%
 \setbox#1\vbox\bgroup\aftergroup\aftergroup@\hsize\maxdimen\leftskip\z@
 \rightskip\z@\hbadness\@M\hfuzz\maxdimen\noindent}
\def\setnonemptybox@#1{\Setnonemptybox@#1\relax}
\def\no{\setnonemptybox@\nobox@}
\def\key{\setnonemptybox@\keybox@\bf}
\def\by{\setnonemptybox@\bybox@}
\def\bysame{\setnonemptybox@\bybox@\leaders\hrule\hskip3em\null}
\def\paper{\setnonemptybox@\paperbox@
 \ifpaperinquotes@\getbibinfo@\paperbox@
 \if\nextv@1\else``\fi\else\it\fi}
\def\paperinfo{\setnonemptybox@\paperinfobox@}
\def\jour{\Setnonemptybox@\jourbox@\prevjour@true}
\def\vol{\setnonemptybox@\volbox@\bf}
\def\issue{\setnonemptybox@\issuebox@}
\def\yr{\setnonemptybox@\yrbox@}

\def\pg{\setnonemptybox@\pgbox@}
\def\pp{\setnonemptybox@\ppbox@}
\def\book{\Setnonemptybox@\bookbox@\prevbook@true
 \ifbookinquotes@\getbibinfo@\bookbox@
 \if\nextv@1\else``\fi\else\it\fi}
\def\inbook{\Setnonemptybox@\inbookbox@\prevbook@true
 \ifininbook@ in \fi\ifbookinquotes@\getbibinfo@\inbookbox@
 \if\nextv@1\else``\fi\fi}
\def\bookinfo{\setnonemptybox@\bookinfobox@}
\def\publ{\setnonemptybox@\publbox@}
\def\publaddr{\setnonemptybox@\publaddrbox@}
\def\ed{\setnonemptybox@\edbox@}
\def\eds{\setnonemptybox@\edsbox@}
\def\lang{\setnonemptybox@\langbox@}
\def\finalinfo{\setnonemptybox@\finalinfobox@}
\def\setboxzl@{\setbox\z@\lastbox}
\def\getbox@#1{\setbox\z@\vbox{\vskip-\@M\p@
 \unvbox#1%
 \setboxzl@
 \global\setbox\@ne\hbox{\unhbox\z@\unskip\unskip\unpenalty}%
 \ifdim\lastskip=-\@M\p@\else
 \loop\ifdim\lastskip=-\@M\p@
 \else\unskip\unpenalty\setboxzl@
 \global\setbox\@ne\hbox{\unhbox\z@\unhbox\@ne}%
 \repeat\fi}%
 \unhbox\@ne}
\def\adjustpunct@#1{\count@\lastkern
 \ifnum\count@=\z@#1\closequotes@\else
 \ifnum\count@>\tw@#1\closequotes@\else
 \ifnum\count@<-\tw@#1\closequotes@\else
  \unkern\unkern\setboxzl@
  \skip@\lastskip\unskip
  \count@@\lastpenalty\unpenalty
  \ifnum\count@=\tw@\unskip\setboxzl@\fi
  \ifdim\skip@=\z@\else\hskip\skip@\fi
  #1\closequotes@
  \ifnum\count@=\tw@\null\hfill\fi
  \penalty\count@@
 \fi\fi\fi}
\def\prepunct@#1#2{\getbibinfo@#2%
 \ifnopunct@
 \else
  \if\nextiii@0\adjustpunct@#1\fi
 \fi
 \closequotes@
 \ifnospace@
 \else
  \if\nextiv@0\space\else\fi
 \fi
 \nopunct@false\nospace@false
 \if\next@1\nopunct@true\fi
 \if\nextii@1\nospace@true\fi}
\def\ppunbox@#1#2{\prepunct@{#1}#2%
 \getbox@#2}
\let\semicolon@;
\def\endbib@{%
 \ifbeginbib@
  \ifvoid\nobox@
   \ifvoid\keybox@\else\hbox to\bibindent@{[\getbox@\keybox@]\hss}\fi
  \else\hbox to\bibindent@{\hss\getbox@\nobox@. }\fi
  \ifvoid\bybox@\else\getbox@\bybox@\fi
 \else
  \nopunct@true
  \ifvoid\bybox@\else\ppunbox@\relax\bybox@\fi
 \fi
 \ifvoid\translbox@\else\ppunbox@,\translbox@\fi
 \ifvoid\paperbox@\else\ppunbox@,\paperbox@\ifpaperinquotes@
  \if\nextv@1\else\openquotes@true\fi\fi
 \fi
 \ifvoid\paperinfobox@\else\ppunbox@,\paperinfobox@\fi
 \test@false
 \ifvoid\jourbox@\else\test@true\ppunbox@,\jourbox@\fi
 \ifprevjour@\test@true\fi
 \iftest@
  \ifvoid\volbox@\else\ppunbox@\relax\volbox@\fi
  \ifvoid\issuebox@
   \else\prepunct@\relax\issuebox@ no.~\getbox@\issuebox@\fi
  \ifvoid\yrbox@\else\prepunct@\relax\yrbox@(\getbox@\yrbox@)\fi
  \ifvoid\ppbox@\else\ppunbox@,\ppbox@\fi
  \ifvoid\pgbox@\else\prepunct@,\pgbox@ p.~\getbox@\pgbox@\fi
 \fi
 \test@false
 \ifvoid\bookbox@\else\test@true\ppunbox@,\bookbox@\ifbookinquotes@
  \if\nextv@1\else\openquotes@true\fi\fi\fi
 \ifvoid\inbookbox@\else\test@true\ppunbox@,\inbookbox@\ifbookinquotes@
  \if\nextv@1\else\openquotes@true\fi\fi\fi
 \ifprevbook@\test@true\fi
 \iftest@
  \ifvoid\edbox@\else\prepunct@\relax\edbox@(\getbox@\edbox@, ed.)\fi
  \ifvoid\edsbox@\else\prepunct@\relax\edsbox@(\getbox@\edsbox@, eds.)\fi
  \ifvoid\bookinfobox@\else\ppunbox@,\bookinfobox@\fi
  \ifvoid\publbox@\else\ppunbox@,\publbox@\fi
  \ifvoid\publaddrbox@\else\ppunbox@,\publaddrbox@\fi
  \ifvoid\yrbox@\else\ppunbox@,\yrbox@\fi
  \ifvoid\ppbox@\else\prepunct@,\ppbox@ pp.~\getbox@\ppbox@\fi
  \ifvoid\pgbox@\else\prepunct@,\pgbox@ p.~\getbox@\pgbox@\fi
 \fi
 \ifvoid\finalinfobox@
  \ifendbib@
   \ifnopunct@\else.\closequotes@\fi
  \else
  \ifvoid\langbox@\else\space(\getbox@\langbox@)\fi
   \/\semicolon@\closequotes@
  \fi
 \else
  \ifendbib@
   \ppunbox@{.\spacefactor3000\relax}\finalinfobox@
    \ifnopunct@\else.\fi
  \else
   \ppunbox@,\finalinfobox@\/\semicolon@\fi
 \fi
 \ifvoid\langbox@\else\space(\getbox@\langbox@)\fi
}
\def\endbib{\unskip\egroup\endbib@true\endbib@\par\endgroup}
\def\morebib{\unskip\egroup
 \endbib@false\endbib@
 \global\let\bibinfo@\empty\beginbib@false
 \bib@}
\def\anotherbib{\unskip\egroup
 \endbib@false\endbib@
 \global\let\bibinfo@\empty\beginbib@false
 \prevjour@false\prevbook@false\bib@}
\def\transl{\unskip
 \xdef\translinfo@{\the\translbox@,\ifnopunct@1\else0\fi
 \ifnospace@1\else0\fi\ifnoprepunct@1\else0\fi\ifnoprespace@1\else0\fi0}%
 \egroup\endbib@false\endbib@
 \global\let\bibinfo@\translinfo@\beginbib@false
 \bib@
 \egroup
 \def\aftergroup@{\ifdim\wd\translbox@=\z@\setbox\translbox@\box\voidb@x\fi}%
 \setbox\translbox@\vbox\bgroup\aftergroup\aftergroup@
 \hsize\maxdimen\leftskip\z@\rightskip\z@\hbadness\@M\hfuzz\maxdimen
 \noindent}
\newwrite\auxwrite@
\newread\bbl@
\def\UseBibTeX{\immediate\openout\auxwrite@=\jobname.aux
 \let\cite\BTcite@
 \def\nocite##1{\immediate\write\auxwrite@{\string\citation{##1}}}%
 \def\bibliographystyle##1{\immediate\write\auxwrite@{\string
  \bibstyle{##1}}}%
 \def\bibliography@W{Bibliography}%
 \def\bibliography##1{\immediate\write\auxwrite@{\string\bibdata{##1}}%
  \immediate\openin\bbl@=\jobname.bbl
  \ifeof\bbl@
   \W@{No .bbl file}%
  \else
   \immediate\closein\bbl@
   \begingroup\input bibtex \input\jobname.bbl \endgroup
  \fi}%
 }
\def\BTcite@{%
 \DNii@(##1)##2{{\rm[}\BTcite@@##2,\BTcite@@{\rm, }{##1\/}{\rm]}%
  \immediate\write\auxwrite@{\string\citation{##2}}}%
 \def\nextiii@##1{{\rm[}\BTcite@@##1,\BTcite@@\/{\rm]}%
  \immediate\write\auxwrite@{\string\citation{##1}}}%
 \DN@{\ifx\next(\expandafter\nextii@\else\expandafter\nextiii@\fi}%
 \FN@\next@}%
\def\BTcite@@#1,{\BTcite@@@{#1}\FN@\BTcite@@@@}
\def\BTcite@@@@{\ifx\next\BTcite@@
 \expandafter\eat@\else{\rm, }\expandafter\BTcite@@\fi}
\catcode`\~=11
\def\BTcite@@@#1{\nolabel@\cite{#1}\relax
 \DNii@##1~##2\nextii@{##1}%
 \csL@{#1}\expandafter\nextii@\Next@\nextii@\fi}
\catcode`\~=\active

\def\beginthebibliography@#1{\rm\setboxz@h{#1\ }\bibindent@\wdz@
 \bigbreak\centerline{\smc\bibliography@W}\nobreak\medskip
 \sfcode`\.=\@m\everypar{}\parindent\z@}

\newif\iffigproofing@
\def\Figureproofing{\figproofing@true}
\def\noFigureproofing{\figproofing@false}
\newif\ifHby@
\def\Hbyw#1{\global\Hby@true\hbyw\vsize{#1}}
\def\hbyw#1#2{%
 \hbox{%
  \ifHby@
  \else
   \iffigproofing@
    \setbox\z@\vbox{\hrule\width5\p@}\ht\z@\z@
    \vbox to#1{\hrule\height5\p@\width.4\p@\vfil\hrule\height5\p@\width.4\p@}%
    \kern-.4\p@\rlap{\copy\z@}\raise#1\hbox{\rlap{\copy\z@}}%
   \fi
  \fi
  \vbox to#1{\hbox to#2{}\vfil}%
  \ifHby@
  \else
   \iffigproofing@
    \vbox to#1{\hrule\height5\p@\width.4\p@\vfil\hrule\height5\p@\width.4\p@}%
    \kern-.4\p@\llap{\copy\z@}\raise#1\hbox{\llap{\boxz@}}%
   \fi
  \fi}}
\newcount\island@C
\let\island@P\empty
\let\island@Q\empty
\def\island@S#1{#1\null.}
\let\island@N\arabic
\def\island@F{\rm}
\def\island@@@P{\csname\exxx@\islandtype@ @P\endcsname}
\def\island@@@Q{\csname\exxx@\islandtype@ @Q\endcsname}
\def\island@@@S{\csname\exxx@\islandtype@ @S\endcsname}
\def\island@@@N{\csname\exxx@\islandtype@ @N\endcsname}
\def\island@@@F{\csname\exxx@\islandtype@ @F\endcsname}
\def\island@@@C{\csname island@C\islandclass@\endcsname}
\newif\ifplace@
\newif\ifisland@
\def\island{%
 \ifplace@
  \DN@{\let\islandclass@\empty\def\islandtype@{\island}\FN@\island@}%
 \else
  \long\DN@##1\endisland{\Err@{\noexpand\island must be used after some
   type of \string\...place}}%
 \fi
 \next@}
\def\island@{\ifx\next\c\let\next@\island@c\else
 \DN@{\FN@\island@@}\fi\next@}
\def\island@@{\ifcat\bgroup\noexpand\next\let\next@\island@@@\else
 \DN@{\Err@{\noexpand\island must be followed by a {prefix} for
 \string\caption's}}\fi\next@}
\newbox\islandbox@
\newcount\captioncount@
\def\island@@@#1{\def\captionprefix@{#1}\captioncount@\z@
 \global\setbox\islandbox@\vbox\bgroup}
\def\island@c\c#1{%
 \ifplace@
 \DN@{\def\islandclass@{#1}%
  \expandafter\ifx\csname island@C#1\endcsname\relax
  \expandafter\newcount@\csname island@C#1\endcsname
   \global\csname island@C#1\endcsname\z@\fi
  \FNSS@\island@c@}%
 \else
 \DN@{\edef\next@{\long\def\noexpand\next@########1\expandafter\noexpand
  \csname end\exxx@\islandtype@\endcsname{\noexpand\Err@{\noexpand\noexpand
  \expandafter\noexpand
  \islandtype@ must be used after some type of \noexpand\string
   \noexpand\...place}}}\next@\next@}%
 \fi
 \next@}
\def\island@c@{%
 \ifcat\bgroup\noexpand\next
  \let\next@\island@c@@
 \else
  \DN@{\Err@{\noexpand\island\string\c{\expandafter\string\islandclass@} must
   be followed by a {prefix} for \string\caption's}}%
 \fi\next@}
\def\island@c@@#1{\def\captionprefix@{#1}%
 \captioncount@\z@\global\setbox\islandbox@\vbox\bgroup}
\rightadd@\caption\to\nofrillslist@
\newbox\captionbox@
\newbox\Captionbox@
\def\caption{%
 \ifnum\captioncount@=\z@
  \ifnopunct@
   \DN@{\egroup\nopunct@true}%
  \else
   \let\next@\egroup
  \fi
 \else
  \let\next@\relax
 \fi
 \next@
 \advance\captioncount@\@ne
 \FN@\caption@}
\def\caption@{\ifx\next"\expandafter\caption@q\else\expandafter\caption@@\fi}
\def\caption@q"#1"{\quoted@true
 {\noexpands@
 \let\pre\island@@@P\let\post\island@@@Q
 \let\style\island@@@S\let\numstyle\island@@@N
 \Qlabel@{#1}\let\style\relax\xdef\Qlabel@@@@{#1}}%
 \finishcaption@}
\def\caption@@{\quoted@false
 \global\advance\island@@@C\@ne
 {\noexpands@
 \xdef\Thelabel@@@{\number\island@@@C}%
 \xdefThelabel@\island@@@N
 \xdef\Thelabel@@@@{\island@@@P\Thelabel@\island@@@Q}%
 \xdefThelabel@@\island@@@S
 \xdef\Thepref@{\Thelabel@@@@}}%
 \finishcaption@}
\long\def\captionformat@#1#2#3{\rm\strut#1 {\island@@@F#2} #3%
 \punct@.\strut}
\long\def\widerthanisland@#1#2#3{\test@true\setbox\z@\vbox{\hsize\maxdimen
 \noindent@@\captionformat@{#1}{#2}{#3}\par\setboxzl@}%
 \ifdim\wdz@=\z@
  \global\setbox\captionbox@\hbox{\noset@\unlabel@
   \captionformat@{#1}{#2}{#3}}%
  \ifdim\wd\captionbox@>\wd\islandbox@\else\test@false\fi
 \fi}
\long\def\captionformat@@#1#2#3{\widerthanisland@{#1}{#2}{#3}%
 \iftest@
  \global\setbox\captionbox@\vbox{\hsize\wd\islandbox@
   \vskip-\parskip\noindent@@\noset@\unlabel@
   \captionformat@{#1}{#2}{#3}\par}%
 \else
  \global\setbox\captionbox@
   \hbox to\wd\islandbox@{\hfil\box\captionbox@\hfil}%
 \fi}
\long\def\finishcaption@#1{\def\entry@{#1}%
 {\locallabel@
 \captionformat@@
  {\expandafter\ignorespaces\captionprefix@\unskip}%
  {\ifx\thelabel@@\empty\unskip\else\thelabel@@\fi}%
  {\ignorespaces#1\unskip}%
 \ifnum\captioncount@=\@ne
  \global\setbox\islandbox@\vbox{\ticwrite@\vbox{\box\islandbox@}}%
  \global\setbox\Captionbox@\vbox{\box\captionbox@}%
 \else
  \global\setbox\islandbox@\vbox{\unvbox\islandbox@\setboxzl@
   \ticwrite@\boxz@}%
  \global\setbox\Captionbox@\vbox{\unvbox\Captionbox@
   \smallskip\box\captionbox@}%
 \fi}%
 \nopunct@false\nospace@false\ignorespaces}
\def\Sixtic@{\ifx\macdef@\empty\else
 \DN@##1##2\next@{\def\macdef@{##1##2}}%
 \expandafter\next@\macdef@\next@
 \edef\next@
  {\noexpand\six@\tic@\macdef@
  \space\space\space\space\space\space\space\space\space\space\space\space
  \noexpand\six@}%
 \next@\let\macdef@\relax\fi}
\def\ticwrite@{%
 \iftoc@
  {\noexpands@\let\style\relax
  \DN@{\island}%
  \edef\next@{\write\tic@{%
   \ifnopunct@\noexpand\noexpand\noexpand\nopunct\fi
   \ifx\islandtype@\next@\noexpand\noexpand\noexpand\island
    \noexpand\string\noexpand\c{\islandclass@}{\captionprefix@}%
     {\QorThelabel@@@@}\else\noexpand\noexpand\expandafter\noexpand
     \islandtype@{\QorThelabel@@@@}}\fi}%
  \next@}%
  \expandafter\unmacro@\meaning\entry@\unmacro@
  \Sixtic@
  \write\tic@{\noexpand\Page{\number\pageno}{\page@N}{\page@P}{\page@Q}^^J}%
 \fi}
\def\Htrim@#1{%
 \ifHby@
  \dimen@\vsize
  \ifnum\captioncount@=\z@
  \else
   \advance\dimen@-\ht\Captionbox@
   \advance\dimen@-#1%
  \fi
  \global\Hby@false
  \dimen@ii\wd\islandbox@
  \global\setbox\islandbox@\vbox
   {\unvbox\islandbox@\setboxzl@
   \vbox to\z@{\vss\boxz@}\nointerlineskip\hbyw\dimen@\dimen@ii}%
  \global\Hby@true
 \fi}
\newif\ifdata@
\def\iclasstest@#1{\DN@{#1}\ifx\next@\islandclass@
 \test@true\else\test@false\fi}
\skipdef\skipi@=1
\def\endisland{\ifnum\captioncount@=\z@\expandafter\egroup\fi
 \ifdata@
 \else
  \iclasstest@{T}%
  \iftest@
   {\rm\global\skipi@-\dp\strutbox}\global\advance\skipi@\bigskipamount
   \Htrim@\skipi@
   \global\setbox\islandbox@\vbox
    {\ifnum\captioncount@=\z@\else
     \box\Captionbox@
     \nointerlineskip
     \vskip\skipi@\fi
     \box\islandbox@}%
  \else
   {\rm\global\skipi@\dp\strutbox}\global\advance\skipi@\medskipamount
   \Htrim@\skipi@
   \global\setbox\islandbox@\vbox
    {\box\islandbox@
     \ifnum\captioncount@=\z@\else
     \nointerlineskip
     \vskip\skipi@
     \box\Captionbox@
     \fi}%
  \fi
  \ifHby@
  \else
   \dimen@\ht\islandbox@\advance\dimen@\dp\islandbox@
   \ifdim\dimen@>\vsize
    \DN@{\island}%
    \Err@{%
     \ifx\islandtype@\next@\noexpand\island\else
      \expandafter\noexpand\islandtype@\fi
     \ifnum\captioncount@=\z@\else
       with \noexpand\caption\fi
      is larger than page}%
     \ht\islandbox@=\vsize
   \fi
  \fi
 \fi
 \global\Hby@false\island@true}
\def\newisland#1\c#2#3{\define#1{}%
 \iftoc@\immediate\write\tic@{\noexpand\newisland\noexpand#1%
  \string\c{#2}{#3}^^J}\fi
 \expandafter\def\csname\exstring@#1@S\endcsname{\island@S}%
 \expandafter\def\csname\exstring@#1@N\endcsname{\island@N}%
 \expandafter\def\csname\exstring@#1@P\endcsname{\island@P}%
 \expandafter\def\csname\exstring@#1@Q\endcsname{\island@Q}%
 \expandafter\def\csname\exstring@#1@F\endcsname{\island@F}%
 \expandafter\def\csname end\exstring@#1\endcsname{\endisland}%
 \expandafter
 \ifx\csname island@C#2\endcsname\relax
  \expandafter\newcount@\csname island@C#2\endcsname
  \global\csname island@C#2\endcsname\z@
 \fi
 \edef\next@{\noexpand\expandafter\noexpand\let\noexpand
  \csname\exstring@#1@C\noexpand\endcsname
  \csname island@C#2\endcsname}%
 \next@
 \def#1{\def\islandtype@{#1}\island@c\c{#2}{#3}}}
\newisland\Figure\c{F}{Figure}
\newisland\Table\c{T}{Table}
\newbox\islandboxi
\newbox\islandboxii
\newbox\islandboxiii
\newbox\captionboxi
\newbox\captionboxii
\newbox\captionboxiii
\long\def\islandpairdata#1#2{{\data@true
 \place@true
 #1%
 \global\setbox\islandboxi\box\islandbox@
 \global\setbox\captionboxi\box\Captionbox@
 #2%
 \global\setbox\islandboxii\box\islandbox@
 \global\setbox\captionboxii\box\Captionbox@
 }}
\long\def\islandpairbox#1#2{\islandpairdata{#1}{#2}%
 \dimen@\ht\captionboxi
 \ifdim\ht\captionboxii>\dimen@\dimen@\ht\captionboxii\fi
 \ifdim\dimen@>\z@
  \ifdim\ht\captionboxi<\dimen@
   \global\setbox\captionboxi\vbox to\dimen@{\unvbox\captionboxi\vfil}\fi
  \ifdim\ht\captionboxii<\dimen@
   \global\setbox\captionboxii\vbox to\dimen@{\unvbox\captionboxii\vfil}\fi
 \fi
 \global\setbox\islandbox@\vbox
 {\hbox to\hsize{\hfil\box\islandboxi\hfil\box\islandboxii\hfil}%
 \ifdim\dimen@>\z@\nointerlineskip
 {\rm\global\skipi@\dp\strutbox}\global\advance\skipi@\medskipamount
  \vskip\skipi@
  \hbox to\hsize{\hfil\box\captionboxi\hfil\box\captionboxii\hfil}\fi}}	
\long\def\islandpairboxa#1#2{\islandpairdata{#1}{#2}%
 \dimen@\ht\captionboxi
 \ifdim\ht\captionboxii>\dimen@\dimen@\ht\captionboxii\fi
 \ifdim\dimen@>\z@
  \ifdim\ht\captionboxi<\dimen@
   \global\setbox\captionboxi\vbox to\dimen@{\vfil\unvbox\captionboxi}\fi
  \ifdim\ht\captionboxii<\dimen@
   \global\setbox\captionboxii\vbox to\dimen@{\vfil\unvbox\captionboxii}\fi
 \fi
 \dimen@ii\ht\islandboxi
 \ifdim\ht\islandboxii>\dimen@ii \dimen@ii\ht\islandboxii\fi
 \ifdim\dimen@ii>\z@
  \ifdim\ht\islandboxi<\dimen@ii
   \global\setbox\islandboxi\vbox to\dimen@ii{\box\islandboxi\vfil}\fi
  \ifdim\ht\islandboxii<\dimen@ii
   \global\setbox\islandboxii\vbox to\dimen@ii{\box\islandboxii\vfil}\fi
 \fi
 \global\setbox\islandbox@\vbox{\ifdim\dimen@>\z@
  \hbox to\hsize{\hfil\box\captionboxi\hfil\box\captionboxii\hfil}%
  \nointerlineskip{\rm\global\skipi@-\dp\strutbox}%
  \global\advance\skipi@\bigskipamount\vskip\skipi@\fi
  \hbox to\hsize{\hfil\box\islandboxi\hfil\box\islandboxii\hfil}}}
\long\def\islandtripledata#1#2#3{{\data@true\place@true
 #1%
 \global\setbox\islandboxi\box\islandbox@
 \global\setbox\captionboxi\box\Captionbox@
 #2%
 \global\setbox\islandboxii\box\islandbox@
 \global\setbox\captionboxii\box\Captionbox@
 #3%
 \global\setbox\islandboxiii\box\islandbox@
 \global\setbox\captionboxiii\box\Captionbox@
 }}
\long\def\islandtriplebox#1#2#3{\islandtripledata{#1}{#2}{#3}%
 \dimen@\ht\captionboxi
 \ifdim\ht\captionboxii>\dimen@ \dimen@\ht\captionboxii\fi
 \ifdim\ht\captionboxiii>\dimen@ \dimen@\ht\captionboxiii\fi
 \ifdim\dimen@>\z@
  \ifdim\ht\captionboxi<\dimen@
   \global\setbox\captionboxi\vbox to\dimen@{\unvbox\captionboxi\vfil}\fi
  \ifdim\ht\captionboxii<\dimen@
   \global\setbox\captionboxii\vbox to\dimen@{\unvbox\captionboxii\vfil}\fi
  \ifdim\ht\captionboxiii<\dimen@
   \global\setbox\captionboxiii\vbox to\dimen@{\unvbox\captionboxiii\vfil}\fi
 \fi
 \global\setbox\islandbox@\vbox
  {\hbox to\hsize{\hfil\box\islandboxi\hfil\box\islandboxii\hfil
   \box\islandboxiii\hfil}%
 \ifdim\dimen@>\z@\nointerlineskip
  {\rm\global\skipi@\dp\strutbox}\global\advance\skipi@\medskipamount
  \vskip\skipi@
  \hbox to\hsize{\hfil\box\captionboxi\hfil\box\captionboxii\hfil
   \box\captionboxiii\hfil}\fi}}
\def\islandtripleboxa#1#2#3{\islandtripledata{#1}{#2}{#3}%
 \dimen@\ht\captionboxi
 \ifdim\ht\captionboxii>\dimen@ \dimen@\ht\captionboxii\fi
 \ifdim\ht\captionboxiii>\dimen@ \dimen@\ht\captionboxiii\fi
 \ifdim\dimen@>\z@
  \ifdim\ht\captionboxi<\dimen@
   \global\setbox\captionboxi\vbox to\dimen@{\vfil\unvbox\captionboxi}\fi
  \ifdim\ht\captionboxii<\dimen@
   \global\setbox\captionboxii\vbox to\dimen@{\vfil\unvbox\captionboxii}\fi
  \ifdim\ht\captionboxiii<\dimen@
   \global\setbox\captionboxiii\vbox to\dimen@{\vfil\unvbox\captionboxiii}\fi
 \fi
 \dimen@ii\ht\islandboxi
 \ifdim\ht\islandboxii>\dimen@ii \dimen@ii\ht\islandboxii\fi
 \ifdim\ht\islandboxiii>\dimen@ii \dimen@ii\ht\islandboxiii\fi
 \ifdim\dimen@ii>\z@
  \ifdim\ht\islandboxi<\dimen@ii
   \global\setbox\islandboxi\vbox to\dimen@ii{\box\islandboxi\vfil}\fi
  \ifdim\ht\islandboxii<\dimen@ii
   \global\setbox\islandboxii\vbox to\dimen@ii{\box\islandboxii\vfil}\fi
  \ifdim\ht\islandboxiii<\dimen@ii
   \global\setbox\islandboxiii\vbox to\dimen@ii{\box\islandboxiii\vfil}\fi
 \fi
 \global\setbox\islandbox@\vbox
  {\ifdim\dimen@>\z@
  \hbox to\hsize{\hfil\box\captionboxi\hfil\box\captionboxii\hfil
   \box\captionboxiii\hfil}%
  \nointerlineskip{\rm\global\skipi@-\dp\strutbox}%
  \global\advance\skipi@\bigskipamount\vskip\skipi@\fi
  \hbox to\hsize{\hfil\box\islandboxi\hfil\box\islandboxii\hfil
   \box\islandboxiii\hfil}}}
\def\Figurepair#1\and#2\endFigurepair{\island@true
 \islandpairbox{\Figure#1\endFigure}{\Figure#2\endFigure}}
\def\Figuretriple#1\and#2\and#3\endFiguretriple{\island@true
 \islandtriplebox{\Figure#1\endFigure}{\Figure#2\endFigure}%
  {\Figure#3\endFigure}}
\def\Tablepair#1\and#2\endTablepair{\island@true
 \islandpairboxa{\Table#1\endTable}{\Table#2\endTable}}
\def\Tabletriple#1\and#2\and#3\endTabletriple{\island@true
 \islandtripleboxa{\Table#1\endTable}{\Table#2\endTable}%
 {\Table#3\endTable}}
\def\place#1{\place@true\island@false
 #1%
 \ifisland@
  \box\islandbox@
 \else
  \Err@{Whoa ... there's no \string\Figure, \string\Table,
   etc., here}%
 \fi
 \place@false}
\newskip\belowtopfigskip
\belowtopfigskip 15\p@ plus 5\p@ minus5\p@
\newskip\abovebotfigskip
\abovebotfigskip 18\p@ plus 6\p@ minus6\p@
\newdimen\minpagesize
\minpagesize 5pc
\dimen@\belowtopfigskip
\advance\dimen@-\abovebotfigskip
\skip\topins\dimen@
\dimen\topins\z@
\newcount\topinscount@
\newbox\topinsdims@
\def\storedim@{\global\setbox\topinsdims@
 \vbox{\hbox to\dimen@{}\unvbox\topinsdims@}}
\def\advancedimtopins@{%
 \ifnum\pageno=\@ne
 \else
   \advance\dimen@\dimen\topins
   \global\dimen\topins\dimen@
 \fi}
\newcount\flipcount@
\def\fliptopins@{%
 \global\flipcount@\z@
 \ifvoid\topins\else
 \setbox\z@\vbox
  {\vskip\p@
   \unvbox\topins
   \global\setbox\topins\vbox{}%
   \loop
    \test@false
    \ifdim\lastskip=\z@\unskip
     \ifdim\lastskip=\z@
      \test@true\fi\fi
    \iftest@
    \global\advance\flipcount@\@ne
    \setboxzl@
    \global\setbox\topins\vbox{\unvbox\topins\boxz@}%
    \unpenalty
   \repeat}\fi}
\newif\ifPar@
\newcount\Parcount@
\newbox\Parbox@
\expandafter\newbox\csname Parfigbox1\endcsname
\expandafter\newbox\csname Parfigbox2\endcsname
\expandafter\newbox\csname Parfigbox3\endcsname
\expandafter\newbox\csname Parfigbox4\endcsname
\expandafter\newbox\csname Parfigbox5\endcsname
\expandafter\newdimen\csname Parprev1\endcsname
\expandafter\newdimen\csname Parprev2\endcsname
\expandafter\newdimen\csname Parprev3\endcsname
\expandafter\newdimen\csname Parprev4\endcsname
\expandafter\newdimen\csname Parprev5\endcsname
\expandafter\newdimen\csname Parprev6\endcsname
\def\Par{\par\global\csname Parprev1\endcsname\prevdepth
 \global\Parcount@\@ne
 \global\Par@true\global\let\Parlist@\empty
 \global\setbox\Parbox@\vbox\bgroup\break}
\def\place@#1#2{%
 \ifisland@
  \ifhmode
   \ifPar@
    \ifnum\Parcount@>5
     \Err@{Only 5 \string\place's allowed per
      \string\Par...\noexpand\endPar paragraph}%
    \else
     \expandafter\expandafter\expandafter
      \global\expandafter\setbox
       \csname Parfigbox\number\Parcount@\endcsname\box\islandbox@
     \global\advance\Parcount@\@ne
     \xdef\Parlist@{\Parlist@#1}%
    \fi
   \else
    \vadjust{#2}%
   \fi
  \else
   #2%
  \fi
 \else
  \Err@{Whoa ... there's no \string\Figure,
   \string\Table, etc., here}%
 \fi
 \place@false}
\long\def\Aplace#1{\prevanish@
 \place@true\island@false
 #1%
 \place@ a\Aplace@
 \postvanish@}
\long\def\AAplace#1{\prevanish@\place@true\island@false
 #1%
 \place@ A\AAplace@
 \postvanish@}
\newif\ifAA@
\def\AAplace@{\AA@true\Aplace@\AA@false}
\let\AAlist@\empty
\def\Aplace@{\allowbreak
 \dimen@=\ht\islandbox@
 \advance\dimen@\abovebotfigskip
 \ht\islandbox@\dimen@
 \advance\dimen@\dp\islandbox@
 \storedim@
 \ifAA@
  \xdef\AAlist@{\AAlist@1}%
  \advancedimtopins@
 \else
  \xdef\AAlist@{\AAlist@0}%
  \ifnum\topinscount@>\@ne\else\advancedimtopins@\fi
 \fi
 \insert\topins{\penalty\z@\splittopskip\z@\floatingpenalty\z@
  \box\islandbox@}%
 \global\advance\topinscount@\@ne}
\long\def\Bplace#1{\prevanish@\place@true\island@false
 #1%
 \place@ b\Bplace@
 \postvanish@}
\def\Bplace@{\allowbreak
 \ifnum\topinscount@=\z@
  \setbox\z@\vbox{\vbox to-\belowtopfigskip{}}%
  \dimen@-\skip\topins
  \ht\z@\dimen@
  \storedim@
  \advancedimtopins@
  \insert\topins{\boxz@}%
  \global\advance\topinscount@\@ne
  \xdef\AAlist@{\AAlist@0}%
 \fi
 \dimen@\ht\islandbox@
 \advance\dimen@\abovebotfigskip
 \ht\islandbox@\dimen@
 \advance\dimen@\dp\islandbox@
 \storedim@
 \xdef\AAlist@{\AAlist@0}%
 \ifnum\topinscount@>\@ne\else\advancedimtopins@\fi
 \insert\topins{\penalty\z@\splittopskip\z@
  \floatingpenalty\z@
  \box\islandbox@}%
 \global\advance\topinscount@\@ne}
\def\breakisland@{\global\setbox\@ne\lastbox\global\skipi@\lastskip\unskip
 \global\setbox\thr@@\lastbox}%
\def\printisland@{\centerline{\box\thr@@}\nobreak\nointerlineskip
 \vskip\skipi@
 \ifdim\ht\@ne<\z@\box\@ne\else\centerline{\box\@ne}\fi}
\def\bottomfigs@{%
 \count@\@ne
 \loop
  \ifnum\count@<\flipcount@
  \nointerlineskip
  \vskip\abovebotfigskip
  \global\setbox\topins\vbox{\unvbox\topins\setboxzl@
   \unvbox\z@
   \breakisland@}%
  \printisland@
  \advance\count@\@ne
  \repeat}
\def\resetdimtopins@{%
 \global\advance\topinscount@-\flipcount@
 \global\setbox\topinsdims@\vbox
  {\unvbox\topinsdims@
   \count@\z@
   \DN@##1##2\next@{\gdef\AAlist@{##2}}%
   \loop
    \ifnum\count@<\flipcount@\setboxzl@
    \expandafter\next@\AAlist@\next@
    \advance\count@\@ne
    \repeat
   \dimen@\z@
   \count@\z@
   \setbox\tw@\vbox{}%
   \edef\nextiii@{\AAlist@}%
   \DN@##1##2\next@{\DNii@{##1}\def\nextiii@{##2}}%
   \loop
    \test@false
    \ifnum\count@<\topinscount@
    \expandafter\next@\nextiii@\next@
     \ifnum\count@<\tw@
      \test@true
     \else
      \if\nextii@ 1\test@true\fi
     \fi
    \fi
    \iftest@
     \setboxzl@
     \advance\dimen@\wdz@
     \setbox\tw@\vbox{\boxz@\unvbox\tw@}%
     \advance\count@\@ne
    \repeat
    \unvbox\tw@
    \global\dimen\topins\dimen@}}
\def\Place@#1#2{%
 \ifisland@
  \ifhmode
   \ifPar@
    \ifnum\Parcount@>5
     \Err@{Only 5 \string\place's allowed per
       \string\Par...\noexpand\endPar paragraph}%
    \else
     \expandafter\expandafter\expandafter\global\expandafter\setbox
      \csname Parfigbox\number\Parcount@\endcsname\box\islandbox@
     \global\advance\Parcount@\@ne
     \xdef\Parlist@{\Parlist@#1}%
     \vadjust{\break}%
    \fi
   \else
    \Err@{\noexpand#2allowed only in a \string\Par...\noexpand\endPar
     paragraph}%
   \fi
  \else
   #2%
  \fi
 \else
  \Err@{Who ... there's no \string\Figure, \string\Table,
   etc., here}%
 \fi
 \place@false}
\newif\ifC@
\newdimen\Cdim@
\long\def\Cplace#1{\prevanish@\place@true\island@false
 #1%
 \Place@ c\Cplace@
 \postvanish@}
\def\Cplace@{\allowbreak
 \ifnum\topinscount@>\z@\else
  \global\C@true\global\Cdim@\pagetotal\fi
 \Aplace@}
\long\def\Mplace#1{\prevanish@\place@true\island@false
 #1%
 \Place@ m\Mplace@
 \postvanish@}
\long\def\MXplace#1{\prevanish@\place@true\island@false
 #1%
 \Place@ M\MXplace@
 \postvanish@}
\newif\ifMX@
\def\MXplace@{\MX@true\Mplace@\MX@false}
\def\Mplace@{\allowbreak
 \dimen@\ht\islandbox@\advance\dimen@\dp\islandbox@
 \ifdim\pagetotal=\z@\else
  \ifdim\lastskip<\abovebotfigskip\advance\dimen@\abovebotfigskip
  \advance\dimen@-\lastskip\fi
 \fi
 \advance\dimen@\pagetotal
 \ifdim\dimen@>\pagegoal
  \Aplace@
 \else
  \nointerlineskip
  \ifdim\lastskip<\abovebotfigskip\removelastskip\vskip\abovebotfigskip\fi
  \setbox\z@\vbox{\unvbox\islandbox@
   \breakisland@}%
  \printisland@
  \ifnum\topinscount@=\z@
   \setbox\z@\vbox{\vbox to-\belowtopfigskip{}}%
   \dimen@-\skip\topins
   \ht\z@\dimen@
   \storedim@
   \advancedimtopins@
   \insert\topins{\boxz@}%
   \global\advance\topinscount@\@ne
   \xdef\AAlist@{\AAlist@0}%
  \fi
  \ifMX@
   \ifnum\topinscount@=\@ne
    \setbox\z@\vbox{\vbox to-\abovebotfigskip{}}%
    \ht\z@\z@
    \dimen@\z@
    \storedim@
    \advancedimtopins@
    \insert\topins{\boxz@}%
    \global\advance\topinscount@\@ne
    \xdef\AAlist@{\AAlist@0}%
   \fi
  \fi
  \nointerlineskip
  \vskip\belowtopfigskip
 \fi}
\expandafter\newbox\csname Parbox1\endcsname
\expandafter\newbox\csname Parbox2\endcsname
\expandafter\newbox\csname Parbox3\endcsname
\expandafter\newbox\csname Parbox4\endcsname
\expandafter\newbox\csname Parbox5\endcsname
\def\endPar{\egroup
 \count@\@ne
 {\vbadness\@M\vfuzz\maxdimen\splitmaxdepth\maxdimen\splittopskip\ht\strutbox
 \setbox\z@\vsplit\Parbox@ to\ht\Parbox@
 \loop
  \ifnum\count@<\Parcount@
  \expandafter\expandafter\expandafter\global\expandafter\setbox
   \csname Parbox\number\count@\endcsname\vsplit\Parbox@ to\ht\Parbox@
  \count@@\count@\advance\count@@\@ne
  \global\csname Parprev\number\count@@\endcsname
   \dp\csname Parbox\number\count@\endcsname
  \advance\count@\@ne
  \repeat}%
 \vskip\parskip
 \count@\@ne
 \def\nextv@##1##2\nextv@{\DN@{##1}\gdef\Parlist@{##2}}%
 \loop
  \ifnum\count@<\Parcount@
   \dimen@\csname Parprev\number\count@\endcsname
   \advance\dimen@\ht\strutbox
   \ifdim\dimen@<\baselineskip
    \advance\dimen@-\baselineskip\vskip-\dimen@
   \else
    \vskip\lineskip
   \fi
   \unvbox\csname Parbox\number\count@\endcsname
   \global\setbox\islandbox@\box\csname Parfigbox\number\count@\endcsname
   \expandafter\nextv@\Parlist@\nextv@
   \if a\next@\Aplace@\else
   \if A\next@\AAplace@\else
   \if b\next@\Bplace@\else
   \if c\next@\Cplace@\else
   \if m\next@\Mplace@\else
   \if M\next@\MXplace@\fi\fi\fi\fi\fi\fi
  \advance\count@\@ne
  \repeat
 \global\Par@false
 \ifvoid\Parbox@
  \prevdepth\csname Parprev\number\count@\endcsname
 \else
  \dimen@\csname Parprev\number\count@\endcsname\advance\dimen@\ht\strutbox
  \ifdim\dimen@<\baselineskip
    \advance\dimen@-\baselineskip\vskip-\dimen@
  \else
    \vskip\lineskip
  \fi
  \dimen@\dp\Parbox@
  \unvbox\Parbox@
  \prevdepth\dimen@
 \fi}
\def\folio{{\page@F\page@S{\page@P\page@N{\number\page@C}\page@Q}}}
\def\advancepageno{\global\advance\pageno\@ne}
\newif\ifspecialsplit@
\newbox\outbox@
\let\shipout@\shipout
\def\plainoutput{\specialsplit@false\ifvoid\topins\else\ifdim\ht\topins=\z@
 \specialsplit@true\advance\minpagesize-\skip\topins\fi\fi
 \fliptopins@
 \setbox\outbox@\vbox{\makeheadline\pagebody\makefootline}%
 {\noexpands@\let\style\relax
 \shipout@\box\outbox@}%
 \advancepageno
 \resetdimtopins@
 \ifvoid\@cclv\else\unvbox\@cclv\penalty\outputpenalty\fi
 \ifnum\outputpenalty>-\@MM\else\dosupereject\fi}
\def\pagebody{\vbox to\vsize{\boxmaxdepth\maxdepth
 \ifvoid\margin@\else
 \rlap{\kern\hsize\vbox to\z@{\kern4\p@\box\margin@\vss}}\fi
 \pagecontents}}
\newif\ifonlytop@
\def\pagecontents{%
 \onlytop@false
 \ifdim\ht\@cclv<\minpagesize\ifnum\flipcount@<\tw@\ifvoid\footins
  \onlytop@true\fi\fi\fi
 \test@false
 \ifC@
  \ifnum\flipcount@=\@ne
   \global\multiply\Cdim@\tw@
   \ifdim\Cdim@>\ht\@cclv
    \test@true
   \fi
  \fi
 \fi
 \global\C@false
 \iftest@
  \dimen@\ht\@cclv
  \advance\dimen@\skip\topins
  {\vfuzz\maxdimen\vbadness\@M
  \splitmaxdepth\maxdepth\splittopskip\topskip
  \setbox\z@\vsplit\@cclv to\dimen@
  \unvbox\z@}%
  \global\setbox\topins\vbox{\unvbox\topins
   \global\setbox\@ne\lastbox}%
  \setbox\z@\vbox{\unvbox\@ne
   \breakisland@}%
  \nointerlineskip
  \vskip\abovebotfigskip
  \printisland@
 \else
  \ifnum\flipcount@>\z@
   \global\setbox\topins\vbox{\unvbox\topins\global\setbox\@ne\lastbox}%
   \setbox\z@\vbox{\unvbox\@ne
    \breakisland@}%
   \printisland@
   \ifonlytop@\kern-\prevdepth\vfill\else\vskip\belowtopfigskip\fi
  \fi
 \fi
 \ifdim\ht\@cclv<\minpagesize
  \ifonlytop@\else\vfill\fi
 \else
  \ifspecialsplit@
   {\vfuzz\maxdimen\vbadness\@M
   \splitmaxdepth\maxdepth\splittopskip\topskip
   \dimen@ii\ht\@cclv \advance\dimen@ii\skip\topins
   \setbox\z@\vsplit\@cclv to\dimen@ii
   \unvbox\z@}%
  \else
   \unvbox\@cclv
  \fi
 \fi
 \bottomfigs@
 \ifvoid\footins\else\vskip\skip\footins\footnoterule\unvbox\footins\fi}
\newread\readdata@
\def\readthedata@#1{\expandafter
 \ifx\csname#1@D\endcsname\relax
  \immediate\openin\readdata@=#1.dat
  \ifeof\readdata@
   \Err@{No file #1.dat}%
  \else
   {\endlinechar\m@ne\gdef\Next@{}%
   \DNii@##1 ##2 ##3pt{\global\data@ht##1\global\data@dp##2%
    \global\data@wd##3pt}%
   \loop
    \ifeof\readdata@
    \else
    \read\readdata@ to\next@
    \ifx\next@\empty\else
     \edef\next@{\expandafter\nextii@\next@}%
     \expandafter\rightadd@\next@\to\Next@
    \fi
    \repeat}%
   \immediate\closein\readdata@
   \expandafter\expandafter\expandafter\global\expandafter
    \let\csname#1@D\endcsname\Next@\global\let\Next@\relax
  \fi
 \fi}
\newdimen\data@ht
\newdimen\data@dp
\newdimen\data@wd
\newif\ifgetdata@
\def\getdata@#1#2{\global\getdata@true\count@#2\relax
 {\let\\\or\xdef\Next@{\ifcase\number\count@#1\else
 \global\noexpand\getdata@false\fi}}\Next@}
\def\paste#1#2{\readthedata@{#1}%
 \getdata@{\csname#1@D\endcsname}{#2}%
 \ifgetdata@
 \dimen@\data@ht \advance\dimen@\data@dp
  \hbox{\special{dvipaste: #1 #2}%
   \lower\data@dp\vbox to\dimen@{\hbox to\data@wd{}\vfil}}%
 \else
  {\lccode`\Z=`\#\lccode`\N=`\N\lccode`\F=`\F%
   \lowercase{\Err@{No data for File [#1], Z#2}}}%
 \fi}
\newdimen\httable
\newdimen\dptable
\newdimen\wdtable
\def\measuretable#1#2{\readthedata@{#1}%
 \getdata@{\csname#1@D\endcsname}{#2}%
 \ifgetdata@
  \httable\data@ht \dptable\data@dp \wdtable\data@wd
 \else
  {\lccode`\Z=`\#\lccode`\N=`\N\lccode`\F=`\F%
  \lowercase{\Err@{No data for File [#1], Z#2}}}%
 \fi}
\def\East#1#2{\setboxz@h{$\m@th\ssize\;{#1}\;\;$}%
 \setbox\tw@\hbox{$\m@th\ssize\;{#2}\;\;$}\setbox4=\hbox{$\m@th#2$}%
 \dimen@\minaw@
 \ifdim\wdz@>\dimen@\dimen@\wdz@\fi\ifdim\wd\tw@>\dimen@\dimen@\wd\tw@\fi
 \ifdim\wd4 >\z@
  \mathrel{\mathop{\hbox to\dimen@{\rightarrowfill}}\limits^{#1}_{#2}}%
 \else
  \mathrel{\mathop{\hbox to\dimen@{\rightarrowfill}}\limits^{#1}}%
 \fi}
\def\West#1#2{\setboxz@h{$\m@th\ssize\;\;{#1}\;$}%
 \setbox\tw@\hbox{$\m@th\ssize\;\;{#2}\;$}\setbox4=\hbox{$\m@th#2$}%
 \dimen@\minaw@
 \ifdim\wdz@>\dimen@\dimen@\wdz@\fi\ifdim\wd\tw@>\dimen@\dimen@\wd\tw@\fi
 \ifdim\wd4 >\z@
  \mathrel{\mathop{\hbox to\dimen@{\leftarrowfill}}\limits^{#1}_{#2}}%
 \else
  \mathrel{\mathop{\hbox to\dimen@{\leftarrowfill}}\limits^{#1}}%
 \fi}
\font\arrow@i=lams1
\font\arrow@ii=lams2
\font\arrow@iii=lams3
\font\arrow@iv=lams4
\font\arrow@v=lams5
\newdimen\standardcgap
\standardcgap40\p@
\newdimen\hunit
\hunit\tw@\p@
\newdimen\standardrgap
\standardrgap32\p@
\newdimen\vunit
\vunit1.6\p@
\def\Cgaps#1{\RIfM@
 \standardcgap#1\standardcgap\relax\hunit#1\hunit\relax
 \else\nonmatherr@\Cgaps\fi}
\def\Rgaps#1{\RIfM@
 \standardrgap#1\standardrgap\relax\vunit#1\vunit\relax
 \else\nonmatherr@\Rgaps\fi}
\newdimen\getdim@
\def\getcgap@#1{\ifcase#1\or\getdim@\z@\else\getdim@\standardcgap\fi}
\def\getrgap@#1{\ifcase#1\getdim@\z@\else\getdim@\standardrgap\fi}
\def\cgaps{\RIfM@\expandafter\cgaps@\else\expandafter\nonmatherr@
 \expandafter\cgaps\fi}
\def\cgaps@{\ifnum\catcode`\;=\active\expandafter\cgapsA@\else
 \expandafter\cgapsO@\fi}
\def\cgapsO@#1{\toks@{\ifcase\i@\or\getdim@=\z@}%
 \gaps@@\standardcgap#1;\gaps@@\gaps@@
 \edef\next@{\the\toks@\noexpand\else\noexpand\getdim@\noexpand\standardcgap
  \noexpand\fi}%
 \toks@=\expandafter{\next@}%
 \edef\getcgap@##1{\i@##1\relax\the\toks@}\toks@{}}
{\catcode`\;=\active
 \gdef\cgapsA@#1{\toks@{\ifcase\i@\or\getdim@=\z@}%
 \gaps@@\standardcgap#1;\gaps@@\gaps@@
 \edef\next@{\the\toks@\noexpand\else\noexpand\getdim@\noexpand\standardcgap
  \noexpand\fi}%
 \toks@=\expandafter{\next@}%
 \edef\getcgap@##1{\i@##1\relax\the\toks@}\toks@{}}
}
\def\Gaps@@{\gaps@@}
\def\gaps@@#1#2;#3{\mgaps@#1#2\mgaps@
 \edef\next@{\the\toks@\noexpand\or\noexpand\getdim@
  \noexpand#1\the\mgapstoks@@}%
 \toks@\expandafter{\next@}%
 \DN@{#3}%
 \ifx\next@\Gaps@@\def\next@##1\gaps@@{}\else
  \def\next@{\gaps@@#1#3}\fi\next@}
{\catcode`\;=\active
 \gdef\rgaps#1{\RIfM@{\ifnum\catcode`\;=\active\def;{\string;}\fi
   \xdef\Next@{\noexpand\rgaps@{#1}}}%
  \Next@\edef\getrgap@##1{\i@##1\relax\the\toks@}\toks@{}\else
  \nonmatherr@\rgaps\fi}
}
\def\rgaps@#1{\toks@{\ifcase\i@\getdim@=\z@}%
 \gaps@@\standardrgap#1;\gaps@@\gaps@@
 \edef\next@{\the\toks@\noexpand\else\noexpand\getdim@\noexpand\standardrgap
  \noexpand\fi}%
 \toks@=\expandafter{\next@}}
\newbox\ZER@
\def\mgaps@#1{\let\mgapsnext@#1\FNSS@\mgaps@@}
\def\mgaps@@{\ifx\next\w\expandafter\mgaps@@@\else
 \expandafter\mgaps@@@@\fi}
\newtoks\mgapstoks@@
\def\mgaps@@@@#1\mgaps@{\getdim@\mgapsnext@\getdim@#1\getdim@
 \edef\next@{\noexpand\getdim@\the\getdim@}%
 \mgapstoks@@\expandafter{\next@}}
\def\mgaps@@@\w#1#2\mgaps@{\mgaps@@@@#2\mgaps@
 \setbox\ZER@\hbox{$\m@th\hskip15\p@\tsize@#1$}%
 \dimen@\wd\ZER@
 \ifdim\dimen@>\getdim@\getdim@\dimen@\fi
 \edef\next@{\noexpand\getdim@\the\getdim@}%
 \mgapstoks@@\expandafter{\next@}}
\def\changewidth#1#2{\setbox\ZER@{$\m@th#2}%
 \hbox to\wd\ZER@{\hss$\m@th#1$\hss}}
\atdef@({\FN@\ARROW@}
\def\ARROW@{\ifx\next)\let\next@\OPTIONS@\else
 \DN@{\csname\string @(\endcsname}\fi\next@}
\newif\ifoptions@
\def\OPTIONS@){\ifoptions@\let\next@\relax\else
 \DN@{\global\options@true\begingroup\optioncodes@}\fi\next@}
\newif\ifN@
\newif\ifE@
\newif\ifNESW@
\newif\ifH@
\newif\ifV@
\newif\ifHshort@
\expandafter\def\csname\string @(\endcsname #1,#2){%
 \ifoptions@\expandafter\endgroup\fi
 \N@false\E@false\H@false\V@false\Hshort@false
 \ifnum#1>\z@\E@true\fi
 \ifnum#1=\z@\V@true\global\tX@false\global\tY@false\global\a@false\fi
 \ifnum#2>\z@\N@true\fi
 \ifnum#2=\z@\H@true\global\tX@false\global\tY@false\global\a@false
  \ifshort@\Hshort@true\fi\fi
 \NESW@false
 \ifN@\ifE@\NESW@true\fi\else\ifE@\else\NESW@true\fi\fi
 \arrow@{#1}{#2}%
 \global\options@false
 \global\scount@\z@\global\tcount@\z@\global\arrcount@\z@
 \global\s@false\global\sxdimen@\z@\global\sydimen@\z@
 \global\tX@false\global\tXdimen@i\z@\global\tXdimen@ii\z@
 \global\tY@false\global\tYdimen@i\z@\global\tYdimen@ii\z@
 \global\a@false\global\exacount@\z@
 \global\x@false\global\xdimen@\z@
 \global\X@false\global\Xdimen@\z@
 \global\y@false\global\ydimen@\z@
 \global\Y@false\global\Ydimen@\z@
 \global\p@false\global\pdimen@\z@
 \global\label@ifalse\global\label@iifalse
 \global\dl@ifalse\global\ldimen@i\z@
 \global\dl@iifalse\global\ldimen@ii\z@
 \global\short@false\global\unshort@false}
\newif\iflabel@i
\newif\iflabel@ii
\newcount\scount@
\newcount\tcount@
\newcount\arrcount@
\newif\ifs@
\newdimen\sxdimen@
\newdimen\sydimen@
\newif\iftX@
\newdimen\tXdimen@i
\newdimen\tXdimen@ii
\newif\iftY@
\newdimen\tYdimen@i
\newdimen\tYdimen@ii
\newif\ifa@
\newcount\exacount@
\newif\ifx@
\newdimen\xdimen@
\newif\ifX@
\newdimen\Xdimen@
\newif\ify@
\newdimen\ydimen@
\newif\ifY@
\newdimen\Ydimen@
\newif\ifp@
\newdimen\pdimen@
\newif\ifdl@i
\newif\ifdl@ii
\newdimen\ldimen@i
\newdimen\ldimen@ii
\newif\ifshort@
\newif\ifunshort@
\def\zero@#1{\ifnum\scount@=\z@
 \if#1e\global\scount@\m@ne\else
 \if#1t\global\scount@\tw@\else
 \if#1h\global\scount@\thr@@\else
 \if#1'\global\scount@6 \else
 \if#1`\global\scount@7 \else
 \if#1(\global\scount@8 \else
 \if#1)\global\scount@9 \else
 \if#1s\global\scount@12 \else
 \if#1H\global\scount@13 \else
 \Err@{\Invalid@@ option \string\0}\fi\fi\fi\fi\fi\fi\fi\fi\fi
 \fi}
\def\one@#1{\ifnum\tcount@=\z@
 \if#1e\global\tcount@\m@ne\else
 \if#1h\global\tcount@\tw@\else
 \if#1t\global\tcount@\thr@@\else
 \if#1'\global\tcount@4 \else
 \if#1`\global\tcount@5 \else
 \if#1(\global\tcount@\ten@ \else
 \if#1)\global\tcount@11 \else
 \if#1s\global\tcount@12 \else
 \if#1H\global\tcount@13 \else
 \Err@{\Invalid@@ option \string\1}\fi\fi\fi\fi\fi\fi\fi\fi\fi
 \fi}
\def\a@#1{\ifnum\arrcount@=\z@
 \if#10\global\arrcount@\m@ne\else
 \if#1+\global\arrcount@\@ne\else
 \if#1-\global\arrcount@\tw@\else
 \if#1=\global\arrcount@\thr@@\else
 \Err@{\Invalid@@ option \string\a}\fi\fi\fi\fi
 \fi}
\def\ds@{\ifnum\catcode`\;=\active\expandafter\dsA@\else
 \expandafter\dsO@\fi}
\def\dsO@(#1;#2){\ds@@{#1}{#2}}
\def\ds@@#1#2{\ifs@\else
 \global\s@true
 \global\sxdimen@\hunit\global\sxdimen@#1\sxdimen@\relax
 \global\sydimen@\vunit\global\sydimen@#2\sydimen@\relax
 \fi}
\def\dtX@{\ifnum\catcode`\;=\active\expandafter\dtXA@\else
 \expandafter\dtXO@\fi}
\def\dtXO@(#1;#2){\dtX@@{#1}{#2}}
\def\dtX@@#1#2{\iftX@\else
 \global\tX@true
 \global\tXdimen@i\hunit\global\tXdimen@i#1\tXdimen@i\relax
 \global\tXdimen@ii\vunit\global\tXdimen@ii#2\tXdimen@ii\relax
 \fi}
\def\dtY@{\ifnum\catcode`\;=\active\expandafter\dtYA@\else
 \expandafter\dtYO@\fi}
\def\dtYO@(#1;#2){\dtY@@{#1}{#2}}
\def\dtY@@#1#2{\iftY@\else
 \global\tY@true
 \global\tYdimen@i\hunit\global\tYdimen@i#1\tYdimen@i\relax
 \global\tYdimen@ii\vunit\global\tYdimen@ii#2\tYdimen@ii\relax
 \fi}
{\catcode`\;=\active
 \gdef\dsA@(#1;#2){\ds@@{#1}{#2}}
 \gdef\dtXA@(#1;#2){\dtX@@{#1}{#2}}
 \gdef\dtYA@(#1;#2){\dtY@@{#1}{#2}}
}
\def\da@#1{\ifa@\else\global\a@true\global\exacount@#1\relax\fi}
\def\dx@#1{\ifx@\else
 \global\x@true
 \global\xdimen@\hunit\global\xdimen@#1\xdimen@\relax
 \fi}
\def\dX@#1{\ifX@\else
 \global\X@true
 \global\Xdimen@\hunit\global\Xdimen@#1\Xdimen@\relax
 \fi}
\def\dy@#1{\ify@\else
 \global\y@true
 \global\ydimen@\vunit\global\ydimen@#1\ydimen@\relax
 \fi}
\def\dY@#1{\ifY@\else
 \global\Y@true
 \global\Ydimen@\vunit\global\Ydimen@#1\Ydimen@\relax
 \fi}
\def\p@@#1{\ifp@\else
 \global\p@true
 \global\pdimen@\hunit\global\divide\pdimen@\tw@
 \global\pdimen@#1\pdimen@\relax
 \fi}
\def\L@#1{\iflabel@i\else
 \global\label@itrue\gdef\label@i{#1}%
 \fi}
\def\l@#1{\iflabel@ii\else
 \global\label@iitrue\gdef\label@ii{#1}%
 \fi}
\def\dL@#1{\ifdl@i\else
 \global\dl@itrue\global\ldimen@i\hunit\global\ldimen@i#1\ldimen@i\relax
 \fi}
\def\dl@#1{\ifdl@ii\else
 \global\dl@iitrue\global\ldimen@ii\hunit\global\ldimen@ii#1\ldimen@ii\relax
 \fi}
\def\s@{\ifunshort@\else\global\short@true\fi}
\def\uns@{\ifshort@\else\global\unshort@true\global\short@false\fi}
\def\optioncodes@{\let\0\zero@\let\1\one@\let\a\a@\let\ds\ds@\let\dtX\dtX@
 \let\dtY\dtY@\let\da\da@\let\dx\dx@\let\dX\dX@\let\dY\dY@\let\dy\dy@
 \let\p\p@@\let\L\L@\let\l\l@\let\dL\dL@\let\dl\dl@\let\s\s@\let\uns\uns@}
\def\slopes@{\\161\\152\\143\\134\\255\\126\\357\\238\\349\\45{10}\\56{11}%
 \\11{12}\\65{13}\\54{14}\\43{15}\\32{16}\\53{17}\\21{18}\\52{19}\\31{20}%
 \\41{21}\\51{22}\\61{23}}
\newcount\tan@i
\newcount\tan@ip
\newcount\tan@ii
\newcount\tan@iip
\newdimen\slope@i
\newdimen\slope@ip
\newdimen\slope@ii
\newdimen\slope@iip
\newcount\angcount@
\newcount\extracount@
\def\slope@{{\slope@i\secondy@\advance\slope@i-\firsty@
 \ifN@\else\multiply\slope@i\m@ne\fi
 \slope@ii\secondx@\advance\slope@ii-\firstx@
 \ifE@\else\multiply\slope@ii\m@ne\fi
 \ifdim\slope@ii<\z@
  \global\tan@i6 \global\tan@ii\@ne\global\angcount@23
 \else
  \dimen@\slope@i\multiply\dimen@6
  \ifdim\dimen@<\slope@ii
   \global\tan@i\@ne\global\tan@ii6 \global\angcount@\@ne
  \else
   \dimen@\slope@ii\multiply\dimen@6
   \ifdim\dimen@<\slope@i
    \global\tan@i6 \global\tan@ii\@ne\global\angcount@23
   \else
    \global\tan@ip\z@\global\tan@iip\@ne
    \def\\##1##2##3{\global\angcount@##3\relax
     \slope@ip\slope@i\slope@iip\slope@ii
     \multiply\slope@iip##1\relax\multiply\slope@ip##2\relax
     \ifdim\slope@iip<\slope@ip
      \global\tan@ip##1\relax\global\tan@iip##2\relax
     \else
      \global\tan@i##1\relax\global\tan@ii##2\relax
      \def\\####1####2####3{}%
     \fi}%
    \slopes@
    \slope@i\secondy@\advance\slope@i-\firsty@
    \ifN@\else\multiply\slope@i\m@ne\fi
    \multiply\slope@i\tan@ii\multiply\slope@i\tan@iip\multiply\slope@i\tw@
    \count@\tan@i\multiply\count@\tan@iip
    \extracount@\tan@ip\multiply\extracount@\tan@ii
    \advance\count@\extracount@
    \slope@ii\secondx@\advance\slope@ii-\firstx@
    \ifE@\else\multiply\slope@ii\m@ne\fi
    \multiply\slope@ii\count@
    \ifdim\slope@i<\slope@ii
     \global\tan@i\tan@ip\global\tan@ii\tan@iip
     \global\advance\angcount@\m@ne
    \fi
   \fi
  \fi
 \fi}%
}
\def\slope@a#1{{\def\\##1##2##3{\ifnum##3=#1\global\tan@i##1\relax
 \global\tan@ii##2\relax\fi}\slopes@}}
\newcount\i@
\newcount\j@
\newcount\colcount@
\newcount\Colcount@
\newcount\tcolcount@
\newdimen\rowht@
\newdimen\rowdp@
\newcount\rowcount@
\newcount\Rowcount@
\newcount\maxcolrow@
\newtoks\colwidthtoks@
\newtoks\Rowheighttoks@
\newtoks\Rowdepthtoks@
\newtoks\widthtoks@
\newtoks\Widthtoks@
\newtoks\heighttoks@
\newtoks\Heighttoks@
\newtoks\depthtoks@
\newtoks\Depthtoks@
\newif\iffirstCDcr@
\def\dotoks@i{%
 \global\widthtoks@\expandafter{\the\widthtoks@\else\getdim@\z@\fi}%
 \global\heighttoks@\expandafter{\the\heighttoks@\else\getdim@\z@\fi}%
 \global\depthtoks@\expandafter{\the\depthtoks@\else\getdim@\z@\fi}}
\def\dotoks@ii{%
 \global\widthtoks@{\ifcase\j@}%
 \global\heighttoks@{\ifcase\j@}%
 \global\depthtoks@{\ifcase\j@}}
\def\preCD@#1\endCD{\setbox\ZER@
 \vbox{%
  \def\arrow@##1##2{{}}%
  \global\rowcount@\m@ne\global\colcount@\z@\global\Colcount@\z@
  \global\firstCDcr@true\toks@{}%
  \global\widthtoks@{\ifcase\j@}%
  \global\Widthtoks@{\ifcase\i@}%
  \global\heighttoks@{\ifcase\j@}%
  \global\Heighttoks@{\ifcase\i@}%
  \global\depthtoks@{\ifcase\j@}%
  \global\Depthtoks@{\ifcase\i@}%
  \global\Rowheighttoks@{\ifcase\i@}%
  \global\Rowdepthtoks@{\ifcase\i@}%
  \Let@
  \everycr{%
   \noalign{%
    \global\advance\rowcount@\@ne
    \ifnum\colcount@<\Colcount@
    \else
     \global\Colcount@\colcount@\global\maxcolrow@\rowcount@
    \fi
    \global\colcount@\z@
    \iffirstCDcr@
     \global\firstCDcr@false
    \else
     \edef\next@{\the\Rowheighttoks@\noexpand\or\noexpand\getdim@\the\rowht@}%
      \global\Rowheighttoks@\expandafter{\next@}%
     \edef\next@{\the\Rowdepthtoks@\noexpand\or\noexpand\getdim@\the\rowdp@}%
      \global\Rowdepthtoks@\expandafter{\next@}%
     \global\rowht@\z@\global\rowdp@\z@
     \dotoks@i
     \edef\next@{\the\Widthtoks@\noexpand\or\the\widthtoks@}%
      \global\Widthtoks@\expandafter{\next@}%
     \edef\next@{\the\Heighttoks@\noexpand\or\the\heighttoks@}%
      \global\Heighttoks@\expandafter{\next@}%
     \edef\next@{\the\Depthtoks@\noexpand\or\the\depthtoks@}%
      \global\Depthtoks@\expandafter{\next@}%
     \dotoks@ii
    \fi}}%
  \tabskip\z@
  \halign{&\setbox\ZER@\hbox{\vrule\height\ten@\p@\width\z@\depth\z@     
   $\m@th\displaystyle{##}$}\copy\ZER@
   \ifdim\ht\ZER@>\rowht@\global\rowht@\ht\ZER@\fi
   \ifdim\dp\ZER@>\rowdp@\global\rowdp@\dp\ZER@\fi
   \global\advance\colcount@\@ne
   \edef\next@{\the\widthtoks@\noexpand\or\noexpand\getdim@\the\wd\ZER@}%
    \global\widthtoks@\expandafter{\next@}%
   \edef\next@{\the\heighttoks@\noexpand\or\noexpand\getdim@\the\ht\ZER@}%
    \global\heighttoks@\expandafter{\next@}%
   \edef\next@{\the\depthtoks@\noexpand\or\noexpand\getdim@\the\dp\ZER@}%
    \global\depthtoks@\expandafter{\next@}%
   \cr#1\crcr}}%
 \Rowcount@\rowcount@
 \global\Widthtoks@\expandafter{\the\Widthtoks@\fi\relax}%
 \edef\Width@##1##2{\i@##1\relax\j@##2\relax\the\Widthtoks@}%
 \global\Heighttoks@\expandafter{\the\Heighttoks@\fi\relax}%
 \edef\Height@##1##2{\i@##1\relax\j@##2\relax\the\Heighttoks@}%
 \global\Depthtoks@\expandafter{\the\Depthtoks@\fi\relax}%
 \edef\Depth@##1##2{\i@##1\relax\j@##2\relax\the\Depthtoks@}%
 \edef\next@{\the\Rowheighttoks@\noexpand\fi\relax}%
 \global\Rowheighttoks@\expandafter{\next@}%
 \edef\Rowheight@##1{\i@##1\relax\the\Rowheighttoks@}%
 \edef\next@{\the\Rowdepthtoks@\noexpand\fi\relax}%
 \global\Rowdepthtoks@\expandafter{\next@}%
 \edef\Rowdepth@##1{\i@##1\relax\the\Rowdepthtoks@}%
 \global\colwidthtoks@{\fi}%
 \setbox\ZER@\vbox{%
  \unvbox\ZER@
  \count@\rowcount@
  \loop
   \unskip\unpenalty
   \setbox\ZER@\lastbox
   \ifnum\count@>\maxcolrow@\advance\count@\m@ne
   \repeat
  \hbox{%
   \unhbox\ZER@
   \count@\z@
   \loop
    \unskip
    \setbox\ZER@\lastbox
    \edef\next@{\noexpand\or\noexpand\getdim@\the\wd\ZER@\the\colwidthtoks@}%
     \global\colwidthtoks@\expandafter{\next@}%
    \advance\count@\@ne
    \ifnum\count@<\Colcount@
    \repeat}}%
 \edef\next@{\noexpand\ifcase\noexpand\i@\the\colwidthtoks@}%
  \global\colwidthtoks@\expandafter{\next@}%
 \edef\Colwidth@##1{\i@##1\relax\the\colwidthtoks@}%
 \global\colwidthtoks@{}\global\Rowheighttoks@{}\global\Rowdepthtoks@{}%
 \global\widthtoks@{}\global\Widthtoks@{}\global\heighttoks@{}%
 \global\Heighttoks@{}\global\depthtoks@{}\global\Depthtoks@{}%
}
\newcount\xoff@
\newcount\yoff@
\newcount\endcount@
\newcount\rcount@
\newdimen\firstx@
\newdimen\firsty@
\newdimen\secondx@
\newdimen\secondy@
\newdimen\tocenter@
\newdimen\charht@
\newdimen\charwd@
\def\outside@{\Err@{This arrow points outside the \string\CD}}
\newif\ifsvertex@
\newif\iftvertex@
\def\arrow@#1#2{\global\xoff@#1\relax\global\yoff@#2\relax
 \count@\rowcount@\advance\count@-\yoff@
 \ifnum\count@<\@ne\outside@\else\ifnum\count@>\Rowcount@\outside@\fi\fi
 \count@\colcount@\advance\count@\xoff@
 \ifnum\count@<\@ne\outside@\else\ifnum\count@>\Colcount@\outside@\fi\fi
 \tcolcount@\colcount@\advance\tcolcount@\xoff@
 \Width@\rowcount@\colcount@\divide\getdim@\tw@\tocenter@-\getdim@
 \ifdim\getdim@=\z@
  \firstx@\z@\firsty@\mathaxis@\svertex@true
 \else
  \svertex@false
  \ifHshort@
   \Colwidth@\colcount@\divide\getdim@\tw@
   \ifE@ \firstx@\getdim@ \else \firstx@-\getdim@ \fi
  \else
   \ifE@ \firstx@\getdim@ \else \firstx@-\getdim@ \fi
  \fi
  \ifE@
   \ifH@ \advance\firstx@\thr@@\p@ \else \advance\firstx@-\thr@@\p@ \fi  
  \else
   \ifH@ \advance\firstx@-\thr@@\p@ \else \advance\firstx@\thr@@\p@ \fi  
  \fi
  \ifN@
   \Height@\rowcount@\colcount@ \firsty@\getdim@                         
   \ifV@ \advance\firsty@\thr@@\p@ \fi                                   
  \else
   \ifV@
    \Depth@\rowcount@\colcount@ \firsty@-\getdim@                        
    \advance\firsty@-\thr@@\p@                                           
   \else
    \firsty@\z@                                                          
   \fi
  \fi
 \fi
 \ifV@
 \else
  \Colwidth@\colcount@\divide\getdim@\tw@
  \ifE@\secondx@\getdim@\else\secondx@-\getdim@\fi
  \ifE@\else\getcgap@\colcount@\advance\secondx@-\getdim@\fi
  \endcount@\colcount@\advance\endcount@\xoff@
  \count@\colcount@
  \ifE@
   \advance\count@\@ne
   \loop
    \ifnum\count@<\endcount@
    \Colwidth@\count@\advance\secondx@\getdim@
    \getcgap@\count@\advance\secondx@\getdim@
    \advance\count@\@ne
    \repeat
  \else
   \advance\count@\m@ne
   \loop
    \ifnum\count@>\endcount@
    \Colwidth@\count@\advance\secondx@-\getdim@
    \getcgap@\count@\advance\secondx@-\getdim@
    \advance\count@\m@ne
    \repeat
  \fi
  \Colwidth@\count@\divide\getdim@\tw@
  \ifHshort@
  \else
   \ifE@\advance\secondx@\getdim@\else\advance\secondx@-\getdim@\fi
  \fi
  \ifE@\getcgap@\count@\advance\secondx@\getdim@\fi
  \rcount@\rowcount@\advance\rcount@-\yoff@
  \Width@\rcount@\count@\divide\getdim@\tw@
  \tvertex@false
  \ifH@\ifdim\getdim@=\z@\tvertex@true\Hshort@false\fi\fi
  \ifHshort@
  \else
   \ifE@\advance\secondx@-\getdim@\else\advance\secondx@\getdim@\fi
  \fi
  \iftvertex@
   \advance\secondx@.4\p@
  \else
   \ifE@\advance\secondx@-\thr@@\p@\else\advance\secondx@\thr@@\p@\fi    
  \fi
 \fi
 \ifH@
 \else
  \ifN@
   \Rowheight@\rowcount@\secondy@\getdim@
  \else
   \Rowdepth@\rowcount@\secondy@-\getdim@
   \getrgap@\rowcount@\advance\secondy@-\getdim@
  \fi
  \endcount@\rowcount@\advance\endcount@-\yoff@
  \count@\rowcount@
  \ifN@
   \advance\count@\m@ne
   \loop
    \ifnum\count@>\endcount@
    \Rowheight@\count@\advance\secondy@\getdim@
    \Rowdepth@\count@\advance\secondy@\getdim@
    \getrgap@\count@\advance\secondy@\getdim@
    \advance\count@\m@ne
    \repeat
  \else
   \advance\count@\@ne
   \loop
    \ifnum\count@<\endcount@
    \Rowheight@\count@\advance\secondy@-\getdim@
    \Rowdepth@\count@\advance\secondy@-\getdim@
    \getrgap@\count@\advance\secondy@-\getdim@
    \advance\count@\@ne
    \repeat
  \fi
  \tvertex@false
  \ifV@\Width@\count@\colcount@\ifdim\getdim@=\z@\tvertex@true\fi\fi
  \ifN@
   \getrgap@\count@\advance\secondy@\getdim@
   \Rowdepth@\count@\advance\secondy@\getdim@
   \iftvertex@
    \advance\secondy@\mathaxis@
   \else
    \Depth@\count@\tcolcount@\advance\secondy@-\getdim@
    \advance\secondy@-\thr@@\p@                                          
   \fi
  \else
   \Rowheight@\count@\advance\secondy@-\getdim@
   \iftvertex@
    \advance\secondy@\mathaxis@
   \else
    \Height@\count@\tcolcount@\advance\secondy@\getdim@
    \advance\secondy@\thr@@\p@                                           
   \fi
  \fi
 \fi
 \ifV@\else\advance\firstx@\sxdimen@\fi
 \ifH@\else\advance\firsty@\sydimen@\fi
 \iftX@
  \advance\secondy@\tXdimen@ii
  \advance\secondx@\tXdimen@i
  \slope@
 \else
  \iftY@
   \advance\secondy@\tYdimen@ii
   \advance\secondx@\tYdimen@i
   \slope@
   \secondy@\secondx@\advance\secondy@-\firstx@
   \ifNESW@\else\multiply\secondy@\m@ne\fi
   \multiply\secondy@\tan@i\divide\secondy@\tan@ii\advance\secondy@\firsty@
  \else
   \ifa@
    \slope@
    \ifNESW@\global\advance\angcount@\exacount@\else
     \global\advance\angcount@-\exacount@\fi
    \ifnum\angcount@>23 \global\angcount@23 \fi
    \ifnum\angcount@<\@ne\global\angcount@\@ne\fi
    \slope@a\angcount@
    \ifY@
     \advance\secondy@\Ydimen@
    \else
     \ifX@
      \advance\secondx@\Xdimen@
      \dimen@\secondx@\advance\dimen@-\firstx@
      \ifNESW@\else\multiply\dimen@\m@ne\fi
      \multiply\dimen@\tan@i\divide\dimen@\tan@ii
      \advance\dimen@\firsty@\secondy@\dimen@
     \fi
    \fi
   \else
    \ifH@\else\ifV@\else\slope@\fi\fi
   \fi
  \fi
 \fi
 \ifH@\else\ifV@\else\ifsvertex@\else
  \dimen@6\p@\multiply\dimen@\tan@ii
  \count@\tan@i\advance\count@\tan@ii\divide\dimen@\count@
  \ifE@\advance\firstx@\dimen@\else\advance\firstx@-\dimen@\fi
  \multiply\dimen@\tan@i\divide\dimen@\tan@ii
  \ifN@\advance\firsty@\dimen@\else\advance\firsty@-\dimen@\fi
 \fi\fi\fi
 \ifp@
  \ifH@\else\ifV@\else
   \getcos@\pdimen@\advance\firsty@\dimen@\advance\secondy@\dimen@
   \ifNESW@\advance\firstx@-\dimen@ii\else\advance\firstx@\dimen@ii\fi
  \fi\fi
 \fi
 \ifH@\else\ifV@\else
  \ifnum\tan@i>\tan@ii
   \charht@\ten@\p@\charwd@\ten@\p@
   \multiply\charwd@\tan@ii\divide\charwd@\tan@i
  \else
   \charwd@\ten@\p@\charht@\ten@\p@
   \divide\charht@\tan@ii\multiply\charht@\tan@i
  \fi
  \ifnum\tcount@=\thr@@
   \ifN@\advance\secondy@-.3\charht@\else\advance\secondy@.3\charht@\fi
  \fi
  \ifnum\scount@=\tw@
   \ifE@\advance\firstx@.3\charht@\else\advance\firstx@-.3\charht@\fi
  \fi
  \ifnum\tcount@=12
   \ifN@\advance\secondy@-\charht@\else\advance\secondy@\charht@\fi
  \fi
  \iftY@
  \else
   \ifa@
    \ifX@
    \else
     \secondx@\secondy@\advance\secondx@-\firsty@
     \ifNESW@\else\multiply\secondx@\m@ne\fi
     \multiply\secondx@\tan@ii\divide\secondx@\tan@i
     \advance\secondx@\firstx@
    \fi
   \fi
  \fi
 \fi\fi
 \ifH@\harrow@\else\ifV@\varrow@\else\arrow@@\fi\fi}
\newdimen\mathaxis@
\mathaxis@90\p@\divide\mathaxis@36
\def\harrow@b{\ifE@\hskip\tocenter@\hskip\firstx@\fi}
\def\harrow@bb{\ifE@\hskip\xdimen@\else\hskip\Xdimen@\fi}
\def\harrow@e{\ifE@\else\hskip-\firstx@\hskip-\tocenter@\fi}
\def\harrow@ee{\ifE@\hskip-\Xdimen@\else\hskip-\xdimen@\fi}
\def\harrow@{\dimen@\secondx@\advance\dimen@-\firstx@
 \ifE@\let\next@\rlap\else\multiply\dimen@\m@ne\let\next@\llap\fi
 \next@{%
  \harrow@b
  \smash{\raise\pdimen@\hbox to\dimen@
   {\harrow@bb\arrow@ii
    \ifnum\arrcount@=\m@ne\else\ifnum\arrcount@=\thr@@\else
     \ifE@
      \ifnum\scount@=\m@ne
      \else
       \ifcase\scount@\or\or\char118 \or\char117 \or\or\or\char119 \or
       \char120 \or\char121 \or\char122 \or\or\or\arrow@i\char125 \or
       \char117 \hskip\thr@@\p@\char117 \hskip-\thr@@\p@\fi
      \fi
     \else
      \ifnum\tcount@=\m@ne
      \else
       \ifcase\tcount@\char117 \or\or\char117 \or\char118 \or\char119 \or
       \char120 \or\or\or\or\or\char121 \or\char122 \or\arrow@i\char125
       \or\char117 \hskip\thr@@\p@\char117 \hskip-\thr@@\p@\fi
      \fi
     \fi
    \fi\fi
    \dimen@\mathaxis@\advance\dimen@.2\p@
    \dimen@ii\mathaxis@\advance\dimen@ii-.2\p@
    \ifnum\arrcount@=\m@ne
     \let\leads@\null
    \else
     \ifcase\arrcount@
      \def\leads@{\hrule\height\dimen@\depth-\dimen@ii}\or
      \def\leads@{\hrule\height\dimen@\depth-\dimen@ii}\or
      \def\leads@{\hbox to\ten@\p@{%
       \leaders\hrule\height\dimen@\depth-\dimen@ii\hfil
       \hfil
      \leaders\hrule\height\dimen@\depth-\dimen@ii\hskip\z@ plus2fil\relax
       \hfil
       \leaders\hrule\height\dimen@\depth-\dimen@ii\hfil}}\or
     \def\leads@{\hbox{\hbox to\ten@\p@{\dimen@\mathaxis@\advance\dimen@1.2\p@
       \dimen@ii\dimen@\advance\dimen@ii-.4\p@
       \leaders\hrule\height\dimen@\depth-\dimen@ii\hfil}%
       \kern-\ten@\p@
       \hbox to\ten@\p@{\dimen@\mathaxis@\advance\dimen@-1.2\p@
       \dimen@ii\dimen@\advance\dimen@ii-.4\p@
       \leaders\hrule\height\dimen@\depth-\dimen@ii\hfil}}}\fi
    \fi
    \cleaders\leads@\hfil
    \ifnum\arrcount@=\m@ne\else\ifnum\arrcount@=\thr@@\else
     \arrow@i
     \ifE@
      \ifnum\tcount@=\m@ne
      \else
       \ifcase\tcount@\char119 \or\or\char119 \or\char120 \or\char121 \or
       \char122 \or \or\or\or\or\char123 \or\char124 \or
       \char125 \or\char119 \hskip-\thr@@\p@\char119 \hskip\thr@@\p@\fi
      \fi
     \else
      \ifcase\scount@\or\or\char120 \or\char119 \or\or\or\char121 \or\char122
      \or\char123 \or\char124 \or\or\or\char125 \or
      \char119 \hskip-\thr@@\p@\char119 \hskip\thr@@\p@\fi
     \fi
    \fi\fi
    \harrow@ee}}%
  \harrow@e}%
 \iflabel@i
  \dimen@ii\z@\setbox\ZER@\hbox{$\m@th\tsize@@\label@i$}%
  \ifnum\arrcount@=\m@ne
  \else
   \advance\dimen@ii\mathaxis@
   \advance\dimen@ii\dp\ZER@\advance\dimen@ii\tw@\p@
   \ifnum\arrcount@=\thr@@\advance\dimen@ii\tw@\p@\fi
  \fi
  \advance\dimen@ii\pdimen@
  \next@{\harrow@b\smash{\raise\dimen@ii\hbox to\dimen@
   {\harrow@bb\hskip\tw@\ldimen@i\hfil\box\ZER@\hfil\harrow@ee}}\harrow@e}%
 \fi
 \iflabel@ii
  \ifnum\arrcount@=\m@ne
  \else
   \setbox\ZER@\hbox{$\m@th\tsize@\label@ii$}%
   \dimen@ii-\ht\ZER@\advance\dimen@ii-\tw@\p@
   \ifnum\arrcount@=\thr@@\advance\dimen@ii-\tw@\p@\fi
   \advance\dimen@ii\mathaxis@\advance\dimen@ii\pdimen@
   \next@{\harrow@b\smash{\raise\dimen@ii\hbox to\dimen@
    {\harrow@bb\hskip\tw@\ldimen@ii\hfil\box\ZER@\hfil\harrow@ee}}\harrow@e}%
  \fi
 \fi}
\let\tsize@\tsize
\def\tsizeCDlabels{\let\tsize@\tsize}
\def\ssizeCDlabels{\let\tsize@\ssize}
\def\tsize@@{\ifnum\arrcount@=\m@ne\else\tsize@\fi}
\def\varrow@{\dimen@\secondy@\advance\dimen@-\firsty@
 \ifN@\else\multiply\dimen@\m@ne\fi
 \setbox\ZER@\vbox to\dimen@
  {\ifN@\vskip-\Ydimen@\else\vskip\ydimen@\fi
   \ifnum\arrcount@=\m@ne\else\ifnum\arrcount@=\thr@@\else
    \hbox{\arrow@iii
     \ifN@
      \ifnum\tcount@=\m@ne
      \else
       \ifcase\tcount@\char117 \or\or\char117 \or\char118 \or\char119 \or
       \char120 \or\or\or\or\or\char121 \or\char122 \or\char123 \or
       \vbox{\hbox{\char117}\nointerlineskip\vskip\thr@@\p@
       \hbox{\char117}\vskip-\thr@@\p@}\fi
      \fi
     \else
      \ifcase\scount@\or\or\char118 \or\char117 \or\or\or\char119 \or
      \char120 \or\char121 \or\char122 \or\or\or\char123 \or
      \vbox{\hbox{\char117}\nointerlineskip\vskip\thr@@\p@
      \hbox{\char117}\vskip-\thr@@\p@}\fi
     \fi}%
    \nointerlineskip
   \fi\fi
   \ifnum\arrcount@=\m@ne
    \let\leads@\null
   \else
    \ifcase\arrcount@\let\leads@\vrule\or\let\leads@\vrule\or
    \def\leads@{\vbox to\ten@\p@{%
     \hrule\height1.67\p@\depth\z@\width.4\p@
     \vfil
     \hrule\height3.33\p@\depth\z@\width.4\p@
     \vfil
     \hrule\height1.67\p@\depth\z@\width.4\p@}}\or
    \def\leads@{\hbox{\vrule\height\p@\hskip\tw@\p@\vrule}}\fi
   \fi
  \cleaders\leads@\vfill\nointerlineskip
   \ifnum\arrcount@=\m@ne\else\ifnum\arrcount@=\thr@@\else
    \hbox{\arrow@iv
     \ifN@
      \ifcase\scount@\or\or\char118 \or\char117 \or\or\or\char119 \or
      \char120 \or\char121 \or\char122 \or\or\or\arrow@iii\char123 \or
      \vbox{\hbox{\char117}\nointerlineskip\vskip-\thr@@\p@
      \hbox{\char117}\vskip\thr@@\p@}\fi
     \else
      \ifnum\tcount@=\m@ne
      \else
       \ifcase\tcount@\char117 \or\or\char117 \or\char118 \or\char119 \or
       \char120 \or\or\or\or\or\char121 \or\char122 \or\arrow@iii\char123 \or
       \vbox{\hbox{\char117}\nointerlineskip\vskip-\thr@@\p@
       \hbox{\char117}\vskip\thr@@\p@}\fi
      \fi
     \fi}%
   \fi\fi
   \ifN@\vskip\ydimen@\else\vskip-\Ydimen@\fi}%
 \ifN@
  \dimen@ii\firsty@
 \else
  \dimen@ii-\firsty@\advance\dimen@ii\ht\ZER@\multiply\dimen@ii\m@ne
 \fi
 \rlap{\smash{\hskip\tocenter@\hskip\pdimen@\raise\dimen@ii\box\ZER@}}%
 \iflabel@i
  \setbox\ZER@\vbox to\dimen@{\vfil
   \hbox{$\m@th\tsize@@\label@i$}\vskip\tw@\ldimen@i\vfil}%
  \rlap{\smash{\hskip\tocenter@\hskip\pdimen@
  \ifnum\arrcount@=\m@ne\let\next@\relax\else\let\next@\llap\fi
  \next@{\raise\dimen@ii\hbox{\ifnum\arrcount@=\m@ne\hskip-.5\wd\ZER@\fi
   \box\ZER@\ifnum\arrcount@=\m@ne\else\hskip\tw@\p@\fi}}}}%
 \fi
 \iflabel@ii
  \ifnum\arrcount@=\m@ne
  \else
   \setbox\ZER@\vbox to\dimen@{\vfil
    \hbox{$\m@th\tsize@\label@ii$}\vskip\tw@\ldimen@ii\vfil}%
   \rlap{\smash{\hskip\tocenter@\hskip\pdimen@
   \rlap{\raise\dimen@ii\hbox{\ifnum\arrcount@=\thr@@\hskip4.5\p@\else
    \hskip2.5\p@\fi\box\ZER@}}}}%
  \fi
 \fi
}
\newdimen\goal@
\newdimen\shifted@
\newcount\Tcount@
\newcount\Scount@
\newbox\shaft@
\newcount\slcount@
\def\getcos@#1{%
 \ifnum\tan@i<\tan@ii
  \dimen@#1%
  \ifnum\slcount@<8 \count@9 \else \ifnum\slcount@<12 \count@8 \else
   \count@7 \fi\fi
  \multiply\dimen@\count@\divide\dimen@\ten@
  \dimen@ii\dimen@\multiply\dimen@ii\tan@i\divide\dimen@ii\tan@ii
 \else
  \dimen@ii#1%
  \count@-\slcount@\advance\count@24
  \ifnum\count@<8 \count@9 \else \ifnum\count@<12 \count@8
   \else\count@7 \fi\fi
  \multiply\dimen@ii\count@\divide\dimen@ii\ten@
  \dimen@\dimen@ii\multiply\dimen@\tan@ii\divide\dimen@\tan@i
 \fi}
\newdimen\adjust@
\def\Nnext@{\ifN@\let\next@\raise\else\let\next@\lower\fi}
\def\arrow@@{\slcount@\angcount@
 \ifNESW@
  \ifnum\angcount@<\ten@
   \let\arrowfont@\arrow@i\global\advance\angcount@\m@ne
   \global\multiply\angcount@13
  \else
   \ifnum\angcount@<19
    \let\arrowfont@\arrow@ii\global\advance\angcount@-\ten@
    \global\multiply\angcount@13
   \else
    \let\arrowfont@\arrow@iii\global\advance\angcount@-19
    \global\multiply\angcount@13
  \fi\fi
  \Tcount@\angcount@
 \else
  \ifnum\angcount@<5
   \let\arrowfont@\arrow@iii\global\advance\angcount@\m@ne
   \global\multiply\angcount@13 \global\advance\angcount@65
  \else
   \ifnum\angcount@<14
    \let\arrowfont@\arrow@iv\global\advance\angcount@-5
    \global\multiply\angcount@13
   \else
    \ifnum\angcount@<23
     \let\arrowfont@\arrow@v\global\advance\angcount@-14
     \global\multiply\angcount@13
    \else
     \let\arrowfont@\arrow@i\global\angcount@117
  \fi\fi\fi
  \ifnum\angcount@=117 \Tcount@115 \else\Tcount@\angcount@\fi
 \fi
 \Scount@\Tcount@
 \ifE@
  \ifnum\tcount@=\z@\advance\Tcount@\tw@\else\ifnum\tcount@=13
   \advance\Tcount@\tw@\else\advance\Tcount@\tcount@\fi\fi
  \ifnum\scount@=\z@\else\ifnum\scount@=13 \advance\Scount@\thr@@\else
   \advance\Scount@\scount@\fi\fi
 \else
  \ifcase\tcount@\advance\Tcount@\thr@@\or\or\advance\Tcount@\thr@@\or
  \advance\Tcount@\tw@\or\advance\Tcount@6 \or\advance\Tcount@7
  \or\or\or\or\or\advance\Tcount@8 \or\advance\Tcount@9 \or
  \advance\Tcount@12 \or\advance\Tcount@\thr@@\fi
  \ifcase\scount@\or\or\advance\Scount@\thr@@\or\advance\Scount@\tw@\or
  \or\or\advance\Scount@4 \or\advance\Scount@5 \or\advance\Scount@\ten@
  \or\advance\Scount@11 \or\or\or\advance\Scount@12 \or\advance
  \Scount@\tw@\fi
 \fi
 \ifcase\arrcount@\or\or\global\advance\angcount@\@ne\else\fi
 \ifN@\shifted@\firsty@\else\shifted@-\firsty@\fi
 \ifE@\else\advance\shifted@\charht@\fi
 \goal@\secondy@\advance\goal@-\firsty@
 \ifN@\else\multiply\goal@\m@ne\fi
 \setbox\shaft@\hbox{\arrowfont@\char\angcount@}%
 \ifnum\arrcount@=\thr@@
  \getcos@{1.5\p@}%
  \setbox\shaft@\hbox to\wd\shaft@{\arrowfont@
   \rlap{\hskip\dimen@ii
    \smash{\ifNESW@\let\next@\lower\else\let\next@\raise\fi
     \next@\dimen@\hbox{\arrowfont@\char\angcount@}}}%
   \rlap{\hskip-\dimen@ii
    \smash{\ifNESW@\let\next@\raise\else\let\next@\lower\fi
      \next@\dimen@\hbox{\arrowfont@\char\angcount@}}}\hfil}%
 \fi
 \rlap{\smash{\hskip\tocenter@\hskip\firstx@
  \ifnum\arrcount@=\m@ne
  \else
   \ifnum\arrcount@=\thr@@
   \else
    \ifnum\scount@=\m@ne
    \else
     \ifnum\scount@=\z@
     \else
      \setbox\ZER@\hbox{\ifnum\angcount@=117 \arrow@v\else\arrowfont@\fi
       \char\Scount@}%
      \ifNESW@
       \ifnum\scount@=\tw@
        \dimen@\shifted@\advance\dimen@-\charht@
        \ifN@\hskip-\wd\ZER@\fi
        \Nnext@
        \next@\dimen@\copy\ZER@
        \ifN@\else\hskip-\wd\ZER@\fi
       \else
        \Nnext@
        \ifN@\else\hskip-\wd\ZER@\fi
        \next@\shifted@\copy\ZER@
        \ifN@\hskip-\wd\ZER@\fi
       \fi
       \ifnum\scount@=12
        \advance\shifted@\charht@\advance\goal@-\charht@
        \ifN@\hskip\wd\ZER@\else\hskip-\wd\ZER@\fi
       \fi
       \ifnum\scount@=13
        \getcos@{\thr@@\p@}%
        \ifN@\hskip\dimen@\else\hskip-\wd\ZER@\hskip-\dimen@\fi
        \adjust@\shifted@\advance\adjust@\dimen@ii
        \Nnext@
        \next@\adjust@\copy\ZER@
        \ifN@\hskip-\dimen@\hskip-\wd\ZER@\else\hskip\dimen@\fi
       \fi
      \else
       \ifN@\hskip-\wd\ZER@\fi
       \ifnum\scount@=\tw@
        \ifN@\hskip\wd\ZER@\else\hskip-\wd\ZER@\fi
        \dimen@\shifted@\advance\dimen@-\charht@
        \Nnext@
        \next@\dimen@\copy\ZER@
        \ifN@\hskip-\wd\ZER@\fi
       \else
        \Nnext@
        \next@\shifted@\copy\ZER@
        \ifN@\else\hskip-\wd\ZER@\fi
       \fi
       \ifnum\scount@=12
        \advance\shifted@\charht@\advance\goal@-\charht@
        \ifN@\hskip-\wd\ZER@\else\hskip\wd\ZER@\fi
       \fi
       \ifnum\scount@=13
        \getcos@{\thr@@\p@}%
        \ifN@\hskip-\wd\ZER@\hskip-\dimen@\else\hskip\dimen@\fi
        \adjust@\shifted@\advance\adjust@\dimen@ii
        \Nnext@
        \next@\adjust@\copy\ZER@
        \ifN@\hskip\dimen@\else\hskip-\dimen@\hskip-\wd\ZER@\fi
       \fi	
      \fi
  \fi\fi\fi\fi
  \ifnum\arrcount@=\m@ne
  \else
   \loop
    \ifdim\goal@>\charht@
    \ifE@\else\hskip-\charwd@\fi
    \Nnext@
    \next@\shifted@\copy\shaft@
    \ifE@\else\hskip-\charwd@\fi
    \advance\shifted@\charht@\advance\goal@-\charht@
    \repeat
   \ifdim\goal@>\z@
    \dimen@\charht@\advance\dimen@-\goal@
    \divide\dimen@\tan@i\multiply\dimen@\tan@ii
    \ifE@\hskip-\dimen@\else\hskip-\charwd@\hskip\dimen@\fi
    \adjust@\shifted@\advance\adjust@-\charht@\advance\adjust@\goal@
    \Nnext@
    \next@\adjust@\copy\shaft@
    \ifE@\else\hskip-\charwd@\fi
   \else
    \adjust@\shifted@\advance\adjust@-\charht@
   \fi
  \fi
  \ifnum\arrcount@=\m@ne
  \else
   \ifnum\arrcount@=\thr@@
   \else
    \ifnum\tcount@=\m@ne
    \else
     \setbox\ZER@
      \hbox{\ifnum\angcount@=117 \arrow@v\else\arrowfont@\fi\char\Tcount@}%
     \ifnum\tcount@=\thr@@
      \advance\adjust@\charht@
      \ifE@\else\ifN@\hskip-\charwd@\else\hskip-\wd\ZER@\fi\fi
     \else
      \ifnum\tcount@=12
       \advance\adjust@\charht@
       \ifE@\else\ifN@\hskip-\charwd@\else\hskip-\wd\ZER@\fi\fi
      \else
       \ifE@\hskip-\wd\ZER@\fi
     \fi\fi
     \Nnext@
     \next@\adjust@\copy\ZER@
     \ifnum\tcount@=13
      \hskip-\wd\ZER@
      \getcos@{\thr@@\p@}%
      \ifE@\hskip-\dimen@\else\hskip\dimen@\fi
      \advance\adjust@-\dimen@ii
      \Nnext@
      \next@\adjust@\box\ZER@
     \fi
  \fi\fi\fi}}%
 \iflabel@i
  \rlap{\hskip\tocenter@
  \dimen@\firstx@\advance\dimen@\secondx@\divide\dimen@\tw@
  \advance\dimen@\ldimen@i
  \dimen@ii\firsty@\advance\dimen@ii\secondy@\divide\dimen@ii\tw@
  \global\multiply\ldimen@i\tan@i\global\divide\ldimen@i\tan@ii
  \ifNESW@\advance\dimen@ii\ldimen@i\else\advance\dimen@ii-\ldimen@i\fi
  \setbox\ZER@\hbox{\ifNESW@\else\ifnum\arrcount@=\thr@@\hskip4\p@\else
   \hskip\tw@\p@\fi\fi
   $\m@th\tsize@@\label@i$\ifNESW@\ifnum\arrcount@=\thr@@\hskip4\p@\else
   \hskip\tw@\p@\fi\fi}%
  \ifnum\arrcount@=\m@ne
   \ifNESW@\advance\dimen@.5\wd\ZER@\advance\dimen@\p@\else
    \advance\dimen@-.5\wd\ZER@\advance\dimen@-\p@\fi
   \advance\dimen@ii-.5\ht\ZER@
  \else
   \advance\dimen@ii\dp\ZER@
   \ifnum\slcount@<6 \advance\dimen@ii\tw@\p@\fi
  \fi
  \hskip\dimen@
  \ifNESW@\let\next@\llap\else\let\next@\rlap\fi
  \next@{\smash{\raise\dimen@ii\box\ZER@}}}%
 \fi
 \iflabel@ii
  \ifnum\arrcount@=\m@ne
  \else
   \rlap{\hskip\tocenter@
   \dimen@\firstx@\advance\dimen@\secondx@\divide\dimen@\tw@
   \ifNESW@\advance\dimen@\ldimen@ii\else\advance\dimen@-\ldimen@ii\fi
   \dimen@ii\firsty@\advance\dimen@ii\secondy@\divide\dimen@ii\tw@
   \global\multiply\ldimen@ii\tan@i\global\divide\ldimen@ii\tan@ii
   \advance\dimen@ii\ldimen@ii
   \setbox\ZER@\hbox{\ifNESW@\ifnum\arrcount@=\thr@@\hskip4\p@\else
    \hskip\tw@\p@\fi\fi
    $\m@th\tsize@\label@ii$\ifNESW@\else\ifnum\arrcount@=\thr@@\hskip4\p@
    \else\hskip\tw@\p@\fi\fi}%
   \advance\dimen@ii-\ht\ZER@
   \ifnum\slcount@<9 \advance\dimen@ii-\thr@@\p@\fi
   \ifNESW@\let\next@\rlap\else\let\next@\llap\fi
   \hskip\dimen@\next@{\smash{\raise\dimen@ii\box\ZER@}}}%
  \fi
 \fi
}
\def\outCD@#1{\def#1{\Err@{\noexpand#1must not be used within \string\CD}}}
\newskip\preCDskip@
\newskip\postCDskip@
\preCDskip@\z@
\postCDskip@\z@
\def\preCDspace#1{\RIfMIfI@
 \onlydmatherr@\preCDspace\else\advance\preCDskip@#1\relax\fi\else
 \onlydmatherr@\preCDspace\fi}
\def\postCDspace#1{\RIfMIfI@
 \onlydmatherr@\postCDspace\else\advance\postCDskip@#1\relax\fi\else
 \onlydmatherr@\postCDspace\fi}
\def\predisplayspace#1{\RIfMIfI@
 \onlydmatherr@\predisplayspace\else
 \advance\abovedisplayskip#1\relax
 \advance\abovedisplayshortskip#1\relax\fi
 \else\onlydmatherr@\preCDspace\fi}
\def\postdisplayspace#1{\RIfMIfI@
 \onlydmatherr@\postdisplayspace\else
 \advance\belowdisplayskip#1\relax
 \advance\belowdisplayshortskip#1\relax\fi
 \else\onlydmatherr@\postdisplayspace\fi}
\def\PreCDSpace#1{\global\preCDskip@#1\relax}
\def\PostCDSpace#1{\global\postCDskip@#1\relax}
\def\CD#1\endCD{%
 \outCD@\cgaps\outCD@\rgaps\outCD@\Cgaps\outCD@\Rgaps
 \preCD@#1\endCD
 \advance\abovedisplayskip\preCDskip@
 \advance\abovedisplayshortskip\preCDskip@
 \advance\belowdisplayskip\postCDskip@
 \advance\belowdisplayshortskip\postCDskip@
 \vcenter{\offinterlineskip
  \vskip\preCDskip@\Let@\global\colcount@\@ne\global\rowcount@\z@
  \everycr{%
   \noalign{%
    \ifnum\rowcount@=\Rowcount@
    \else
     \getrgap@\rowcount@\vskip\getdim@
     \global\advance\rowcount@\@ne\global\colcount@\@ne
    \fi}}%
  \tabskip\z@
  \halign{&\global\xoff@\z@\global\yoff@\z@
   \getcgap@\colcount@\hskip\getdim@
   \hfil\vrule\height\ten@\p@\width\z@\depth\z@
   $\m@th\displaystyle{##}$\hfil
   \global\advance\colcount@\@ne\cr
   #1\crcr}\vskip\postCDskip@}%
 \preCDskip@\z@\postCDskip@\z@
 \def\getcgap@##1{\ifcase##1\or\getdim@\z@\else\getdim@\standardcgap\fi}%
 \def\getrgap@##1{\ifcase##1\getdim@\z@\else\getdim@\standardrgap\fi}%
 \let\Width@\relax\let\Height@\relax\let\Depth@\relax\let\Rowheight@\relax
 \let\Rowdepth@\relax\let\Colwidth@\relax
}
\let\enddocument\bye
\def\alloc@#1#2#3#4#5{\global\advance\count1#1by\@ne
  \ch@ck#1#4#2%
  \allocationnumber=\count1#1%
  \global#3#5=\allocationnumber
  \wlog{\string#5=\string#2\the\allocationnumber}}
\catcode`\@=\active